\input amstex.tex
\documentstyle{amsppt}
\UseAMSsymbols\TagsAsMath\widestnumber\key{ASSSS}
\magnification=\magstephalf\pagewidth{6.2in}\vsize8.0in
\parindent=6mm\parskip=3pt\baselineskip=16pt\tolerance=10000\hbadness=500
\NoRunningHeads\loadbold\NoBlackBoxes\nologo
\def\today{\ifcase\month\or January\or February\or March\or April\or May\or June\or
     July\or August\or September\or October\or November\or
    December\fi \space\number\day, \number\year}

\def\ns#1#2{\mskip -#1.#2mu}\def\operator#1{\ns07\operatorname{#1}\ns12}
\def\pa{\partial}\def\na{\nabla\!}
\def\dt{{\ns18{\Cal D}_{\ns23 t\!} }}
\def\tr{\operatorname{tr}}
\def\div{\operator{div}}\def\curl{\operator{curl}}
\def\dist{\operatorname{dist} }
\def\leftalignspace{\mskip -45mu}\def\rightalignspace{\mskip -90mu}

\topmatter
\title
Well posedness for the motion of a compressible liquid
with free surface boundary.
\endtitle
\author  Hans Lindblad
\endauthor
\thanks  The author was supported in part by the National
Science Foundation.  \endthanks
\address
University of California at San Diego
\endaddress
\email
lindblad\@math.ucsd.edu
\endemail

\abstract We study the motion of a compressible perfect liquid
body in vacuum. This can be thought of as a model for the motion
of the ocean or a star. The free surface moves with the velocity
of the liquid and the pressure vanishes on the free surface. This
leads to a free boundary problem for Euler's equations, where the
regularity of the boundary enters to highest order. We prove local
existence in Sobolev spaces assuming a "physical condition",
related to the fact that the pressure of a fluid has to be
positive.
\endabstract

\endtopmatter

\document

\head{1. Introduction}\endhead
We consider Euler's equations
$$
\rho \big(\partial_t +V^k\partial_{k} \big) v_j+\partial_{j\,} p=0,\quad j=1,...,n
\quad \text{in}\quad {\Cal D},\qquad \text{where}\quad
\partial_i =\partial/\partial x^i,\tag 1.1
$$
describing the motion of a perfect compressible fluid body in vacuum:
$$
 (\partial_t+V^k\partial_{k}) \rho+ \rho \div V=0 ,\qquad \div V=\partial_{k} V^k
\qquad \text{in}\quad {\Cal D},\tag 1.2
$$
where  $V^{k \!}=\delta^{ki} v_i=v_k$ and we use the summation
convention over repeated upper and lower indices. Here the
velocity $V\!=(V^1\!,...,V^n)$, the density $\rho$ and the domain
${\Cal D}\!=\cup_{0\leq t\leq T} \{t\}\!\times {\Cal D}_t$, ${\Cal
D}_t\!\subset\Bbb R^n$ are to be determined. The pressure
$p=p\,(\rho)$ is assumed to be a given strictly increasing smooth
function of the density. The boundary $\pa{\Cal D}_t$ moves with
the velocity of the fluid particles at the boundary. The fluid
body moves in vacuum so the pressure vanishes in the exterior and
hence on the boundary. We therefore also require the boundary
conditions on
$\pa {\Cal D}=\cup_{0\leq t\leq T} \{t\}\!\times \pa {\Cal D}_t$:
$$
\align   (\partial_t+V^k\partial_{k})|_{\partial {\Cal D}}&\in T(\partial {\Cal
D}),\tag 1.3\\
p=0,\quad & \text{on}\quad \partial {\Cal D}.\tag 1.4
\endalign
$$

Constant pressure on the boundary leads to energy conservation and it is
 needed for the linearized equations
to be well posed. Since the pressure is assumed to be a strictly
increasing function of the density we can alternatively think of
the density as a function of the pressure and for physical reasons
this function has to be non negative. Therefore
 the density has to be a non negative constant
$\overline{\rho}_0$ on the boundary and we will in fact assume
that $\overline{\rho}_0\!>\!0$, which is the case of liquid.
We hence assume that
$$
p(\overline{\rho}_0)=0\qquad\text{and}\qquad   p^\prime(\rho)>0,
\qquad\text{for}\quad \rho\geq \overline{\rho}_0,
\qquad\text{where}\quad \overline{\rho}_0>0\tag 1.5
$$
From a physical point of view one can alternatively think of
the pressure as a small positive constant on the boundary.
By thinking of the density as function of the pressure
the incompressible case can be thought of
as the special case of constant density function.

The motion of the surface of the ocean is described by the above
model. Free boundary problems for compressible fluids are also
of fundamental importance in astrophysics since they describe
stars. The model also describes the case of one fluid surrounded
by and moving inside another fluid.  For large massive bodies like
stars gravity helps holding it together and for
smaller bodies like water drops surface tension helps
holding it together. Here we neglect the influence of gravity
which will just contribute with a lower order term and we neglect
surface tension which has a regularizing effect.

Given a bounded domain
${\Cal D}_0\subset \Bbb R^n$, that is homeomorphic to the unit ball,
and initial data $V_0$ and $\rho_0$,
we want to find a set
${\Cal D}\subset [0,T]\times \Bbb R^n$,
a vector field $V$ and a function $\rho$,
solving (1.1)-(1.4) and satisfying the initial conditions
$$
\align \{x;\, (0,x)\in {\Cal D}\}&={\Cal  D}_0,\tag 1.6\\
V=V_0,\quad\rho=\rho_0\quad &\text{on}\quad \{0\}\times {\Cal D}_0.\tag 1.7
\endalign
$$
In order for the initial-boundary value problem (1.1)-(1.7) to be
solvable initial data (1.7) has to satisfy certain compatibility
conditions at the boundary. By (1.2), (1.4) also implies that
$\div V\big|_{\pa{\Cal D}}=0$. We must therefore have
$\rho_0\big|_{\pa{\Cal D}_0}=\overline{\rho}_0$ and $\div
V_0\big|_{\pa{\Cal D}_0}=0$.
Furthermore, taking the divergence of (1.1) gives an equation for
$(\pa_t+V^k\pa_k )\div V$ in terms of only space derivatives of
$V$ and $\rho$, which leads to further compatibility conditions.
In general we say that initial data satisfy the compatibility
condition of order $m$ if there is a formal power series solution
in $t$, of (1.1)-(1.7) $(\tilde{\rho},\tilde{V})$, satisfying
$$
(\pa_t +\tilde{V}^k\pa_k)^j (\tilde{\rho}-\overline{\rho}_0)\big|_{\{0\}\times
\pa{\Cal D}_0}=0,\qquad\quad j=0,..,m-1\tag 1.8
$$

Let ${\Cal N}$ be the exterior unit normal to the
free surface $\pa{\Cal D}_t$.
Christodoulou\cite{C2} conjectured the initial value problem (1.1)-(1.8),
is well posed in  Sobolev spaces under the assumption
$$
\na_{\Cal N}\, p\leq -c_0 <0,\quad\text{on}\quad \partial {\Cal D},
\qquad \text{where}\quad \na_{\Cal N}={\Cal N}^i\partial_{x^i}. \tag 1.9
$$
Condition (1.9) is a natural physical condition.
It says that the pressure and hence the density is larger
in the interior than at the boundary.
Since we have assumed that the pressure vanishes or is close to zero
at the boundary this is therefore related to the fact that the pressure of a
fluid has to be positive.

In general it is possible to prove local existence for analytic data
for the free interface between two fluids.
However, this type of problem might be subject to instability in Sobolev norms,
in particular Rayleigh-Taylor instability, which occurs when a heavier fluid is
on top of a lighter fluid.
Condition (1.9) prevents Rayleigh-Taylor instability from occurring.
Indeed,  if this condition is violated  Rayleigh-Taylor instability occurs
in a linearized analysis.

In the irrotational incompressible case the physical condition
(1.9) always hold, see \cite{W1,2,CL}, and \cite{W1,2} proved
local existence in Sobolev spaces in that case. \cite{W1,2}
studied the classical water wave problem describing the motion of
the surface of the ocean and showed that the water wave is not
unstable when it turns over. Ebin\cite{E1} showed that
the general incompressible problem is ill posed in Sobolev spaces when the
pressure is negative in the interior and the physical condition is
not satisfied. Ebin\cite{E2} also announced a local
existence result for the incompressible problem with
surface tension on the boundary which has a regularizing effect
so (1.9) is not needed then.

In \cite{CL}, together with Christodoulou, we proved {\it
a priori} bounds in Sobolev spaces in the general incompressible
case of non vanishing curl, assuming the physical condition
(1.9) for the pressure. We also showed that the Sobolev
norms remain bounded as long as the physical condition hold and
the second fundamental form of the
free surface and the first order derivatives of the velocity are bounded.
Usually, existence follows from similar bounds
for some iteration scheme, but the bounds in \cite{CL}
used all the symmetries of the equation
and so only hold for modifications that preserve all the symmetries.
In \cite{L1} we showed existence for
the  linearized equations and in \cite{L3} we proved
local existence for the nonlinear incompressible problem with non
vanishing curl, assuming that (1.9) holds initially.

For the corresponding compressible free boundary problem with
non-vanishing density on the boundary, there are however in
general no previous existence or well-posedness results.
Relativistic versions of these problems have been studied in
\cite{C1,DN,F,FN,R} but solved only in special cases. The methods
used for the irrotational incompressible case use that the
components of the velocity are harmonic to reduce the equations to
equations on the boundary and this does not work in the
compressible case since the divergence is non vanishing and the
pressure satisfies a wave equation in the interior. To be able to
deal with the compressible case one therefore needs to use
interior estimates as in \cite{CL,L1}.
Let us also point out that in nature one expects fluids to be compressible,
e.g.${}_{\!\!}$  water satisfies (1.5), see \cite{CF}.
For the general relativistic equations there is no special case
corresponding to the incompressible case.
In \cite{L2} we showed existence for
the  linearized equations in the compressible case and here we prove local existence
for the nonlinear compressible problem:

\proclaim{Theorem {1.}1} Suppose that $p=p(\rho)$ is a smooth function satisfying
(1.5). Suppose also that initial data $v_0$, $\rho_0$
and ${\Cal D}_0$  are smooth satisfying the compatibility conditions (1.8)
to all orders $m$ and  ${\Cal D}_0$ is diffeomorphic
to the unit ball. Then if the physical condition (1.9) hold at $t=0$
there is a $T>0$ such that (1.1)-(1.4) and (1.6)-(1.7) has a smooth solution for
$0\leq t\leq T$. Furthermore (1.9) hold for $0\leq t\leq T$, with $c_0$
replaced by $c_0/2$.
\endproclaim

 A few remarks are in order. The existence of smooth
solutions implies existence of solutions in Sobolev spaces if one
has {\it a priori} bounds in Sobolev spaces. In the incompressible
case we had already proven {\it a priori} energy bounds in Sobolev
spaces spaces in \cite{CL} as well as a continuation result, that
the solution remains smooth as long as the physical condition is
satisfies and the solution is in $C^2$. Similar bounds in Sobolev
spaces hold also in the compressible case.
We also remark that there are initial
data satisfying the compatibility conditions to all orders, see
section 16. If only finitely many compatibility
conditions are satisfied then we get existence in $C^k$ for some
$k$.  What then is essential is that the physical condition hold
and this and  the existence time only depends on a bound of
finitely many derivatives of initial data. This makes it possible
to construct a sequence of smooth solutions converging to a
solution in Sobolev norms since we have a uniform lower
bound for the existence times, in terms of the Sobolev norm.

A few remarks about the proof are also in order. As in the incompressible case \cite{L3}
we will use the Nash-Moser technique to prove local existence.
However, because of the presence of the
boundary problem for the wave equation for the enthalpy one has
to take as many time derivatives of the equations as space
derivatives. Therefore one has to use interpolation in time as
well and one might just as well do smoothing also in time
in the application of the Nash-Moser technique although there is no
loss of regularity in the time direction. Because of this
all our constants will depend on a lower bound of the time
interval. This will make it a bit more delicate since we will also
need to choose a small time, in order that the physical and
coordinate conditions should hold. However, at the same time certain estimates
are more natural when one includes time derivatives up to the highest order.

The plan of the paper is as follows. We will assume that the reader is
somewhat familiar with the notation in \cite{L3}. We will also
assume the existence proofs given in \cite{L1,L2,L3} for the
inverse of the linearized operator so we will only prove improved
estimates here. In section 2 we formulate Euler's equations in the
Lagrangian coordinates and derive the linearized equations in
these coordinates. Here we also define a modified linearized
operator which is easier to first consider. In section 3 we define
the orthogonal projection onto divergence free vector fields, the
normal operator and decompose the equation onto a divergence free
part and a wave equation for the divergence. In section 4 we
construct the families of tangential vector fields, define the
modified Lie derivatives with respect to these and calculate its
commutators with the normal operator and other operators that
occur in the linearized equation. In section 5 we derive estimates
of derivatives of a vector field in terms of the curl the
divergence and tangential derivatives or the normal operator. In
section 6 we give tame estimates for the Dirichlet problem, in
section 7 we give tame estimates for the wave equation and in
section 8 we give tame estimates for the divergence free part.
Then in section 9 we put these estimates together to get tame
estimates for the inverse of the modified linearized operator. In
section 10 we give estimates of the enthalpy in terms of the
coordinate. In section 11 we show that the physical and coordinate
conditions can be satisfied for small times if they hold
initially. In section 12 we then get tame estimates for the
inverse of the linearized operator. In section 13 we give tame
estimates for the second variational derivative. In section 14 we
construct the smoothing operators needed for the Nash-Moser
iteration. Finally, in section 15 we construct the Nash-Moser
iteration that proves local existence of a smooth solution. In
section 16 we show that one can construct a large class of initial
data that satisfy the compatibility conditions to all orders.
Most of the steps of the proof above works for the pressure any
smooth strictly increasing functions of the density.
However, the estimates for the enthalpy in terms of the coordinate
simplified if one assume that the pressure is a linear function of the
density and we first do the proof in this case and then in section 17
give the additional estimates needed for the general case.

 \head 2. Lagrangian coordinates and the linearized
equation.\endhead Let us introduce Lagrangian coordinates in which
the boundary becomes fixed. Let $\Omega$ be a the unit ball in
$\bold{R}^n$ and let $f_0:\Omega\to {\Cal D}_0$ be a
diffeomorphism. By a theorem in \cite{DM} the volume form
$\kappa_0=\det{(\pa f_0/\pa y)}$ can be arbitrarily prescribed up
to a multiplicative constant and by a scaling of the equations we
can also assume that the volume of ${\Cal D}_0$ is that of the
unit ball.  Assume that $v(t,x)$, $p(t,x)$, $(t,x)\in {\Cal D}$
are given satisfying the boundary conditions (1.3)-(1.4). The
Lagrangian coordinates $x=x(t,y)=f_t(y)$ are given by solving
$$
{d x}/{dt}=V(t,x(t,y)), \qquad x(0,y)=f_0(y),\quad y \in
\Omega\tag 2.1
$$
Then $f_t:\Omega\to {\Cal D}_t$ is a diffeomorphism, and the
boundary becomes fixed in the new $y$ coordinates. Let us
introduce the notation
  $$\align
 \quad D_t&=\frac{\partial }{\partial t}\Big|_{y=constant}=
\frac{\partial }{\partial t}\Big|_{x=constant}
+\,  V^k\frac{\partial}{\partial x^k},\tag 2.2\\
\endalign
$$
for the material derivative. The partial derivatives
$\pa_i=\pa/\pa x^i$ can then be expressed in terms of partial
derivatives $\pa_a=\pa/\pa y^a$ in the Lagrangian coordinates. We
will use letters $a,b,c,...,f$ to denote partial differentiation
in the Lagrangian coordinates and $i,j,k,...$ to denote partial
differentiation in the Eulerian frame.

In  these coordinates Euler's equation (1.1) become
$$
\rho D_t^2 x_i+\pa_i p=0, \qquad (t,y)\in[0,T]\times\Omega\tag 2.3
$$
and the continuity equation (1.2) become
$$
D_t \rho+\rho\div V=0,\qquad (t,y)\in[0,T]\times\Omega\tag 2.4
$$
Here the pressure $p=p(\rho)$ is assumed to be smooth strictly
increasing function of the density $\rho$. With $h$, the enthalpy,
i.e. $h^\prime(\rho)=p^\prime(\rho)/\rho$ and $h=0$ when $p=0$,
(2.3) becomes
$$
D_t^2 x_i+\pa_i h=0,\quad (t,y)\in[0,T]\times\Omega\tag 2.5
$$
Since $h$ is a strictly increasing function of $\rho$ we can solve
for $\rho=\rho(h)$ as a function of $h$ and with
$e(h)=\ln{\rho(h)}$ (2.4) become
$$
D_t \, e(h)+\div V=0\tag 2.6
$$

Euler's equations are now replaced by
$$
D_t^2 x_i+\pa_i h=0,\quad (t,y)\in[0,T]\times\Omega,\qquad
\text{where}\quad \pa_i=\frac{\pa y^a}{\pa x^i}\frac{\pa}{\pa
y^a}\tag 2.7
$$
and taking the divergence of (2.7) using that $[D_t,\pa_i]=-(\pa_i
V^k)\pa_k$ we obtain a wave equation for the enthalpy
$$
D_t^2 e(h)-\triangle h-(\pa_i V^k)\pa_k V^i=0,\qquad h\Big|_{\pa
\Omega}\!\!\!=0\tag 2.8
$$
Here $e(h)$ is a given smooth strictly increasing function and
$$
\triangle h\!=\!\sum_{i}\pa_i^2 h\!=\! \kappa^{-1}\pa_a\big(\kappa
g^{ab}\pa_b h\big)\qquad\text{where}\quad
g_{ab}=\delta_{ij}\frac{\pa x^i}{\pa y^a}\frac{\pa x^j}{\pa
y^b}\tag 2.9
$$
and $g^{ab}$ is the inverse of the metric $g_{ab}$ and
$\kappa=\det{(\pa x/\pa y)} =\sqrt{\det{g}}$. The initial
conditions are
 $$\align
 x\big|_{t=0}&=f_0,\qquad \qquad \qquad D_t x\big|_{t=0}=v_0,\tag 2.10\\
 \qquad h\big|_{t=0} &=e^{-1}(\ln \rho_0),\qquad \quad
 D_t h\big|_{t=0}=-\div V_0\,/e^\prime\big(e^{-1}(\ln \rho_0)\big)\tag 2.11
\endalign
 $$
 where $e^{-1}$ is the inverse function of $e(h)$.
 If $(x,h)$ satisfies (2.7)-(2.8) with initial data of the form
 (2.10)-(2.11) then $(x,h)$ also satisfies (2.6).
 By a theorem in \cite{DM} the volume form can be arbitrarily
 prescribed so we can in fact choose it so
 $\kappa_0=\det{(\pa f_0/\pa y)}=1/\rho_0$, in which case $e(h)=-\ln\kappa$,
 since this is true when $t=0$ and since $D_t\ln\kappa=\div V$.
Hence we are left we two independent initial data, $f_0$ and
$v_0$.

In order for (2.8) to be solvable we must have the following condition on $e(h)$ and
coordinate condition:
$$\align
c_1^{-1}&\leq e^\prime(h)\leq c_1,\tag 2.12\\
\sum_{a,b} |g^{ab}|+|g_{ab}|\leq c_1^2, &\qquad\qquad
|\pa x/\pa y|^2+|\pa y/\pa x|^2\leq c_1^2 \tag 2.13
\endalign
$$
for some constant $0<c_1<\infty$. In order for (2.7) to be solvable we must have
the physical condition:
$$
\na_{\Cal N}\, h\leq -c_0 <0,\quad\text{on}\quad \partial
{\Omega}, \qquad \text{where}\quad \na_{\Cal N}={\Cal
N}^i\partial_{x^i}. \tag 2.14
$$
We assume that $e(h)$ is given smooth function of
$h$ that satisfies (2.12). The condition on the metric is true
initially since $\Omega$ is diffeomorphic to some set in
${\bold{R}^n}$ but we have to assume that the condition (2.14) is
true initially. By continuity, the conditions on the inverse of
the metric and on the enthalpy are then true also for small
times, with $c_1$ replaced by $2c_1$ and $c_0$ replaced by
$c_0/2$. In the iteration we construct we have to make sure that
the iterates are small enough and the time is small enough that
these conditions remain small for all the iterates. This will be
discussed in section 11. In order to solve the wave equation
on a bounded domain one also needs compatibility conditions on
initial data. These conditions are so that the initial conditions
are compatible with the boundary condition $h\big|_{\pa\Omega}=0$.
These conditions are:
$$
D_t^k \, h\big|_{\pa\Omega}=0,\qquad\text{when }\quad t=0,\qquad
\text{for}\quad k=0,...,m-1\tag 2.15
$$
For $k=0,1$ this is simply conditions on initial data and for
$k\geq 2$ one can use (2.8) and (2.7) to express it in terms of
lower time derivatives of $h$, so it can be calculated in terms of
the initial conditions. (2.15) is called the $m^{th}$ order
compatibility condition. We will assume that our initial data are
smooth and satisfies the $m^{th}$ order condition for all $m$.
This will be used to construct an approximate solution that
satisfy the equation to all orders as $t\to 0$, and in fact the
initial conditions will be turned into an inhomogeneous term that
vanishes to all orders as $t\to 0$. This will also be discussed
in following sections. We prove the following theorem:

\proclaim{Theorem {2.}1} Suppose that initial data $(v_0,f_0,\rho_0)$ in (2.10)-(2.11)
are smooth and the compatibility conditions
(2.15) hold for all orders $m$. Suppose also that (2.13) and (2.14)
hold when $t=0$. Then there is $T>0$ such that (2.7)-(2.8) have a
smooth solution $x\in C^\infty([0,T]\times\overline{\Omega})$, if
$e(h)$ is a smooth strictly increasing function with $e(0)=0$.
Furthermore (2.13) and (2.14) hold for
$0\leq t\leq T$ with $c_0$ replaced by $c_0/2$ and $c_1$ replaced
by $2c_1$.
\endproclaim

We remark that if $x\in
C^\infty([0,T]\times\overline{\Omega})$ and the conditions in
the theorem hold then (2.8) has a solution
$h\in C^\infty([0,T]\times\overline{\Omega})$. Theorem {1.}1 follows
from Theorem {2.}1 since initially we assumed that ${\Cal D}_0$ is
diffeomorphic to the unit ball, so (2.13) holds initially.
We will first prove Theorem {2.}1 when $e(h)=ch$ is a
linear functions of $h$. This corresponds to
$p=p(\rho)=c_0(\rho-\rho_0)$, where $c,\rho_0>0$ are positive
constants. The result is true in general when the pressure is any strictly
increasing smooth function of the density.
Only the estimates for the enthalpy in terms of the coordinates
have to be modified and we show how to do this in section 17.
The reason we first pick a linear function is
that in this case we get a linear wave equation for the enthalpy
and that simplifies the estimates and makes it more similar to the
incompressible case where we have a linear elliptic equation for
the pressure.

Let us now, define the Euler map, that will be used to find the solution of
Euler's equations. A solution of Euler's equations is given by
$\Theta(x,h)=0$, where $\Theta=(\Theta_0,\dots,\Theta_n)$ is given by
$$
\Theta_i(x,h)=D_t^2 x_i+\pa_i h,\qquad i=1,...,n, \tag 2.16
$$
and
$$
\Theta_0(x,h)=D_t^2 e(h)-\triangle h-(\pa_i V^k)\pa_k V^i, \qquad
h\big|_{\pa\Omega}=0. \tag 2.17
$$
We have assumed that our initial conditions satisfy compatibility
to all orders, i.e. there are smooth functions $(x_0,h_0)$
satisfying the initial conditions (2.10)-(2.11) and
 $$
 D_t^k\Theta (x_0,h_0)\big|_{t=0}=0,\qquad\quad \text{for
 all}\quad k\geq 0,
 \qquad\quad h_0\big|_{\pa\Omega}=0.\tag 2.18
 $$
In the process of solving $\Theta(x,h)=0$ we will only consider
functions $(x,h)$ which has the same time derivatives as
$(x_0,h_0)$ when $t=0$ and which satisfy $h\big|_{\pa\Omega}=0$.
Let us therefore introduce the notation
 $$
 C_{00}^\infty([0,T]\times\overline{\Omega})=
 \{ u\in C^\infty([0,T]\times\overline{\Omega}); \,\, D_t^k u\big|_{t=0}=0,
 \,\, \, \text{for all}\,\,\, k\geq 0\}\tag 2.19
 $$
 and for short, let
 $C_{00}^\infty=C_{00}^\infty([0,T]\times\overline{\Omega})$ and
 $C^\infty=C^\infty([0,T]\times\overline{\Omega})$.
 (2.18) then says that $\Theta(x_0,h_0)\in C_{00}^\infty$ and
  if $x-x_0\in C_{00}^\infty$ then it follows that also
 $\Theta_0(x,h_0)\in C_{00}^\infty$.
 Then by \cite{H1} we can solve the equation $\Theta_0(x,h)=0$, with
initial conditions (2.10)-(2.11) and boundary conditions
$h\big|_{\pa\Omega}=0$. The result in \cite{H1} is formulated
for vanishing initial conditions and instead an inhomogeneous term
that vanishes to all orders as $t\to 0$.
 However, we can turn the problem into this by considering
 $\tilde{\Theta}_0(\tilde{h})=\Theta_0(x,\tilde{h}+h_0)-\Theta_0(x,h_0)=-\Theta_0(x,h_0)$, with vanishing
 initial data for $\tilde{h}$. This gives a solution in
 $\tilde{h}\in C_{00}^\infty$ satisfying the boundary condition
 $\tilde{h}\big|_{\pa\Omega}=0$.

Let $(x_0,h_0)$ be the formal solution given above. For $x$
satisfying $x-x_0\in C_{00}^\infty$, we now the define the Euler
map $\Phi=(\Phi_1,\dots,\Phi_n)$ to be
$$
\Phi_i(x)=D_t^2 x_i+\pa_i h,\qquad i=1,...,n\tag 2.20
$$
where $h=\Psi(x)$ is given implicitly by solving
$$
D_t^2 e(h)-\triangle h-(\pa_i V^j)(\pa_j V^i) =0, \qquad
h\big|_{\pa\Omega}=0\tag 2.21
$$
with initial conditions (2.11). A solution of Euler's equations is given by
$$
\Phi(x)=0\tag 2.22
$$
By the preceding argument we can, in fact, find a solution $h$ to (2.21) such that $h-h_0\in C_{00}^\infty$ if $x-x_0\in C_{00}^\infty$. The reason we choose to consider
the map $\Phi(x)$ instead of $\Theta(x,h)$ is that we must make sure that $h\big|_{\pa\Omega}=0$ and
that the physical condition is satisfied, since the linearized
operator is not invertible otherwise. Alternatively, one could
also have tried to only consider $h$ satisfying these conditions,
but it seems much more difficult to preserve these conditions in
the smoothing process used in the Nash-Moser iteration. The main
work will now be to prove tame estimates for the inverse of the
linearized operator.

In the Nash-Moser iteration we will in fact solve for
$$
\tilde{\Phi}(u)=\Phi(u+x_0)-\Phi(x_0)\tag 2.23
$$
Let $F=\Phi(x_0)$, when $t\geq 0$ and $F=0$ when $t<0$, and for
$\delta\geq 0$ let $F_\delta(t,y)=F(t-\delta,y)$. Then $F_\delta\in
C_{00}^\infty$. We will solve for
 $$
 \tilde{\Phi}(u)=F_\delta-F_0\tag 2.24
 $$
The Nash-Moser theorem says that if the linearized operator is
invertible and we have tame estimates for its inverse and for the
second variation of the operator then in fact we have a solution
of (2.24) if the right hand side is small in $C^\infty_{00}$.
 But the right hand side of (2.24) tends to zero in $C^\infty$ when
$\delta\to 0$ so (2.24) has a solution for some $\delta>0$ and
hence
 $$
 \Phi(u+x_0)=0,\qquad 0\leq t\leq \delta\tag 2.25
 $$
As pointed out above at each step of the iteration we will only
have functions $u$ that vanish to all orders as $t\to 0$. This
condition can in fact be preserved by smoothing operators.

In order to solve $\Phi(x)=0$ we must show that the linearized operator is invertible.
Let us therefore calculate the linearized equations. If $\delta$ is a variation
in the Lagrangian coordinates, i.e. derivative w.r.t. a parameter
when $t$ and $y$ are fixed.
Since $[\delta, \pa /\pa y^a]=0$ it follows that
$$
[\delta,\pa_i] =\Big(\delta \frac{\pa y^a}{\pa x^i}\Big)\frac{\pa}{\pa y^a}
=-(\pa_i \delta x^l)\pa_l,\tag 2.26
$$
where we used the formula for the derivative of the inverse of a matrix
$\delta A^{-1}=-A^{-1}(\delta A)A^{-1}$.
It follows that $[\delta\! -\delta x^k\pa_k,\pa_i]\!=\!0$ and
by (2.20)
$$
D_t^2 \delta x_i-(\pa_k D_t^2 x_i)\delta x^k - \pa_i \big( \delta
x^k\pa_k h- \delta h \big) =\Phi^\prime (x)\delta x_i -\delta
x^k\pa_k \Phi_i \tag 2.27
$$
It follows from this and (2.21) that we have:
\proclaim{Lemma {2.}2} Let
$\overline{x}=\overline{x}(r,t,y)$ be a smooth function of
$(r,t,y)\in K=[-\varepsilon,\varepsilon]\times [0,T]\times\overline{\Omega} $,
$\varepsilon>0$,
such that $\overline{x}\big|_{r=0}=x$.
Then $\Phi(\overline{x})$ is a smooth function of $(r,t,y)\in K$, such that
$\pa \Phi(\overline{x})/\pa r\big|_{r=0}=\Phi^\prime( x) \delta x$,
where $\delta x=\pa \overline{x}/\pa r\big|_{r=0}$
and the linear map $L_0=\Phi^\prime(x)$ is given by
$$
\Phi^\prime (x)\delta x_i =D_t^2 \delta x_i+(\pa_i\pa_k h)\delta
x^k- \pa_i \big( \delta x^k\pa_k h- \delta h \big)\tag 2.28
$$
where
$$
D_t^2 (e^\prime(h)\delta h)-\delta x^k\pa_k D_t^2 e(h)
-\triangle\big(\delta h-\delta x^k\pa_k h\big) -2(\pa_i
V^k)\pa_k(\delta V^i\!-\delta x^l\pa_l  V^i)=0 ,\qquad \delta
h\Big|_{\pa \Omega}\!\!\!=0\tag 2.29
$$
\endproclaim

In order to use the Nash-Moser iteration scheme to obtain a solution of (2.13) we
must  show that linearized operator is invertible and that the inverse
satisfies tame estimates:

\proclaim{Theorem {2.}3} Suppose that $x_0,h_0\in C^\infty\big([0,T]\times\overline{\Omega}\big)$ is a formal solution at $t=0$,
i.e. (2.18) hold. Suppose that (2.13) and (2.14) hold when $t=0$.    Let
 $$
\||u\||_{r,\infty}
=\sup_{0\leq t\leq T} \sup_{y\in\Omega}\sum_{|\alpha|\leq r}
|\pa_{t,y}^\alpha u(t,y)|
\tag 2.30
$$

Then there is a  $T_0=T(x_0,h_0)>0$, depending only on upper bounds for
$\||x_0\||_{r_0+4,\infty}+\||h_0\||_{r_0+4,\infty}$, where $r_0=[n/2]+1$, $c_0^{-1}$ and $c_1$, such
that the following hold. If
$x-x_0\in C_{00}^\infty$, $h$ is the solution to (2.21) with $h-h_0\in C_{00}^\infty$ and
$$
T\leq T_0,\qquad \quad \||x-x_0\||_{r_0+4,\infty}\leq 1,\tag 2.31
$$
then (2.13) and (2.14) hold for $0\leq t\leq T$ with $c_0$ replaced by $c_0/2$
and $c_1$ replaced by $2c_1$. Furthermore, linearized equations
$$
\Phi^\prime(x)\delta x=\delta \Phi,\qquad \text{in}\quad
[0,T]\times\overline{\Omega},\tag 2.32
$$
where $\delta \Phi\in C_{00}^\infty$
has a solution
$\delta x\in C_{00}^\infty$.
The solution satisfies the estimates
$$
\||\delta x\||_{a,\infty}
\leq C_a\big( \||\delta \Phi\||_{a+r_0,\infty}
+\||\delta \Phi\||_{0,\infty}\, \||x-x_0\||_{a+2r_0+4,\infty} \big),\qquad a\geq 0\tag 2.33
$$
where $C_a=C_a(x_0,h_0,c_0^{-1},c_1)$ is bounded when $a$ is bounded.

Furthermore $\Phi$ is twice differentiable and the second derivative
satisfies the estimates
$$
\multline
\||{\Phi}^{\prime\prime}(x)(\delta x,\epsilon x)\||_{a,\infty}\\
\leq C_a \Big( \||\delta x\||_{a+2r_0+4,\infty}
\||\epsilon{x}\||_{0,\infty}+
 \||\delta x\||_{0,\infty}
\||\epsilon{x}\||_{a+2r_0+4,\infty}
+\||x-x_0\||_{a+3r_0+6,\infty} \||\delta x\||_{0,\infty}
\||\epsilon {x}\||_{0,\infty}\Big)
\endmultline\tag 2.34
$$
\endproclaim
Theorem {2.}1 follows from Theorem {2.}3 (or rather a version with H\"older norms.) and the Nash-Moser theorem; Theorem {15.}1.
Theorem {2.}3 follows from Lemma {11}.3 and Theorem {12}.2.
More precisely, we will first prove an $L^2$ estimate
$$
 \|\delta\dot{ x}(t)\|_r+\|\delta x(t)\|_r\leq C_r\sum_{s=1}^r \||x\||_{r+r_0+4-s,\infty}
 \int_0^t\|\delta \Phi(\tau)\|_s\, d\tau,\tag 2.35
$$
where
$$
 \|u(t)\|_r=\sum_{|\alpha|\leq r }\|\pa^\alpha_{t,y} u(t,\cdot)\|_{L^2(\Omega)}
 \tag 2.36
 $$
and we will first prove this estimate for a lower order modification,
$L_1$, of the linearized operator to be described below.
Furthermore in Theorem {9.}1 we will first prove the estimate for $L_1$ expressed in the Lagrangian coordinates, i.e. when also the vector field is expressed in the Lagrangian frame to be described below.

It will follow from estimating the solution of the wave equation
(2.29) that $\delta h$ will have the same regularity as $\delta
x$.
Since
$$
[\pa_i ,D_t^2]\delta x^i=2(\pa_i V^k)\pa_k(\delta
V^i-\delta x^l\pa_l V^i) +(\pa_k D_t V^i)\pa_i\delta x^k+ \delta
x^l\pa_l\big((\pa_i V^k)\pa_k V^i \big)\tag 2.37
$$
 we get by taking the
divergence of (2.28):
$$\multline
D_t^2 \div \delta x-\delta x^l\pa_l\big(\div D_t V\!-(\pa_i V^k)\pa_k V^i\big)
-\triangle  \big(\delta x^k\pa_k h- \delta h \big)
+2(\pa_i V^k)\pa_k(\delta V^i\!-\delta x^l\pa_l V^i)\\
=\div\,(\Phi^\prime(x)\delta x)-(\pa_i\delta x^l)\pa_l
\Phi^i-\delta x^l\pa_l \div\Phi
\endmultline \tag 2.38
$$
Hence if we add (2.29) and (2.38) we get
$$
\div\,(\Phi^\prime(x)\delta x)=D_t^2\big(\div \delta
x+e^\prime(h)\delta h\big)+(\pa_i\delta x^l)\pa_l \Phi^i\tag 2.39
$$
There is also a similar identity for the curl, see \cite{L3};
$$
\curl\,( \Phi^\prime(x)\delta x)={\Cal L}_{D_t}\curl\,\big( D_t
\,\delta x-\delta x^k\pa v_k \big)+(\pa_i\delta
x^k)\pa_j\Phi_k-(\pa_j\delta x^k)\pa_i\Phi_k\tag 2.40
$$
where ${\Cal L}_{D_t}$ is the space time Lie derivative with
respect to the vector field $D_t=(1,V)$:
$$
{\Cal L}_{D_t} \sigma_{ij}=D_t \,\sigma_{ij} +(\pa_i V^l)
\sigma_{lj}+(\pa_j V^l)\sigma_{il}\tag 2.41
$$
restricted to the space components. It can be integrated along
characteristics since it is invariant under changes of
coordinates:
$$
D_t \big( a^i_a a^j_b\sigma_{ij}\big) =a^i_a a^j_b{\Cal L}_{D_t}
\sigma_{ij},\qquad \text{where} \quad a^i_a=\pa x^i/\pa y^a \tag
2.42
$$

We now want to modify the linearized operator by adding a lower
order term so as to remove the last term on the right in (2.39) or
else replace it by a lower order term proportional to $\div \delta x$
and $ \delta h$ so that the equation
$\div\,( \Phi^\prime(x)\delta x)=0$ gives an estimate for $\div\delta x$
in terms of $\delta h$ which we have better control of than $\pa \delta x$.
In \cite{L3} we used the
modification
$$\align
L_1 \delta x^i &=L_0\delta x^i-\delta x^l\pa_l \Phi^i+\delta x^i
\div \Phi\tag 2.43\\
&=D_t^2 \delta x_i-(\pa_k D_t^2 x_i)\delta x^k - \pa_i \big(
\delta x^k\pa_k h- \delta h \big) +\delta x^i \div \Phi
\endalign
$$
where
$$ \div\Phi=D_t\div V+\triangle h+(\pa_i V^j)\pa_j
V^i=D_t\div V+D_t^2 e(h)\tag 2.44
$$
 Then
we get
$$
 \div\,(L_1 \delta x)=D_t^2\big(\div \delta
x+e^\prime(h)\delta h\big)+(\div\delta x)\div \Phi \tag 2.45
$$
$L_1$ is a lower order modification of the linearized operator,
that reduces to the linearized operator at a solution of $\Phi(x)=0$.
We will first prove that $L_1$ is invertible using that it has a
nice equation for the divergence. Then once we inverted $L_1$ we
will obtain estimates for $L_1$ that are so good that they can be
used to iterate and deal with lower order modifications like
$L_0$. However, as pointed out in \cite{Ha} one do not need to
invert the operator exactly but only do so up to a quadratic
error. But this operator $L_0$ or some modification of it is
likely to show up at other places so it seems important to show
that a more general class of operators are invertible and to have
good estimates for them.

Things are somewhat easier to see if we express the vector field in the Lagrangian frame:
$$
W^a=\frac{\pa y^a}{\pa x^i} \delta x^i, \qquad
\omega_{ab}=\frac{\pa x^i}{\pa y^a} \frac{\pa x^i}{\pa y^b}\big(
\pa_i v_j-\pa_j v_i\big). \tag 2.46
$$
Then
$$
\align D_t W^a&=\frac{\pa y^a}{\pa x^k } \big( \delta V^k -\delta
x^l\pa_l V^k\big),\tag 2.47\\
\leftalignspace
 D_t^2 W^a&= \frac{\pa y^a}{\pa x^i }
\Big(D_t^2\delta x_i -\delta x^k\pa_k D_t^2 x_i -2\big( \delta
V^k\! -\delta x^l\pa_l V^k\big)\pa_k v_i\Big).\rightalignspace\tag 2.48
\endalign
$$
Since also
$$
D_t g_{ab}=\frac{\pa x^i}{\pa y^a} \frac{\pa x^i}{\pa y^b}\big(
\pa_i v_j+\pa_j v_i\big)\tag 2.49
$$
the modified linearized operator (2.43) in the Lagrangian frame
become
$$
L_1 W^a= g_{ab} D_t^2 W^b-\pa_a\big( (\pa_c h) W^c -\delta h\big)
+\big(D_t {{g}}_{ac}-\omega_{ac} \big) D_t W^c+\div\Phi \, W^a.
\tag 2.50
$$

Also, the divergence is invariant,
$$
 \div \delta x=\div W=\kappa^{-1}\pa_a \big(\kappa W^a\big) \tag
 2.51
$$
Now, $D_t$ does not commute with taking the divergence, if
$\kappa^{-1}D_t\kappa=\div V\neq 0$, so we will replace it by a
modified time derivative that does:
$$
\hat{D}_t=D_t+\div V,\qquad\text{ie}\qquad \hat{D}_t W^a=D_t
W^a+(\div V)W^a=\kappa^{-1} D_t(\kappa W^a)\tag 2.52
$$
It then follows that
 $$
\hat{D}_t \div W=\div \hat{D}_t W\tag 2.53
 $$
  if $\sigma=\ln\kappa$ then
$\dot{\sigma}=D_t\sigma=\div V$, see \cite{L1}. With this notation
we have $\hat{D}_t^2=(D_t+\div V)(D_t+\div V)= D_t^2+2\dot{\sigma}
D_t +\dot{\sigma}^2+\ddot{\sigma}= D_t^2+2\dot{\sigma} \hat{D}_t
+\ddot{\sigma}-\dot{\sigma}^2$ so
$$
D_t^2=\hat{D}_t^2-2\dot{\sigma}\hat{D}_t+\dot{\sigma}^2-\ddot{\sigma},
\qquad D_t=\hat{D}_t-\dot{\sigma}\tag 2.54
$$
Hence, with $\dot{W}=\hat{D}_t W$ and $\ddot{W}=\hat{D}_t^2 W$, we
can write the equation (2.50) as
$$
 L_1 W^a=\ddot{W}^a-g^{ab}\pa_b\big( (\pa_c h) W^c-\delta h\big)
 -B_1 \dot{W}^a-B_0 W^a
 \tag 2.55
$$
where
$$\align
B_1\dot{W}^a&= -g^{ab}(D_t{{g}}_{bc}-\omega_{bc}) \dot{W}^c
+2\dot{\sigma}\dot{W}^a ,
\tag 2.56\\
B_0{W}^a&= g^{ab}(D_t {g}_{bc}-\omega_{bc})\dot{\sigma} W^c
-(D_t^2 e(h)+\dot{\sigma}^2)W^a , \tag 2.57
\endalign
$$

Taking the divergence of (2.55) gives
$$
\div L_1 W=\hat{D}_t^2 \div W-\triangle\big( (\pa_c h) W^c-\delta
h\big) -\div B_1 \dot W-\div B_0 W\tag 2.58
$$
On the other hand by (2.29)
$$
D_t^2 (e^\prime(h)\delta h) -\triangle\big(\delta h-(\pa_c h)
W^c\big)+\div B_1\dot{W}-2\dot{\sigma}\div \dot{W} +\div B_0
W+(D_t^2 e(h)+\dot{\sigma}^2)\div W =0\tag 2.59
$$
If we add (2.59) to (2.58) we get
$$
\div L_1 W= D_t^2\big( \div W+e^\prime(h) \delta h\big)
+\div\Phi\, \div W\tag 2.60
$$
as it should be, by (2.45). With $\varphi=\div W+e^\prime(h)\delta h$
we can hence alternatively write (2.59):
$$
\hat{D}_t^2 (e^\prime(h)\delta h) -\triangle\big(\delta h-(\pa_c
h) W^c\big)+\div B_1\dot{W}+\div B_0 W-2\dot{\sigma}\dot{\varphi}
+(D_t^2 e(h)+\dot{\sigma}^2)\varphi-\div\Phi\, e^\prime(h)\delta h
=0\tag 2.61
$$
or
$$
D_t^2 (e^\prime(h)\delta h)-\triangle\delta h+ \pa_i \big( ( \pa^i
\delta x^k)\pa_k h\big)+(\pa^i\delta x^k)\pa_i \pa_k h -2(\pa_i
V^k)\pa_k\delta V^i=0 ,\qquad \delta h\Big|_{\pa
\Omega}\!\!\!=0\tag 2.62
$$

 We now, also want to
express $L_0=\Phi^\prime(x)$ is these coordinates. In order to do
this we must transform the term $\delta x^k\pa_k \Phi^i$ in (2.25)
to the Lagrangian frame. If $\Phi^a=\Phi^i \pa y^a/\pa x^i$, then
 $(\delta x^k\pa_k \Phi^i) \pa y^a/\pa x^i=W^c\na_c \Phi^a$, where $\na_c$ is covariant
differentiation, see e.g. \cite{CL}. Hence by (2.43)
$$
L_0 W^a=L_1 W^a -B_3 W^a ,\qquad\text{where}\qquad B_3
W^a=-W^c\na_c \Phi^a+W^a\div \Phi\tag 2.63
$$

 \comment We now want to modify the time derivative and
replace it by a derivative that preserves the divergence free
condition. Let
$$
{\Cal L}_{D_t}\delta x^i=D_t\delta x^i-(\pa_k V^i)\delta
x^k,\qquad \text{and}\qquad \hat{\Cal L}_{D_t}\delta x^i={\Cal
L}_{D_t}\delta x^i+\div V\, \delta x^i \tag 2.24
$$
Then $\div \hat{\Cal L}_{D_t} \delta x=\hat{D}_t \div\delta x$,
where $\hat{D}_t=D_t+\div V$, and
$$\align
\leftalignspace {\Cal L}_{D_t}^2 \delta x^i &=D_t^2\delta x^i
-\delta x^k\pa_k D_t^2 x^i-2(\pa_k V^i) {\Cal L}_{D_t}\delta x^k
\rightalignspace
\tag 2.25\\
\hat{\Cal L}_{D_t}^2 \delta x^i&={\Cal L}_{D_t}^2 \delta x^i
+2\div V {\Cal L}_{D_t}\delta x^i+\big(D_t\div V+(\div V)^2\big)
\delta x^i \tag 2.27
\endalign
$$
Since $D_t\div V+D_t^2 e(h)=\div\Phi$ we can hence write
$$
L_1 \delta x^i=\hat{\Cal L}_{D_t}^2 \delta x^i+2(\pa_k V^i){\Cal
L}_{D_t}\delta x^k-2\div V\, {\Cal L}_{D_t} \delta x^i -(\div
V)^2\delta x^i +\pa_i \big( (\pa_k h)\, \delta x^k-\delta
h\big)-(D_t^2 e(h))\delta x^i
$$
\endcomment

\head{3. The projection and the normal operator.\\
The energy and curl estimates.}\endhead
Let us now also define the projection $P$ onto divergence free
vector fields by
$$
PU^a=U^a-g^{ab}\pa_b p_{U},\qquad \triangle p_U=\div U, \qquad
p_U\big|_{\pa\Omega}=0\tag 3.1
$$
(Here $\triangle q=\kappa^{-1}\pa_a\big( \kappa g^{ab}\pa_b
q\big)$. ) $P$ is the orthogonal projection in the inner product
$$
\langle U,W\rangle =\int_\Omega g_{ab} U^a W^b\kappa dy\tag 3.2
$$
and its operator norm is one:
$$
\|PW\|\leq \|W\|\tag 3.3
$$

For a function $f$ that vanishes on the boundary define $A_f
W^a=g^{ab}\underline{A}_f W_b $, where
$$
\underline{A}_f W_a=-\pa_a\big( (\pa_c f) W^c-q\big),\qquad
\triangle \big( (\pa_c f) W^c-q\big)=0, \qquad
q\big|_{\pa\Omega}=0\tag 3.4
$$
This is defined for general vector fields but it is symmetric in
the divergence free class. If $U$ and $W$ are divergence free then
$$
\langle U,A_f W\rangle =\int_{\pa\Omega} n_a \, U^a (-\pa_c f)
W^c\, dS \tag 3.5
$$
where $n$ is the unit conormal. If $f\big|_{\pa\Omega}=0$ then
$-\pa_c f\big|_{\pa\Omega} =(-\na_N f) n_c$. It follows that $A_f$
is a symmetric operator on divergence free vector fields, and in
particular
$$
A=A_h\tag 3.6
$$
is positive since we assumed that $-\na_N h\geq c>0$ on the
boundary. We have
$$
|\langle U,A_f W\rangle|\leq \| \na_N f/\na_N
p\|_{L^\infty(\pa\Omega)}\langle U, A U\rangle^{1/2} \langle W, A
W\rangle^{1/2}\tag 3.7
$$
Note also that $A_f$ only depends on $\na_N f$ on the boundary so
we can replace $f$ by something with the same first order
derivative that is supported in a neighborhood of the boundary. We
now want to estimate the norm of $A_f$. Now the projection has
norm one so we can drop $q$ in (3.6). If $S$ is a tangential
vector field then $S^a\pa_a \big((\pa_c f) W^c)=(\pa_c S f) W^c
+(\pa_c f){\Cal L}_S W^c$, where the Lie derivative, ${\Cal L}_S$
is defined in the next section. Furthermore if $R$ is the normal
vector fields then we can replace $f$ by the distance $d$ to the
boundary times the value of $\na_N f$ at the boundary extended to
be constant along the normal. Replacing $f$ by this function we
see that $\pa_c f$ is then independent of the radial variable so
$R^a\pa_a \big((\pa_c f) W^c)=(\pa_c f)R W^c$. In conclusion we
get that
$$
\|A_f W\|\leq C\sum_{S\in{\Cal S}}\|\na_N
Sf\|_{L^\infty(\pa\Omega)}\|W\| + C\|\na_N
f\|_{L^\infty(\pa\Omega)}\|\pa W\|\tag 3.8
$$
where ${\Cal S}$ is a spanning set of tangential vector fields.

Let us now also define the projected multiplication operators
$M_\beta$ with a two form $\beta$ by
$$
\underline{M}_\beta W_a=\underline{P}(\beta_{ab} W^b)\tag 3.9
$$
Since the projection has norm one it follows that
$$
\|M_\beta W\|\leq \|\beta\|_{\infty}\|W\|\tag 3.10
$$
Furthermore we define the operator taking vector fields to one
forms
$$
\underline{G} W_a=\underline{M}_{g} W_a=P(g_{ab} W^b)\tag 3.11
$$
Then $G$ acting on divergence free vector fields is just the
identity $I$.

Let $L_1$ be the modified linearized operator in (2.55). We now
want to project the equation
 $$
 L_1 W=F\tag 3.12
 $$
to the divergence free vector fields:
  We will decompose $L_1$ into 3 parts. We write
$$
W=W_0+W_1,\qquad W_0=PW, \qquad W_1^a=g^{ab}\pa_b q_1,\tag 3.13
$$
Then if $\dot{g}_{ab}=\check{D}_t g_{ab}$, where
$\check{D}_t=D_t-\dot{\sigma} $, we have $ \pa_a D_t \,q_1=D_t(
g_{ab} W^b)=\dot{g}_{ab} W_1^b+g_{ab}\dot{W}_1^b$ and $\pa_a D_t^2
q_1= \ddot{g}_{ab} W_1^b +2\dot{g}_{ab}
\dot{W}_1^b+g_{ab}\ddot{W}_1^b$. Hence
$$
\ddot{W}_1^a =g^{ab}\pa_b D_t^2 q_1-2 g^{ab}\dot{g}_{bc}
\dot{W}_1^c-g^{ab}\ddot{g}_{bc} W_1^c\tag 3.14
$$
Since the projection of a gradient of a function that vanishes on
the boundary vanishes it follows that
$$
P\ddot{W}_1^a=B_2(W_1,\dot{W}_1)^a,\qquad\text{where}
\qquad B_2(W_1,\dot{W}_1)^a=-P\big(2 g^{ab}\dot{g}_{bc}
\dot{W}_1^c+g^{ab}\ddot{g}_{bc} W_1^c\big)\tag 3.15
$$
Let $L_{10}=PL_1$. Since $\div \ddot{W}_0=0$ and $\delta h$
vanishes the boundary it follows from projecting (2.55) that,
$$
\align
L_{10} W_0&=\ddot{W}_0+A W_0 -B_{10} \dot{W}_0-B_{00} W_0\tag 3.16\\
L_{10} W_1&=AW_1-B_{11} \dot{W}_1-B_{01} W_1 \tag 3.17
\endalign
$$
where
$$
B_{i0}W\!=P B_{i} W,\qquad
 B_{11}{W}^a\!=PB_{1} W^a\! +2 P\big(
g^{ab}\dot{g}_{bc} {W}^c)\,\qquad
 B_{01} W^a\!=PB_{0} W^a\!
+P\big(g^{ab}\ddot{g}_{bc} W^c\big).\tag 3.18
$$
Hence the projection of (3.12) becomes
$$
L_{10} W_{0}=-L_{10}W_1+PF,\tag 3.19
$$
Here, by (2.60)
$$
W_1^{a}=g^{ab}\pa_b q_1,\qquad \triangle\,
q_1=\varphi-e^\prime(h)\delta h,\qquad q_1\big|_{\pa\Omega}=0\tag
3.20
$$
where
$$
D_t^2\varphi +\div\Phi\, \varphi=\div F+\div\Phi\, e^\prime(h)\delta
h,\tag 3.21
$$
and by (2.59)
$$
D_t^2 (e^\prime(h)\delta h) -\triangle\delta h =\triangle\big(
(\pa_c h) W^c\big) -\div B_1\dot{W}+2\dot{\sigma}\div \dot{W}
-\div B_0 W-(D_t^2 e(h)+\dot{\sigma}^2)\div W \tag 3.22
$$
Hence we have obtained a system of equations for
$(W_0,W_1,\varphi,\delta h)$.

We now want to show how to show the main idea of how to obtain
estimates for the divergence free equation.
 Let $\dot{W}=\hat{D}_t
W= D_t W+ (\div V) W= \kappa^{-1}D_t(\kappa W) =D_t
W+\dot{\sigma}W$ Then $\dot{W}$ is divergence free if $W$ is
divergence free. Let us now derive the basic energy estimate for
equations of the form
$$
\ddot{W}+AW=H\tag 3.23
$$
where $A$ is either the normal operator or the smoothed out normal
operator. Now, for any symmetric operator $B$ we have
$$
\frac{d}{dt}\langle W, BW\rangle =2\langle \dot{W}, BW\rangle
+\langle W,\dot{B}W\rangle\tag 3.24
$$
where $\dot{B}$ is the time derivative of the operator $B$
considered as an operator from the divergence free vector fields
to the one forms, see section 4. $\dot{B}$ is defined by (4.24)
with $T=D_t$. For the two operators that we will consider here
this is given by (4.39) and (4.40). Note that, the projection in
(4.24) comes up here since we take the inner product with a
divergence free vector field in (3.24). Let the lowest order
energy $E_0=E(W)$ be defined by
$$
E(W)=\langle \dot{W},\dot{W}\rangle +\langle W,(A+I)W\rangle.\tag
3.25
$$
Since $\langle W,W\rangle =\langle W,GW\rangle$, it follows that
$$
\dot{E}_0=2\langle \dot{W},\ddot{W}+(A+I)W\rangle +\langle
\dot{W},\dot{G}\dot{W}\rangle+\langle W,(\dot{A}+\dot{G})W\rangle
\tag 3.26
$$
It follows from (4.24)-(4.36) and (3.7) and (3.10) that
 $$
 |\langle
W,\dot{A}W\rangle| \leq \|\dot{h}/h\|_\infty \langle W,AW\rangle
,\qquad |\langle W,\dot{G}W\rangle |\leq \|\dot{g}\|_\infty
\langle W,W\rangle\tag 3.27
$$
The last two terms are hence bounded by a constant times the
energy so it follows that
$$
|\dot{E}_0|\leq \sqrt{E_0}\big( 2\|H\| +c\sqrt{E_0}\big), \qquad
c=\|\dot{h}/h\|_\infty+\|\dot{g}\|_\infty+2 \tag 3.28
$$

Now, let $w=\underline{W}$, $\dot{w}=\underline{\dot{W}}$ and
$\ddot{w}=\underline{\ddot{W}}$,
 i.e. $w_a=g_{ab} W^b$, $\dot{w}_a=g_{ab} \dot{W}^b$ and
$\ddot{w}_a=g_{ab}\ddot{W}^b$. Observe that $\dot{w}\neq D_t w$
etc. since $\dot{W}=\hat{D}_t W$. In fact
$$
D_t w_a= \dot{g}_{ab}W^b+\dot{w}_a, \qquad D_t
\dot{w}_a=\dot{g}_{ab}\dot{W}^b+\ddot{w}_a\tag 3.29
$$
where $\dot{g}_{ab}=\check{D}_t g_{ab}$ and $\dot{W}^a=\hat{D}_t
W^a$. Now our equation says that
$$
\ddot{w}+\underline{A} W=\underline{H},\qquad H=-B_0
W-B_1\dot{W}+F\tag 3.30
$$
Since $\curl \underline{A}W=0$ it follows that
$$
|\curl \ddot{w}|\leq C\big(|\pa W|+|W| +|\pa \dot{W}|
+|\dot{W}|+|\curl \underline{F}|\big)\tag 3.32
$$
and hence
$$\align
|D_t\curl w|&\leq C\big(|\curl\dot{w}|+|\pa W| +|W|\big),\tag 3.32\\
|D_t \curl \dot{w}|&\leq C\big(|\pa W|+|W| +|\pa \dot{W}|
+|\dot{W}|+|\curl \underline{F}|\big)\tag 3.33
\endalign
$$

\head 4. The tangential vector fields and Lie derivatives\endhead
Following \cite{L1}, we now construct the tangential vector
fields, that are time independent expressed in the Lagrangian
coordinates, i.e. that commute with $D_t$. This means that in the
Lagrangian coordinates they are of the form $S^a(y)\pa/\pa y^a$.
Furthermore, they will satisfy,
$$
\pa_a S^a=0,\tag 4.1
$$
Since $\Omega$ is the unit ball in $\bold{R}^n$ the vector fields
can be explicitly given. The vector fields
$$
y^a\pa/\pa y^b-y^b\pa/\pa y^a\tag 4.2
$$
corresponding to rotations, span the tangent space of the boundary
and are divergence free in the interior. Furthermore they span the
tangent space of the level sets of the distance function from the
boundary in the Lagrangian coordinates
$$
d(y)=\dist{(y,\pa \Omega)}=1-|y|\tag 4.3
$$
away from the origin $y\neq 0$. We will denote this set of vector
fields by ${\Cal S}_0$ We also construct a set of divergence free
vector fields that span the full tangent space at distance
$d(y)\geq d_0$ and that are compactly supported in the interior at
a fixed distance $d_0/2$ from the boundary. The basic one is
$$
h(y^3,...,y^n)\Big(f(y^1)g^\prime(y^2)\pa/\pa y^1-
f^\prime(y^1)g(y^2)\pa/\pa y^2\Big),\tag 4.4
$$
which satisfies (4.1). Furthermore we can choose $f,g,h$ such that
it is equal to $\pa/\pa y^1 $ when $|y^i|\leq 1/4$, for
$i=1,...,n$ and so that it is $0$ when $|y^i|\geq 1/2$ for some
$i$. In fact let $f$ and $g$ be smooth functions such that
$f(s)=1$ when $|s|\leq 1/4$ and $f(s)=0$ when $|s|\geq 1/2$ and
$g^\prime(s)=1$ when $|s|\leq 1/4$ and $g(s)=0$ when $|s|\geq
1/2$. Finally let $h(y^3,...,y^n)=f(y^3)\cdot\cdot\cdot f(y^n)$.
By scaling, translation and rotation of these vector fields we can
obviously construct a finite set of vector fields that span the
tangent space when $d\geq d_0$ and are compactly supported in the
set where $d\geq d_0/2$. We will denote this set of vector fields
by ${\Cal S}_1$. Let ${\Cal S}={\Cal S}_0\cup {\Cal S}_1$ denote
the family of tangential space vector fields  and let ${\Cal
T}={\Cal S}\cup \{D_t\}$ denote the family of space time
tangential vector fields.

Let the radial vector field be
$$
R=y^a\pa/\pa y^a. \tag 4.5
$$
Now,
$$
\pa_a R^a=n\tag 4.6
$$
 is not $0$ but for our purposes it suffices that
it is constant. Let ${\Cal R}={\Cal S}\cup\{R\}$. Note that ${\Cal
R}$ span the full tangent space of the space everywhere. Let
${\Cal U}={\Cal S}\cup \{R\}\cup\{D_t\}$ denote the family of all
vector fields. Note also that the radial vector field commutes
with the rotations;
$$
[R,S]=0,\qquad S\in {\Cal S}_0\tag 4.7
$$
Furthermore, the commutators of two vector fields in ${\Cal S}_0$
is just $\pm$ another vector field in ${\Cal S}_0$. Therefore, for
$i=0,1$, let ${\Cal R}_i={\Cal S}_i\cup\{R\}$, ${\Cal T}_i={\Cal
S}_i\cup\{D_t\}$ and ${\Cal U}_i={\Cal S}_i\cup\{R\}\cup\{ D_t\}$.

Let ${\Cal U}=\{ U_i\}_{i=1}^M $ be some labeling of our family of
vector fields. We will also use multindices $I=(i_1,...,i_r)$ of
length $|I|=r$. so $U^I=U_{i_1}\cdot\cdot\cdot U_{i_r}$ and ${\Cal
L}_U^I={\Cal L}_{U_{i_1}}\cdot\cdot\cdot {\Cal L}_{U_{i_r}}$,
where ${\Cal L}_U$ is the Lie derivative, defined below. Sometimes
we will write ${\Cal L}_U^I$, where $U\in {\Cal S}_0$ or $I\in
{\Cal S}_0$, meaning that $U_{i_k}\in {\Cal S}_0$ for all of the
indices in $I$.

Let us now introduce the Lie derivative of the vector field $W$
with respect to the vector field $T$;
$$
{\Cal L}_T W^a=TW^a -(\pa_c T^a) W^c\tag 4.8
$$
We will only deal with Lie derivatives with respect to the vector
fields $T$ constructed in the previous section. For those vector
fields $T$ we have
$$
 [D_t,T],\qquad\text{and}\qquad [ D_t , {\Cal L}_{T} ]=0\tag 4.9
$$
The Lie derivative of a one form is defined by
$$
{\Cal L}_T \alpha_a =T\alpha_a+(\pa_a T^c) \alpha_c,\tag 4.10
$$
The Lie derivative also commute with exterior differentiation,
$[{\Cal L}_T, d]=0$ so
$$
{\Cal L}_T \pa_a q=\pa_a T q\tag 4.11
$$
if $q$ is a function. The Lie derivative of a two form is given by
$$
{\Cal L}_T \beta_{ab} =T\beta_{ab} +(\pa_a T^c) \beta_{cb}+(\pa_b
T^c) \beta_{ac}\tag 4.12
$$
Furthermore if $w$ is a one form and $\curl w_{ab}= dw_{ab}=\pa_a
w_b-\pa_b w_a$ then since the Lie derivative commutes with
exterior differentiation:
$$
{\Cal L}_T \curl w_{ab}=\curl {\Cal L}_T w_{ab}\tag 4.13
$$
We will also use that the Lie derivative satisfies Leibniz rule,
e.g.
$$
{\Cal L}_T (\alpha_{c} W^c)= ({\Cal L}_T \alpha_{c})
W^c+\alpha_{c} {\Cal L}_T W^c,\qquad {\Cal L}_T (\beta_{ac} W^c)=
({\Cal L}_T\beta_{ac}) W^c+\beta_{ac} {\Cal L}_T W^c. \tag 4.14
$$
Furthermore, we will also treat $D_t$ as if it was a Lie
derivative and set
$$
{\Cal L}_{D_t}=D_t\tag 4.15
$$
Now of course this is not a space Lie derivative. It can however
be interpreted as a space time Lie derivative. But the important
thing is that it satisfies all the properties of the other Lie
derivatives we are considering.
 The reason we want to call it ${\Cal L}_{D_t}$
is simply a matter of that we will apply products of Lie
derivatives and $D_t$ applied to the equation and since they
behave in exactly the same way it is more efficient to have one
notation for them.



In \cite{L1} we used extensively that the Lie derivatives with
respect to the vector fields above preserved the divergence free
condition. This is no longer true if $\kappa$ is not a constant.
 since $\div U=\kappa^{-1} \pa_a(\kappa U^a)$.
This is no longer the case if $U$ is not divergence free. One
could modify the vector fields by multiplying them by
$\kappa^{-1}$. However, instead we will essentially multiply the
vector field we apply them to with $\kappa$. The modified Lie
derivative is now for any of our tangential vector fields defined
by
$$
\tilde{\Cal L}_U W={\Cal L}_U W+ (\div U) W,\qquad \tag 4.16
$$
They preserves the divergence free condition, in fact
$$
\div \tilde{\Cal L}_U W=\tilde{U}\div W, \qquad\text{where}\qquad
\tilde{U}f=Uf+(\div U)f/\tag 4.17
$$
if $f$ is a function. This definition is invariant and (4.17)
holds for any vector field $U$. However, in general, since we are
considering Lie derivatives only with respect to the vector fields
constructed above and only expressed in the Lagrangian coordinates
it is simpler to us the definition
$$
\hat{\Cal L}_U W=\kappa^{-1} {\Cal L}_U (\kappa W)= {\Cal L}_U
W+(U\sigma)W, \qquad \text{where}\qquad \sigma=\ln \kappa\tag 4.18
$$
Due to (4.1), $\div S=S\sigma$ if $S$ is any of the tangential
vector fields and $\div R=R\sigma+n $, if $R$ is the radial vector
field. For any of the tangential vector fields it then follows
that
$$
\div \hat{\Cal L}_U W=\hat{U}\div W, \qquad\text{where}\qquad
\hat{U}f=Uf+(U\sigma)f. \tag 4.19
$$
This has several advantages. The commutators satisfy $[\hat{\Cal
L}_U,\hat{\Cal L}_T] =\hat{\Cal L}_{[U,T]}$, since this is true
for the usual Lie derivative. Furthermore, this definition is
constant with our previous definition of $\hat{D}_t$.

However, when applied to one forms we want to use the regular
definition of the Lie derivative. Also, when applied to two forms
most of the time we use the regular definition: However, when
applied to two forms it turns out to be sometimes convenient to
use the opposite  modification:
$$
\check{\Cal L}_T \beta_{ab}={\Cal L}_T \beta_{ab}-(U\sigma)
\beta_{ab},\qquad \check{U}=U-(U\sigma) \tag 4.20
$$
We will most of the time apply the Lie derivative to products of
the form $\alpha_a=\beta_{ab} W^b$:
$$
{\Cal L}_T \big(\beta_{ab} W^b\big)= (\check{\Cal L}_T
\beta_{ab})W^b+\beta_{ab} \hat{\Cal L}_T W\tag 4.21
$$
since the usual Lie derivative satisfies Leibniz rule.
 Using the modified Lie derivative
we indicated in \cite{L2} how to extend the existence theorem in
\cite{L1} to the case when $\kappa $ is no longer constant, i.e.
$D_t\sigma=\div V\neq 0$. This will be carried out in more detail
here.

We will now calculate commutator between Lie derivatives and the
operator defined in the previous section, i.e. the normal operator
and the multiplication operators. It is easier to calculate the
commutator with Lie derivatives of these operators considered as
operators with values in the one forms. The one form $w$
corresponding to the vector fields $W$ is given by lowering the
indices
$$
w_a=\underline{W}_a=g_{ab}W^b\tag 4.22
$$
For an operator $B$ on vector fields we denote the corresponding
operator with values in the one forms by $\underline{B}$. These
are related by
$$
\underline{B} W_a=g_{ab} BW^b,\qquad
BW^a=g^{ab}\underline{B}_a\tag 4.23
$$
Most operators that we consider will map onto the divergence free
vector fields so we will project the result afterwards to stay in
this class. Furthermore, in order to preserve the divergence free
condition we will use the modified Lie derivative. If the modified
Lie derivative is applied to a divergence free vector field then
the result is divergence free so projecting after commuting does
not change the result. As pointed out above, for our operators it
is easier to commute Lie derivatives with the corresponding
operators from the divergence free vector fields to the one forms.
Let $B_T$ be defined by
$$
{B}_T W^a=P\big(g^{ab}\big({\Cal L}_T \underline{B}W_b
-\underline{B}_b\hat {\Cal L}_{T}W\big)\big)\tag 4.24
$$
In particular if $B$ is a projected multiplication operator
 $\underline{B}_a W=P(\beta_{ab} W^b)=\beta_{ab}W^b-\pa_a q$, where $q$ vanishes
 on the boundary is chosen so that $\div BW=0$ then
$$
{\Cal L}_T \underline{B}_{a} W= \beta_{ab}\hat{\Cal L}_T W^b+
(\check{\Cal L}_T \beta_{ab})W^b+\pa_a Tq\tag 4.25
$$
and if we project to the divergence free vector fields then the
term $\pa_a Tq$ vanishes since if $T$ is a tangential vector field
then $Tq=0$ as well. It therefore follows that $B_T$ is another
projected multiplication operator:
$$
\underline{B}_T W_a =P\big((\check{\Cal L}_T
\beta_{ab})W^b\big)\tag 4.26
$$
In particular, we will denote the time derivative of an operator
by $\dot{B}=B_{D_t}$ and for a projected multiplication operator
this is
$$
\dot{B}W=B_{D_t} W= P (\big( \check{D}_t \beta_{ab}) W^b\big)\tag
4.27
$$

If $B$ maps on to the divergence free vector fields
$$
\hat{\Cal L}_T BW^a=\hat{\Cal L}_T (g^{ab}\underline{B}_a W)
=(\hat{\Cal L}_T g^{ab})\underline{B}_a W+ g^{ab}{\Cal L}_T
\underline{B}_a W\tag 4.28
$$
Here $\hat{\Cal L}_T g^{ab}=-g^{ac} g^{bd}\check{\Cal L}_T
g_{cd}$, and if $B$ maps onto the divergence free vector fields
then $\hat{\Cal L}_T B$ is also divergence free so the left hand
side is unchanged if we do so and we get:
$$
\hat{\Cal L}_T BW^a=-P\big(g^{ab}(\check{\Cal L}_T
g_{bc})\underline{B}W^c\big) +P\big(g^{ab}\big({\Cal L}_T
\underline{B}_a W-\underline{B}_a \hat{\Cal L}_T W\big) \big)
+B\hat{\Cal L}_T W^a \tag 4.29
$$
By (4.26) applied the $\underline{G}_{ab}=P(g_{ab}W^b)$ we see
that $G_T W=P\big((g^{ab}\check{\Cal L}_T g_{bc}) W^c)\big)$ so
the first term in the right of (4.29) is $G_T BW^a$. The second
term is by definition (4.24) $B_T W$ so we get
$$
\hat{\Cal L}_T BW=B\hat{\Cal L}_T W +B_T W -G_T BW\tag 4.30
$$

 The most important property of the projection
is that it almost commutes with Lie derivatives with respect to
tangential vector fields: I.e. let $\underline{P} u_a=u_a-\pa_a
p_U$. Then
$$
\underline{P}{\Cal L}_T \underline{P} u_a=\underline{P} {\Cal
L}_{T} u_a \tag 4.31
$$
since ${\Cal L}_T \pa_a p_{U}= \pa_a T p_{U}$ vanishes when we
project again since $T p_U$ vanishes on the boundary. We have just
used this fact above. We have already calculated commutators
between Lie derivatives and the projected multiplication operators
so let us now also calculate the commutator between the Lie
derivative with respect to tangential vector fields and the normal
operator. Recall that the normal operator is defined by $A_f
W^a=g^{ab}\underline{A}_f W_b $, where
$$
\underline{A}_f W_a=-\pa_a\big( (\pa_c f) W^c-q\big),\qquad
\triangle \big( (\pa_c f) W^c-q\big)=0, \qquad
q\big|_{\pa\Omega}=0\tag 4.32
$$
and $f$ is function that vanished on the boundary. Hence since the
Lie derivative commutes with exterior differentiation:
$$
{\Cal L}_T \underline{A}_f W_a= -\pa_a\big( (\pa_c f)\hat{\Cal
L}_T W^c + (\pa_c \check{T} f)W^c+(\pa_c T\sigma) f
W^c-Tq\big)\tag 4.33
$$
However, now since $q$ vanishes on the boundary it follows that
$Tq$  also vanish on the boundary and so does $(\pa_c T\sigma) f
W^c$. Therefore the last two terms vanish when we project again so
we get
$$
P\big(g^{ab}{\Cal L}_T \underline{A}_f W_b\big)
=P\big(g^{ab}\underline{A}_f \hat{\Cal L}_T W_b\big) +
P\big(g^{ab}\underline{A}_{\check{T}f} W_b\big)\tag 4.34
$$
Let us now change notation so $A=A_h$, where $h$ is the enthalpy,
see section 3. Then we have just calculated $A_T$ defined by
(4.24) to be $A_T=A_{\check{T}h}$i, i.e.
$$
A_T=A_{\check{T} h},\qquad \text{if}\qquad A=A_{h}\tag 4.35
$$
In particular, if $T=D_t$ is the time derivative we will use the
notation $\dot{A}=A_{D_t}$ which then is
 $$
 \dot{A} W=A_{D_t }W= A_{ \check{D}_t h}W\tag 4.36
$$

We can now also calculate higher order commutators:
\demo{Definition {4.}1} If $T$ is a vector fields let $B_T$ be
defined by (4.24). If $T$ and $S$ are two tangential vector fields
we define $B_{TS}=(B_S)_T$ to be the operator obtained by first
using (4.24) to define $B_S$ and then define $(B_S)_T $ to be the
operator obtained from (4.24) with $B_S$ in place of $B$.
Similarly if $S^I=S^{i_2}\cdot\cdot\cdot S^{i_r}$ is a product of
$r=|I|$ vector fields then we define
$$
B_I= \big(\cdot\cdot\cdot
(B_{S^{i_1}})\cdot\cdot\cdot\big)_{S^{i_k}}\tag 4.37
$$
\enddemo

If $B$ is a multiplication operator
$BW^a=P\big(g^{ab}\beta_{bc}W^c\big)$ then
$$
B_I W=P\big(g^{ab}(\check{\Cal L}_T^I \beta_{bc}) W^c\big).\tag
4.38
$$
In particular if $G$ is the identity operator
$GW^a=P\big(g^{ab}g_{bc}W^c\big)$ then
$$
G_I W=P\big(g^{ab}(\check{\Cal L}_T^I g_{bc}) W^c\big).\tag 4.39
$$
If $A$ is the normal operator then
$$
{A}_I W^a=P\big( g^{ab}\pa_b\big( (\pa_c \check{T}^I p)
W^c\big)\,\big)\tag 4.40
$$

With $B_T$ as in (4.4) we have proven that if $B$ maps onto the
divergence free vector fields then
$$
\hat{\Cal L}_T BW=B W_T+B_T W-G_T BW,\qquad W_T=\hat{\Cal L}_T
W\tag 4.41
$$
Repeating this, gives for a product of modified Lie derivatives:
$$
\hat{\Cal L}_T^{I} BW=c^{I_1...I_k I} G_{I_3}\cdot\cdot\cdot
G_{I_k} B_{I_1} W_{I_2} \qquad W_J=\hat{\Cal L}_T^J W\tag 4.42
$$
where the sum is over all combinations of $I=I_1+...+I_k$, and
$c^{I_1...I_k I}$ are some constants, such that $c^{I_1...I_k
I}=1$ if $I_1+I_2=I$. Let us then also introduce the notation
$$
G^{I_1 I_2 I}=c^{I_1...I_k I} G_{I_3}\cdot\cdot\cdot G_{I_k},\tag
4.43
$$
where the sum is over all combination such that
$I_3+...I_k=I-I_1-I_2$. With this notation we can write (4.42)
$$
\hat{\Cal L}_T^{I} BW=G^{I_1 I_2 I} B_{I_1} W_{I_2}\tag 4.44
$$
where again $G_{I_1 I_2 I}=1$ if $I_1+I_2=I$. Also let
$$
\tilde{G}^{I_1...I_k I}=0,\quad \text{if} \quad I_2=I, \qquad
\text{and}\qquad \tilde{G}^{I_1...I_k I}=G^{I_1...I_k
I},\quad\text{otherwise}. \tag 4.45
$$
Then we also have
$$
\hat{\Cal L}_T^{I} BW=B W_I+\tilde{G}^{I_1 I_2 I} B_{I_1}
W_{I_2}\tag 4.46
$$

\head 5. Estimates for derivatives of a vector field in terms of
the curl, the divergence and tangential derivatives or the normal
operator.\endhead

The first part of the lemma below says that one can get a point
wise estimate of any first order derivative of a vector field by
the curl, the divergence and derivatives that are tangential at
the boundary. The second part say that one can also get estimates
in $L^2$ with a normal derivative instead of tangential
derivatives. The last part says that we can get the estimate for
the normal derivative from the normal operator. The lemma is
formulated in terms in the Eulerian frame, i.e. in terms the
original Euclidean coordinates. Later we will reformulate it in
the Lagrangian frame and then also get estimates for higher
derivatives in similar fashion.

\proclaim{Lemma {5.}1} Let $\tilde{N}$ be a vector field that is
equal to the  normal ${\Cal N}$ at the boundary $\pa{\Cal D}_t$
and satisfies $ |\tilde{N}|\leq 1$ and $|\pa\tilde{N}|\leq K$. Let
$q^{ij}=\delta^{ij}-\tilde{\Cal N}^i\tilde{\Cal N}^j$. Then
$$\align
\!\!\!\!\!\!\!\!\!\!\!\!\! |\pa \beta|^2 dx&\leq C\big(
q^{kl}\delta^{ij} \pa_k \beta_i\, \pa_l \beta_j +|\curl \beta|^2
+|\div \beta|^2 \big)
\!\!\!\!\!\!\!\!\!\!\!\!\!\!\!\!\!\!\!\!\!\!\!\!\!\tag 5.1\\
\!\!\!\!\!\!\!\!\!\!\!\!\!\int_\dt\! |\pa \beta|^2 dx&\leq
C\int_\dt\!\!\big( \delta^{ij} \tilde{\Cal N}^k \tilde{\Cal N}^l
\pa_i \beta_k\, \pa_j \beta_l +|\curl \beta|^2\!+|\div
\beta|^2\!+K^2|\beta|^2\big)\, dx
\!\!\!\!\!\!\!\!\!\!\!\!\!\!\!\!\!\!\!\!\!\!\!\!\!\!\!\!\!\!\!\!\!\!\!\!
\!\!\!\!\!\!\!\!\!\tag 5.2
\endalign
$$
Suppose that $\delta^{ij} \alpha_j$ is another vector field that
is normal at the boundary and let $A\beta_i=\pa_i( \alpha_k
\beta^k-q)$ and $q$ is chosen so that $\div A\beta=0$ and
$q|_{\pa\Omega}=0$. Then
$$
\int_\dt\!\! \delta^{ij} \alpha_k \alpha_l\,
 \pa_i \beta^k\, \pa_j \beta^l
\leq C\int_\dt\!\!\big( \delta^{ij} A\beta_i\, A\beta_j
+|\alpha|^2\big(|\curl \beta|^2\!+|\div \beta|^2\big)
\!+|\pa\alpha|^2 |\beta|^2\big)\, dx\tag 5.3
$$
\endproclaim
This lemma was proven in \cite{L3}.
 We will now state some results that were proven in \cite{L3}
 and give some further definitions.

 \demo{Definition {5.}1} For ${\Cal V}$ and ${\Cal V}^\prime$ any of the family
of vector fields introduced in section 4, and for $\beta$ a two
form, a one form, a function or a vector field we define
$$
|\beta|_r^{\Cal V}=\sum_{|I|\leq r ,\, I\in {\Cal V}} |{\Cal
L}_U^I \beta\,|,\qquad\qquad |\beta|_{r,\, \, s}^{{\Cal
V}}=\sum_{|I|\leq r ,\, I \in {\Cal V},\,  j\leq s } |{\Cal L}_U^I
D_t^j \beta\, | ,\tag 5.4
$$
where ${\Cal L}_U^I$ is a product of Lie derivatives. Furthermore
let
$$
|\beta|_r=\sum_{|\alpha|+k\leq r}|\pa_y^\alpha D_t^k\beta|, \qquad\qquad
|\beta|_{r,s}=\sum_{|\alpha|\leq r,\, k\leq s}|\pa_y^\alpha
D_t^k\beta| \tag 5.5
$$
where $|\beta|$ is the point wise norm.
Furthermore let
 $\big[\beta\,\big]_0^{\Cal V}=1$, $\big[\beta\,\big]_0=1$ and
$$
\big[\beta\,\big]_s^{\Cal V}=\sum_{s_1+...+s_k\leq s,s_i\geq 1} |\beta|_{s_1}^{\Cal
V} \cdot\cdot\cdot |\beta|_{s_k}^{\Cal V},
\qquad\quad
\big[\beta\,\big]_s=\sum_{s_1+...+s_k\leq s,s_i\geq 1} |\beta|_{s_1}
\cdot\cdot\cdot |\beta|_{s_k},
\tag 5.6
$$
\enddemo
We note that if $\beta$ is a function then ${\Cal L}_U \beta
=U\beta$ and in general it is equal to this plus terms
proportional to $\beta$. Hence (5.4) is equivalent to
$\sum_{|I|\leq r,\,I\in {\Cal V}} |U^I \beta| $ . In particular if
${\Cal R}$ denotes the family of space vector fields then
$|\beta|_r^{\Cal R}$ is equivalent to $|\beta|_r$ with a constant
of equivalence independent of the metric. Note also that if
$\beta$ is the one form $\beta_a=\pa_a q$ then ${\Cal L}_U^I\beta=
\pa U^I q$ so $|\pa q|_r^{\Cal V}=\sum_{|I|\leq r,\, I\in {\Cal
V}}|\pa U^I q|$.

\demo{Definition {5.}2 } Let $c_1$ be a constant such that
$$
\sum_{a,b}\big(|g_{ab}|+|g^{ab}|\big)\leq c_1^2,\qquad\qquad |\pa
x/\pa y|^2+|\pa y/\pa x|^2\leq c_1^2\tag 5.7
$$
Let $K_{i0}$ respectively $K_i$, for $i\geq 1$, denote continuous
increasing functions of
 $$
 c_1+\||x\||_{i,0,\infty}
 ,\qquad \text{respectively}\qquad  \||x\||_{i,\infty}+c_1+T^{-1},\tag 5.8
 $$
which also depends on the order $r$ of differentiation. Here, the
norms are as in Definition {5.}3.
\enddemo
We note that the bounds for the metric $g$ and its inverse follows
from the bounds for the the Jacobian of the coordinate and its
inverse since $g_{ab}=\delta_{ij}(\pa x^i/\pa y^a)(\pa x^j/\pa
y^b)$ and $g^{ab}=\delta^{ij}(\pa y^a/\pa x^i)(\pa y^b/\pa x^j)$.
It also follows that we have a bound for $\kappa=\det{(\pa x/\pa
y)}=\sqrt{\det{g}}$ and its inverse. Furthermore, if we have a
bound for $\kappa^{-1}$ the bound for the inverse of the metric
follows from the bound for the metric, since an upper bound for
the eigenvalues and a lower bound for their product gives a lower
bound for the eigenvalues, since they are all positive by
assumption. Moreover we note that also a  lower bound
$|\pa x/\pa y|^2\geq c_1^{-2}$ follows.

In what follows it will be convenient to consider the norms of
$\hat{\Cal L}_U^I W=\kappa^{-1}{\Cal L}_U^I (\kappa W)$ if $W$ is
a vector field and of $\check{\Cal L}_U^I g=\kappa {\Cal
L}_U^I(\kappa^{-1} g)$, if $g$ is the metric. The reason for this
is simply that $\div (\hat{\Cal L}_U^I W) =\hat{U}^I \div W$ and
${\Cal L}_U^I\curl w=\curl\, ({\Cal L}_U^I W)$ and when we lower
indices $w_a=g_{ab} W^b=(\kappa^{-1} g_{ab})(\kappa W^b)$ and
apply the Lie derivative to the product we get ${\Cal L}_U w_a
=(\check{\Cal L}_U g_{ab})W^b +g_{ab}\hat{\Cal L}_U W^b$.

The following Lemma was proven in \cite{L3} for the cases without
time derivatives. The cases with time derivatives follows from
these by applying them to $W$ replaced by $\hat{D}_t^k W$, noting
that $\div \hat{D}_t W=\hat{D}_t \div W$ and $\curl \,
(\hat{D}_t^k W)_\flat$ is equal to $D_t^k \curl w$ plus terms that
are lower order in time. This is similar to the calculation in
\cite{L1,L3} that $\curl\, (\hat{\Cal L}_S^I W)_\flat$, where
$(\hat{\Cal L}_S^I W)_{\flat b}=g_{ab} \hat{\Cal L}_S^I W^b$, is
equal to ${\Cal L}_S^I \curl w$ plus terms of lower order. We
have:
 \proclaim{Lemma {5.}2}
 Let $W$ be a vector field and let $w_a=g_{ab} W^b$ be the
corresponding one form. Let $\kappa=\det{(\pa x/\pa y)}=
\sqrt{\det{g}}$. Then
$$
|\kappa|+|\kappa^{-1}|\leq K_{10}, \qquad\qquad |U^I\kappa|+|U^I
\kappa^{-1}|\leq K_{10}\, c_{I_1...I_k} |U^{I_1}g|\!
\cdot\cdot\cdot |U^{I_k} g|\tag 5.9
$$
where the sum is over all $I_1+...+I_k=I$ and $K_{10}$ is as in
Definition {5.}2.

With notation as in Definition {5.}1 and section 4 we have
$$
|\kappa W|_r^{\Cal R} \leq K_{10}\big( |\curl w|_{r-1}^{\Cal
R}+|\kappa\div W|_{r-1}^{\Cal R} +|\kappa W|_r^{\Cal
S}+\sum_{s=0}^{r-1} |g/\kappa|_{r-s}^{\Cal R}
 | \kappa W|_{s}^{\Cal R}\big)\tag 5.10
$$
We also have
$$
|\kappa W|_{r}^{\Cal R}\leq K_{10}\sum_{s=0}^{r} \big[\, g/\kappa\big]_{r-s}^{\Cal
R} \big(|\curl w|_{s-1}^{\Cal R}+|\kappa\div
W|_{s-1}^{\Cal R}+|\kappa W|_{s}^{\Cal S}\big),\tag 5.11
$$
where the first two terms in the sum should be interpreted as $0$
if $s=0$.

The inequalities (5.10)-(5.11) also hold with $({\Cal R},{\Cal S})$
replaced by $({\Cal U},{\Cal T})$ and moreover,
$$
|\kappa W|_{r,\, s}^{\Cal R} \leq K_{10}\big( |\curl
w|_{r-1,s}^{\Cal R}+|\kappa\div W|_{r-1,s}^{\Cal R} +|\kappa
W|_{r,s}^{\Cal S}+\!\!\!\!\!\!\!\! \sum_{i+j\leq r+s-1,i\leq
r,j\leq s} \!\!\!\!\!\!\!\!\!\!\!\!\! |g/\kappa|_{r+s-i-j}^{\Cal
U}
 | \kappa W|_{i,j}^{{\Cal R}}\big)\tag 5.12
$$
and
$$
|\kappa W|_{r,s}^{\Cal R}\leq K_{10}\!\!\!\sum_{i\leq r, j\leq
s}\!\!\!\! \big[\,g/\kappa\big]_{r+s-i-j}^{\Cal U} \big(|\curl
w|_{i-1,j}^{\Cal R}+|\kappa\div W|_{i-1,j}^{\Cal R }+|\kappa
W|_{i,j}^{\Cal S}\big),\tag 5.13
$$
where the sums are only over positive indices and if $i=0$ then
the first two terms should be interpreted as $0$.
\endproclaim

\demo{Definition {5.}3} For ${\Cal V}$ any of the family of vector
fields in section 4 let
$$
\| W\|_{r}^{\Cal V}=\sum_{|I|\leq r,\, I\in{\Cal V}} \|{\Cal
L}_U^I W\|, \qquad\qquad \| W\|_{r,\infty}^{\Cal V}= \sum_{|I|\leq
r,\, I\in{\Cal V}} \|{\Cal L}_U^I W\|_{\infty} \tag 5.14
$$
and let
$$
\|W\|_r=\sum_{|\alpha|+k\leq r}\|D_t^k\pa_y^\alpha W\|\qquad\qquad
\|W\|_{r,\infty}=\sum_{|\alpha|+k\leq r}\| D_t^k \pa_y^\alpha
W\|_\infty \tag 5.15
$$
where $\|W\|=\|W(t)\|=\|W(t,\cdot)\|_{L^2(\Omega)}$,
$\|W\|_\infty=\|W(t)\|_\infty=\|W(t,\cdot)\|_{L^\infty(\Omega)}$.
Let the mixed norms be defined by
$$
\|W\|_{r,s}^{\Cal V}=\sum_{j=0}^s \|D_t^j\, W\|_{r}^{\Cal V},
\qquad \|W\|_{r,s}=\sum_{j=0}^s\sum_{|\alpha|\leq r} \|D_t^j\,
\pa_y^\alpha W\|\tag 5.16
$$
and
$$
\|W\|_{r,s,\infty}^{\Cal V}=\sum_{j=0}^s \|D_t^j\,
W\|_{r,\infty}^{\Cal V}, \qquad
\|W\|_{r,s,\infty}=\sum_{j=0}^s\sum_{|\alpha|\leq r} \|D_t^j\,
\pa_y^\alpha W\|_\infty\tag 5.17
$$
 Furthermore, let
$$
 \||\alpha\||_{r,\infty}=\sup_{0\leq t\leq T}
 \|\alpha(t,\cdot)\|_{r,\infty},
 \qquad \||\alpha\||_{r,s,\infty}=\sup_{0\leq t\leq T}
 \|\alpha(t,\cdot)\|_{r,s,\infty},
 \tag 5.18
 $$
and
$$
\big[\big[\big[\beta\big]\big]\big]_{r,\infty}
 =\sum_{r_1+...+r_k\leq r,\, r_i\geq 1}\||\beta\||_{r_1,\infty}\cdot\cdot\cdot
\||\beta\||_{r_k,\infty}\tag 5.19
$$
\enddemo
It follows from the discussion after Definition {5.}1 and the
beginning of Lemma {5.}3 that $\|W\|_r$ is equivalent to
$\|W\|_r^{\Cal R}$ with a constant of equivalence depending only
on the dimension. As with the point wise estimates it will
sometimes be convenient to instead use $\|\hat{\Cal L}_U^I W\|
=\|\kappa^{-1}{\Cal L}_U^I(\kappa W)\|$. This in particular true
for the family of space tangential vector fields ${\Cal S}$.
However instead of introducing a special notation we then write
$\|\kappa W\|_r^{\Cal S}$, since $\kappa$ is bounded from above
and below by a constant $K_1$ this is equivalent with a constant
of equivalence $K_1$. Furthermore, by interpolation $\|\kappa
W\|^{\Cal S}_r\leq K_1 (\||g\||_r\|W\|+\|W\|_r^{\Cal S})$ and
$\|W\|^{\Cal S}_r\leq K_1 (\||g\||_r\|W\|+\|\kappa W\|_r^{\Cal
S})$, and in our inequalities we will have lower order terms of
this form anyway. Note also that it follows from the interpolation inequalities in Lemma {5.}7
that
$$
\||\big[\, g\big]_r\||_{0,\infty}\leq \big[\big[\big[ \, g\big]\big]\big]_{r,\infty}\leq K_1 \||g\||_{r,\infty}\tag 5.20
$$
Here the constant $K_1$ depends on a lower bound for $T>0$ and in fact tends to infinity
as $T\to 0$. Most of the time this does not matter, but when estimating lower norms
of the enthalpy we can not use this estimate.

The following lemmas were proven in \cite{L3} for the cases
without time derivatives. Their generalizations to including time
derivatives is immediate.
 \proclaim{Lemma {5.}3} We have with a
constant $K_{10}$ as in Definition {5.}2, we have
$$
\|\kappa W\|_{r,0}\leq K_{10}\big(\|\curl w\|_{r-1,0}+\|\kappa
\div W\|_{r-1,0}+\|\kappa W\|_{r,0}^{\Cal S} +K_{10}
\sum_{s=0}^{r-1}\|g\|_{r-s,0,\infty}\|\kappa W\|_{s,0}\big)\tag
5.21
$$
and
$$
\|\kappa W\|_{r,0}\leq K_{10}\sum_{s=0}^r \|g\|_{r-s,0,\infty}
\big(\|\curl w\|_{s-1,0}+\|\kappa \div W\|_{s-1,0}+\|\kappa
W\|_{s,0}^{\Cal S}\big)\tag 5.22
$$
where the two terms with $s-1$ should be interpreted as $0$ when
$s=0$. We also have with a constant $K_1$ as in Definition {5.}2,
we have
$$
\|\kappa W\|_{r}\leq K_{10}\big(\|\curl w\|_{r-1}+\|\kappa \div
W\|_{r-1}+\|\kappa W\|_{r}^{\Cal T} +K_{10}
\sum_{s=0}^{r-1}\big[\big[\big[ g\big]\big]\big]_{r-s,\infty}\|\kappa W\|_{s}\big)\tag 5.23
$$
and
$$
\|\kappa W\|_{r}\leq K_{10}\sum_{s=0}^r\big[\big[\big[ g \big]\big]\big]_{r-s,\infty}
\big(\|\curl w\|_{s-1}+\|\kappa \div W\|_{s-1}+\|\kappa
W\|_{s}^{\Cal T}\big)\tag 5.24
$$
where the two terms with $s-1$ should be interpreted as $0$ when
$s=0$. Furthermore
$$
\|\kappa W\|_{r,s}\leq K_{10}\big(\|\curl w\|_{r-1,s}\!+\|\kappa \div
W\|_{r-1,s}\!+\|\kappa W\|_{r,s}^{\Cal S}\! +
\!\!\!\!\!\!\!\!\!\sum_{i\leq r,\, j\leq s,\, i+j\leq r+s-1
}\!\!\!\!\!\!\!\!\!\!\!\!\!\big[\big[\big[ g\big]\big]\big]_{r+s-i-j,\infty}\|\kappa
W\|_{i,j}\big)\tag 5.25
$$
and
$$
\|\kappa W\|_{r,s}\leq K_{10}\!\!\!\!\!\!\sum_{i\leq r,\, j\leq s
}\!\!\!\!\!\!\big[\big[\big[ g\big]\big]\big]_{r+s-i-j,\infty}\ \big(\|\curl
w\|_{i-1,j}+\|\kappa \div W\|_{i-1,j}+\|\kappa W\|_{i,j}^{\Cal
S}\big) \tag 5.26
$$
where the two terms with $i-1$ should be interpreted as $0$ if
$i=0$.
\endproclaim

The following lemma, proven in \cite{L3}, gives a bounds of
derivatives of a vector field by the curl, the divergence and the
normal operator. We have:
  \proclaim{Lemma {5.}4} Let $c_0>0$ be a
constant such that $|\na_N h|\geq c_0>0$, let $K_1^\prime $ and
$K_2^\prime$ be constants such that $\|\na_N
h\|_{L^\infty(\pa\Omega)}+c_1\leq K_1^\prime$ and $\sum_{S\in{\Cal
S}}\|\na_N S h\|_{L^\infty(\pa\Omega)}+|\pa g|+c_1+K_1^\prime\leq
K_2^\prime$ Then
$$
c_0\| \pa W\|\leq \big( \|AW\|+ K_1^\prime(\|\curl w\|+\|\div
W\|)+K_2^\prime\|W\|\big)\tag 5.27
$$
\endproclaim
By Lemma {5.}4 we also have
$$
c_0\|\pa \hat{\Cal L}_T^J W\| \leq K_2^\prime\big(\|\curl
\,(\hat{\Cal L}_T^J W)_\flat\|+\|\div\hat{\Cal L}_T^J W\| +\|
A\hat{\Cal L}_T^J W\|+\|\hat{\Cal L}_T^J W\|\big)\tag 5.28
$$
Here $(\hat{\Cal L}_T^J W)_{\flat b}=g_{ab} \hat{\Cal L}_T^J W^b$,
and as in the proof of Lemma {5.}3 in \cite{L3} we see that
$\curl\, (\hat{\Cal L}_T^J W)_{\flat}$ is equal to $\curl {\Cal
L}_T^J w={\Cal L}_T^J \curl w$, plus terms of lower order, where
$w_a=g_{ab}W^b$. In particular we see that we can get any space
tangential derivative in this way so we also get:
 \proclaim{Lemma {5.}5} With $K_2^\prime$ as in Lemma {5.}4 we have
$$
c_0 \|W\|_{r,0}\leq  K_2^\prime\big(\|\curl w\|_{r-1,0} +\|\div
W\|_{r-1,0}+\| W\|_{r-1,A}^{\Cal S} +\sum_{s=0}^{r-1}
\|g\|_{\infty,r-s}\|W\|_{s,0}\big)\tag 5.29
$$
where
$$
\|W\|_{s,A}^{\Cal S} =\sum_{|I|=s,I\in{\Cal S}}\| A\hat{\Cal
L}_S^I W\|\tag 5.30
$$
We also have
$$
c_0 \|W\|_{r}\leq  K_2^\prime\big(\|\curl w\|_{r-1} +\|\div
W\|_{r-1}+\| W\|_{r-1,A}^{\Cal T}+\|\dot{W}\|_{r-1}^{\Cal T}
+\sum_{s=0}^{r-1} \||g\||_{\infty,r-s}\|W\|_s \big)\tag 5.31
$$
where
$$
\|W\|_{s,A}^{\Cal T} =\sum_{|I|=s,I\in{\Cal T}}\| A\hat{\Cal
L}_T^I W\|\tag 5.32
$$
\endproclaim

Let us now state the interpolation inequalities that we will use.
For a proof of Lemma {5.}6 see e.g. \cite{H1,H2} for the
$L^\infty$ estimate and \cite{CL} for the $L^2$ estimate.
\proclaim{Lemma {5.}6} There are constants $C_{r}$ respectively
$C_{r,T}$ depending on $r$ respectively on $r$ and $T$, such that
if $\beta,W\in C^\infty\big([0,T]\times\overline{\Omega}\big)$ is
a two form, a function or a vector field, then
$$\align
\||\beta\||_{s,\infty}&\leq
 C_{r,T}\||\beta\||_{0,\infty}^{1-s/r}\||\beta\||_{r,\infty}^{s/r}\tag 5.33\\
\|W\|_{s}&\leq C_r\|W\|_{0}^{1-s/r}\|W\|_{r}^{s/r}\tag 5.34
\endalign
$$
Here $\|W\|_s=\|W(t)\|_s$. Furthermore, for $\beta\in
C_{00}^\infty([0,T]\times\overline{\Omega})$, see (2.19), (5.33)
holds with constants independent of $T$.
\endproclaim
The consequence we will use is that
 \proclaim{Lemma {5.}7} There are constants $C$ depending
on $T$ and $r$, such that if  $f,f_i,\alpha,\beta,W\in
 C^\infty\big([0,T]\times\overline{\Omega}\big)$ are functions,
 two forms or vector fields, then
$$
\align\leftalignspace
 \||f\||_{s_1,\infty}\!\cdot\cdot\cdot
\||f\||_{s_k,\infty}&\leq C\||f\||_{0,\infty}^{k-1}
\||f\||_{r,\infty}\tag 5.35\\
\leftalignspace \||\alpha\||_{r-s,\infty}\||\beta\||_{r,\infty}&
\leq C\big(\||\alpha\||_{0,\infty}\||\beta\||_{s,\infty}
+\||\beta\||_{0,\infty}\||\alpha\||_{r,\infty}\big)\tag 5.36\\
\||\beta\||_{r-s,\infty}\|W\|_{s}&\leq C\big(
\||\beta\||_{0,\infty}\|W\|_r
+\||\beta\||_{r,\infty}\|W\|_0\big)\tag
5.37\\
\leftalignspace \||f_1\||_{s_1,\infty}\!\cdot\cdot\cdot
\||f_k\||_{s_k,\infty}&\leq
C\sum_{i=1}^k\||f_1\||_{0,\infty}\!\cdot\cdot\cdot
\||f_{i-1}\||_{0,\infty} \||f_i\||_{r,\infty}
\||f_{i+1}\||_{0,\infty}\!\cdot\cdot\cdot
\||f_k\||_{0,\infty}\rightalignspace \tag 5.38
\endalign
$$
where $r=s_1+...+s_k$. Here $\|W\|_s=\|W(t)\|_s$.
 Furthermore, for $\alpha,\beta,f,f_i\in
C_{00}^\infty([0,T]\times\overline{\Omega})$, see (2.19),
(5.35)-(5.38) holds with constants independent of $T$.
\endproclaim
\demo{Proof} This follows from using Lemma {5.}6 on each factor
and then using the inequality $A^{s/r} B^{1-s/r}\leq A+B$, see
\cite{L3}. \qed\enddemo

\head 6. Tame estimates for the Dirichlet problem.\endhead

In this section we will give tame estimates for the Dirichlet
problem
 $$
 \triangle q= f, \qquad\text{in}\qquad [0,T]\times\Omega,
 \qquad \quad q\big|_{\pa\Omega}=0\tag 6.1
 $$
Let us recall the definition of the mixed norms
$$
\| q\|_{r,s}=\sum_{|\alpha|\leq r,\, k \leq s}\|{D}_t^k
\pa_y^\alpha q\| \tag 6.2
$$
In the proofs that follow we will use the interpolations
$$
\||g\||_{s,\infty}\||g\||_{r,\infty}\leq K_{1} \||g\||_{s+r,\infty},
\qquad\qquad
\||g\||_{s+1,\infty}\||g\||_{r+1,\infty} \leq K_{2} \||g\||_{s+r+1,\infty} \tag 6.3
$$
where $K_{1}$ and $K_2$ are as in Definition {5.}2.

\proclaim{Theorem {6.}1 } Suppose that $q$ is a solution of the
Dirichlet problem, $q\big|_{\pa\Omega}=0$ and $W_1^a=g^{ab}\pa_b q$
then
$$
\|W_1\|_{r+1,s}+\|\pa q\|_{r+1,s}\leq K_{10}\!\!\!\!\sum_{i\leq r,\, j\leq s}\!\!\!\!
\big[\big[\big[ g\big]\big]\big]_{r+s-i-j,\infty}\,  \|\triangle q\|_{i,j}
+K_{10}\sum_{j=0}^s\big[\big[\big[ g\big]\big]\big]_{r+s+1-j,\infty}
\|W_1\|_{0,j}  \tag 6.4
$$
if $r,s\geq 0$, and
$$
\|W_1\|_{r+1,s}+\| q\|_{r+2,s}\leq K_{20}\sum_{i\leq r,\, j\leq s}
\big[\big[\big[ (\pa g,D_t g)\big]\big]\big]_{r+s-i-j}\,  \|\triangle q\|_{i,j} \tag 6.5
$$
 Moreover if $P$ is the orthogonal projection onto
divergence free vector fields and $W$ is any vector field then,
$$
\|PW\|_{r,s} \leq K_{10} \sum_{i\leq r,\, j\leq s}
\big[\big[\big[ g\big]\big]\big]_{r+s-i-j,\infty}\, \|W\|_{i,j}\tag 6.6
$$

Here $K_{10}$ and $K_{20}$ are as in Definition {5.}2,
$$
\big[\big[\big[ g\big]\big]\big]_{r,\infty}=\sum_{r_1+...+r_k\leq r,\, r_i\geq 1}
 \||g\||_{r_1,\infty}\cdot\cdot\cdot
\||g\||_{r_k,\infty}\leq K_1 \|| g\||_{r,\infty} \tag 6.7
$$
and
$$
\big[\big[\big[ (\pa g,D_t g)\big]\big]\big]_{r,\infty}
=\sum_{r_1+...+r_k\leq r,\, r_i\geq 1}
 \|| g\||_{r_1+1,\infty}\cdot\cdot\cdot
\||g\||_{r_k+1,\infty}\leq K_2 \|| g\||_{r+1,\infty} \tag 6.8
$$
\endproclaim
\demo{Remark} Note that the constant in (6.4)-(6.6) are independent of $T$
whereas those in (6.7)-(6.8) depends on a lower bound for
$T>0$.
\enddemo

First two useful lemmas: \proclaim{Lemma {6.}2} Suppose that $S\in
{\Cal S}$ and $q\big|_{\pa\Omega}=0$, and
$$
\hat{\Cal L}_S W^a=g^{ab}\pa_b q+F^a .\tag 6.9
$$
Then
$$
\| \hat{\Cal L}_S W\|\leq K_{10}(\| \div W\|+\|F\|)\tag 6.10
$$
\endproclaim
\demo{Proof} Let $W_S=\hat{\Cal L}_S W$.
$$
\int_\Omega g_{ab}  W_S^a W_S^b \, \kappa \, dy =\int_\Omega W_S^a
\pa_a  q\, \,\kappa dy +\int_\Omega W_S^a  g_{ab} F^b\kappa dy
\tag 6.11
$$
If we integrate by parts in the first integral on the right, using
that $q$ vanishes on the boundary we get
$$
\int_\Omega W_S^a\, \pa_a q\, \kappa\, dy =-\int_\Omega \div\,(
W_S)\, q \,\kappa \, dy\tag 6.12
$$
However $\div W_S=\hat{S} \div W$.  Then we can integrate by parts
in the angular direction. $S=S^a\pa_a$, $\hat{S}=S+\div S$ so we
get $\int_\Omega ( \hat{S} f)\, \kappa dy= \int_\Omega \pa_a \big(
S^a f\kappa )\, dy=0$, where $\pa_a S^a=0$.
 Hence we get
$$
\int_\Omega W_S^a\, \pa_a q\, \kappa\, dy =\int_\Omega \div\,(
W)\, ({S} q )\kappa \, dy\tag 6.13
$$
Here $| {S}  q|\leq |\pa q|\leq K_{10}(|W_S|+|F|)$ so it follows that
$$
\|W_S\|^2\leq K_{10} \|\div W\|\big( \|W_S\| +\|F\| \big) +K_{10}
\|W_S\| \| F\|\tag 6.14
$$
and using that the inequality $A^2\leq K(A+B)B$ implies that $A\leq (2K+1)B$ we get for some other constant $K_{10}$:
$$
\|W_S\|\leq K_{10} (\| \div W\|+\|F\|)\qed \tag 6.15
$$
\enddemo
\proclaim{Lemma {6.}3} Let $W^a=g^{ab} w_b $ Then
$$
\hat{D}_t^k W^a =g^{ab} D_t^k w_b-\sum_{i=0}^{k-1} \binom{k}{i}
g^{ab}(\check{D}_t^{k-i} g_{bc})\hat{D}_t^i W^c\tag 6.16
$$
\endproclaim
\demo{Proof} We have $D_t^k w_b=D_t^k\big( \kappa^{-1} g_{bc} \,
\kappa W^c\big)=\sum_{i=0}^k\binom{k}{i} \big( D_t^{k-i}
(\kappa^{-1} g_{bc})\big)\hat{D}_t^{i} (\kappa
W^{c})$, which proves the lemma.
 \qed\enddemo

\demo{Proof of Theorem {6.}1 }
To simplify notation let us denote $W^a=W_1^a=g^{ab} \pa_b q$.
If we apply ${\Cal L}_T^I$ to $w_a=g_{ab}W^b$
we get
$$
\pa_a T^I q=g_{ab}W_I^b+ \tilde{c}^{I_1 I_2} {g}_{I_1\,ab}
W_{I_2}^b, \qquad \quad W_I=\hat{\Cal L}_T^I W, \qquad
{g}_{I\,ab}=\check{\Cal L}_T^{I} g_{ab} \tag 6.17
$$
and the sum is over all combinations $I=I_1+I_2$, $\tilde{c}^{I_1
I_2}$ are constants such that $\tilde{c}^{I_1 I_2}=0$ if $I_2=I$.
We assume that $\hat{\Cal L}_T^I=\hat{\Cal L}_S^J \hat{\Cal
L}_{D_t}^s$, where $J\in{\Cal S}$ and $|J|=r+1$. If we write
$\hat{T}^I=\hat{S} \hat{T}^K$, $W_I= \hat{\Cal L}_S W_K$, and use Lemma {6.}2 we get
since $\div W_K =\hat{T}^K \div W= \hat{T}^K \triangle q
=\kappa^{-1} T^K (\kappa \triangle q)$
$$
\|W_I\|\leq K_{10}\|\hat{T}^K \triangle q\| + K_{10}\tilde{c}^{I_1
I_2}\|{g}_{I_1}\|_\infty \|W_{I_2}\|\tag 6.18
$$
or if we sum over all of them
$$
\|\kappa W\|^{\Cal S}_{r+1,s}\leq K_{10}\|\kappa \triangle  q\|_{r,s}  +K_{10}
\sum_{i\leq r+1,\, j\leq s, \,i+j\leq
r+s}\big[\big[\big[ g\big]\big]\big]_{r+s+1-i-j,\infty} \|\kappa W\|_{i,j}^{\Cal S}, \tag
6.19
$$
 We now want to apply Lemma {5.}3 to $W^a=g^{ab}\pa_b q$. Then
$\curl w=0$ and $\div W=\triangle q$.
$$ \|\kappa W\|_{r+1,s}\leq K_{10}\|\kappa\triangle q\|_{r,s}+K_1
\|W\|_{r+1,s}^{\Cal S} +K_{10}\!\!\!\!\! \sum_{i\leq r+1,\, j\leq s, \,i+j\leq
r+s}\!\!\!\!\!\big[\big[\big[ g\big]\big]\big]_{r+s+1-i-j,\infty}
 \|\kappa W\|_{i,j}, \tag 6.20
$$
(6.19) and (6.20) therefore gives for $r,s\geq 0$,
$$ \|\kappa W\|_{r+1,s}\leq K_{10}\|\kappa\triangle q\|_{r,s}
+K_{10} \sum_{i\leq r+1,\, j\leq s, \,i+j\leq
r+s}\big[\big[\big[ g \big]\big]\big]_{r+s+1-i-j,\infty} \|\kappa W\|_{i,j}, \tag 6.21
$$
Using induction it follows that
$$
\|\kappa W\|_{r+1,s}\leq K_{10} \sum_{i\leq r,\, j\leq s}
\big[\big[\big[ g\big]\big]\big]_{r+s-i-j,\infty}\|\kappa\triangle q\|_{i,j}
 +K_{10} \sum_{j=0}^s
\big[\big[\big[ g\big]\big]\big]_{r+s+1-j,\infty}\|\kappa W\|_{0,j}\tag 6.22
$$
(6.4) follows from this and $\pa_a q=g_{ab} W^b$.
(Replacing $\kappa$ by $1$ just causes more terms of the same form.)

To prove (6.5) we need estimates for the last terms on the right in (6.4).
By (6.4) with $r=0$ we have
$$
\|\pa^2 q\|_{0,s}\leq K_{10}\sum_{j=0}^s
\big[\big[\big[ g\big]\big]\big]_{s-j,\infty}\|\triangle q\|_{0,j}
+K_{10}\sum_{j=0}^s \big[\big[\big[ g\big]\big]\big]_{s+1-j,\infty}\|\pa q\|_{0,j}
.\tag 6.23
$$
Since $\hat{D}_t^s \triangle q=\kappa^{-1}\pa_a D_t^s \big(\kappa g^{ab}\pa_b
q\big)$ we have
$$
\|\triangle D_t^s q\|\leq \| \hat{D}_t^s \triangle q\|+ K_{10}
\sum_{j=0}^{s-1} \big( \big[\big[\big[ g\big]\big]\big]_{s-j,\infty} \|\pa^2 D_t^j q\|
+\big[\big[\big[
g\big]\big]\big]_{s+1-j,\infty}\|\pa D_t^j q\|\big),\tag 6.24
$$
and using (6.23) it follows that
$$
\|\triangle D_t^s q\|
\leq K_{10}\sum_{j=0}^s
\big[\big[\big[ g\big]\big]\big]_{s-j,\infty}\|\triangle q\|_{0,j}
+K_{10}\sum_{j=0}^{s-1} \big[\big[\big[ g\big]\big]\big]_{s+1-j,\infty}\|\pa q\|_{0,j}
.\tag 6.25
$$
We have $\int_\Omega g^{ab}(\pa_a q)(\pa_b q)\, \kappa \, dy
=-\int_{\Omega}( \triangle q) q\, \kappa \, dy$ and there is
a constant $K_{10}$, depending just on the volume of $\Omega$, i.e.
$\int_{\Omega} \kappa \, dy$, such that $\|q\|\leq K_{10}\|\triangle
q\|$, see \cite{SY}. Hence in addition we have
$$
\|q\|+\|\pa q\|\leq K_{10}\|\triangle q\|,
\qquad\quad \text{if}\qquad q\big|_{\pa\Omega}=0
 \tag 6.26
$$
(6.26) applied to $D_t^s q$ in place of $q$ together with (6.25) gives
$$
\|\pa q\|_{0,s}
\leq K_{10}\sum_{j=0}^s
\big[\big[\big[ g\big]\big]\big]_{s-j,\infty}\|\triangle q\|_{0,j}
+K_{10}\sum_{j=0}^{s-1} \big[\big[\big[ g\big]\big]\big]_{s+1-j,\infty}\|\pa q\|_{0,j}
,\tag 6.27
$$
where the last term should be interpreted as $0$ if $s=0$,
and inductively it follows that
$$
\|\pa q\|_{0,s}
\leq K_{20}\sum_{j=0}^s
\big[\big[\big[(\pa g,D_t g)\big]\big]\big]_{s-j,\infty}\|\triangle q\|_{0,j}
\tag 6.28
$$

\comment
and hence if we also use (6.23) we obtain
$$
\|q\|_{0,k}+\|\pa q\|_{0,k}\leq K_{10} \|\hat{D}_t^k \triangle q\| +K_{10}
\sum_{j=0}^{k-1} \big( \big[\big[\big[g\big]\big]\big]_{k-j,\infty} \|\pa^2 q\|_{0,j}
+\big[\big[\big[g\big]\big]\big]_{k+1-j,\infty} \|\pa q\|_{0,j}\big) \tag 6.25
$$
and
$$
\|\pa^2 q\|_{0,k}\leq K_{10}\|g\|_{1,\infty}\|\hat{D}_t^k \triangle q\|
+K_{10} \sum_{j=0}^{k-1} \big( \big[\big[\big[g\big]\big]\big]_{k+1-j,\infty} \|\pa^2
q\|_{0,j} +\big[\big[\big[g\big]\big]\big]_{k+2-j,\infty} \|\pa q\|_{0,j}\big)
 \tag 6.26
$$
Inductively it follows that
$$
\| q\|_{0,s}+\|\pa q\|_{0,s}\leq K_1\sum_{j=0}^s
\big[\big[\big[g\big]\big]\big]_{s-j,\infty}\|\triangle q\|_j\tag 6.27
$$
Since $W^a=g^{ab} \pa_b q$ we also get
$$
\|W\|_{0,s}\leq K_1\sum_{j=0}^s \big[\big[\big[g\big]\big]\big]_{s-j,\infty} \|\pa
q\|_j\leq K_1\sum_{j=0}^s \big[\big[\big[g\big]\big]\big]_{s-j,\infty}\|\triangle
q\|_j\tag 6.28
$$
where we also used interpolation. (6.28) together with (6.22)
proves (6.4) for $W$. The estimate for $\pa q$ follows from the
same argument since $\pa_a q=g_{ab}W^b$.
\endcomment

It remains to prove the estimates for the projection (6.6). We
have $W=W_0+W_1$, where $W_0=PW$, and $W_1=g^{ab}\pa_b q$ where
$\triangle q=\div W$ and $q\big|_{\pa\Omega}=0$. Proving
(6.6) for $r\geq 1$ reduces to proving it for $r=0$
by using (6.4), since $\hat{R}^I \triangle
q=\div\,( \hat{\Cal L}_R^I W)$ and replacing $\kappa$ by $1$ just produces
more terms of the same form . (6.6), for $r=0$ and $s=0$ follows since
the projection has norm $1$, $\|PW\|\leq \|W\|$.  Since the
projection of $g^{ab}D_t^k w_{1b}=g^{ab}\pa_b D_t^k q$ vanishes we obtain from Lemma
{6.}3:
$$
\| P \hat{D}_t^k W_1\| \leq K_1 \sum_{i=0}^{k-1} \big[\big[\big[g\big]\big]\big]_{k-i}
\|\hat{D}_t^i W_1\|\tag 6.29
$$
Since also $P \hat{D}_t^k W_0=\hat{D}_t^k W_0$ we have
$$
(I-P)\hat{D}_t^k W_1=(I-P) D_t^k W\tag 6.30
$$
and hence since the projection has norm one
$$
\| \hat{D}_t^k W_1\|+\| \hat{D}_t^k W_0\|\leq K_1 \|\hat{D}_t^k W\|+ K_1
\sum_{i=0}^{k-1} \big[\big[\big[g\big]\big]\big]_{k-i} \|\hat{D}_t^i W_1\|\tag 6.31
$$
Using induction we therefore obtain
$$
\|\hat{D}_t^k W_0\|+\|\hat{D}_t^k W_1\|\leq K_1 \sum_{j=0}^k
\big[\big[\big[g\big]\big]\big]_{k-j,\infty} \|\hat{D}_t^{j} W\|\tag 6.33
$$
Since as before replacing $\kappa$ by $1$ just produces
more terms of the same form
this proves (6.6) also for $r=0$.
(6.7)-(6.8) follows from interpolation.\qed\enddemo

\head 7. Tame estimates for the wave equation.\endhead
 Existence of solutions for a wave equation with Dirichlet
boundary conditions and initial conditions satisfying some
compatibility conditions is well known, see e.g. \cite{H3}. The
result in \cite{H3} is stated with vanishing initial condition but
we will reduce to that case by subtracting of a function that
solves the equation to all orders as $t\to 0$. We consider the
Cauchy problem for the wave equation on a bounded domain with
Dirichlet boundary conditions:
$$\align
{D}_t^2 (e^\prime\psi)-\triangle \psi &={f},\qquad
\text{in}\quad [0,T]\times\Omega,\qquad \psi\big|_{\pa\Omega}=0,\tag 7.1 \\
\psi\big|_{t=0}\!\! &=\psi_0,\quad {D}_t\psi\big|_{t=0}\!\!=\psi_1
\tag 7.2
\endalign
$$
Here
$$
\triangle\psi =\frac{1}{\sqrt{\det{g}}}\pa_a \Big( \sqrt{\det{g}}
g^{ab}\pa_b \psi\Big),\tag 7.3
$$
where $g^{ab}$ is the inverse of the metric $g_{ab}$ and
$\det{g}=\det\{g_{ab}\}=\kappa^2$, in our earlier notation. We
assumed that  $e^\prime$ is positive, $g^{ab}$ is symmetric and positive
definite (since the metric is), and that $g^{ab}$ and $e^\prime$ are
smooth satisfying:
$$
e^\prime+1/e^\prime\leq c_1^\prime,\qquad \sum_{a,b}(|g^{ab}|+|g_{ab}|)\leq
c_1^2, \tag 7.4
$$
for some constants $0<c_1<c_1^\prime<\infty$.

\comment  With vanishing initial conditions and a smooth
inhomogeneous term that vanishes to all orders as $t\to 0$, the
equation (7.1) has a smooth solution by \cite{H1}. If initial
conditions vanish, the inhomogeneous term is smooth but only
vanishes to order $m-2$ as $t\to 0$ then we can replace the right
hand side by ${f}_\varepsilon =\chi_\varepsilon {f}$, where
$\chi_\varepsilon(t)=\chi(\varepsilon t)$, and  $\chi\in C^\infty$
is $0$ for $t\leq 0$ and $1$ for $t\geq 1$. It then follows that
the solution $\psi_\varepsilon $ will converge to a solution of
(7.1) and satisfy an estimate
 $ \|\psi(t,\cdot)\|_{k}\leq C\int_0^t \|\dot{{f}}\|_{k-2}\, d\tau$
  for $k\leq m+1$,
  where $\dot{{f}}=D_t {f}$.
Finally, if the inhomogeneous term is not smooth we can regularize
it and pass to the limit to obtain a solution in any of the spaces
defined by the above estimates.
\endcomment

 We will apply this theorem to
$$
D_t^2 (e^\prime(h)\delta h)-\triangle\delta h=-\pa_i \big( ( \pa^i
\delta x^k)\pa_k h\big)-(\pa^i\delta x^k)\pa_i \pa_k h +2(\pa_i
V^k)\pa_k\delta V^i ,\qquad \delta h\Big|_{\pa \Omega}\!\!\!=0\tag
7.5
$$
with vanishing initial conditions and the right hand side vanishing
to all orders as $t\to0$. We will also need some estimates for
the equation for the enthalpy itself, which is best dealt with by
writing $h=\tilde{h}+h_0$ where $h_0$ is the smooth approximate
solution. In particular, in case $e(h)=c^{-2} h$ we get the equation
$$
c^{-2} D_t^2 \tilde{h}-\triangle \tilde{h}= -\big(c^{-2}D_t^2
h_0-\triangle h_0- (\pa_i V^j)(\pa_j V^i)\big)\tag 7.6
$$
with vanishing initial conditions but where the right hand side
vanishes to all orders as $t\to 0$ and is smooth for $t\geq 0$. In
both cases we have vanishing initial conditions.

 In the first case, (7.5), we then also want to use the special form of
${f}$ in (7.1).
$$
{f}={f}_1+\kappa^{-1} \pa_a (\kappa F_1^a )\tag 7.7
$$
For the lowest order energy estimate it will be useful to take
advantage of the special form (7.7) because it will allow us to
estimate $F$ with one less space derivative and one more time
derivative instead and in the estimate for the divergence equation
we have one more time derivative than space derivatives.
 Let
$$
g_s=\|| g\||_{s+1,\infty}, \qquad h_s=\||h\||_{s+2,\infty}, \tag
7.8
$$
and let $K_2^\prime$ denote a continuous function of
 $$
h_0+g_0+c_1^\prime+T^{-1}+r\tag 7.9
$$
which in what follows also depends on the order $r$ of
differentiation.
 We will use that by interpolation;
$$
(g_{s}+h_{s})(g_{r}+h_{r})\leq K_2^\prime (g_{r+s}+h_{r+s}).\tag
7.10
 $$
\proclaim{Theorem {7.}1} Suppose that initial data in (7.2)
vanishes and ${f}$ in (7.1) or ${f}_1$ and $F_1$ in (7.7) are smooth
and  vanish to all orders as $t\to 0$. Suppose also that $e^\prime=e^\prime(h)$
in (7.1) is a smooth function of $h$ and that $h$ and $g$ are
smooth. Then for $r\geq 1$ the solution of (7.1) satisfies the
estimates
$$
\|\psi\|_{0,r}+\|\psi\|_{1,r-1}\leq K_2^\prime\sum_{s=1}^r
(h_{r-s}+g_{r-s}) \int_0^t (\|{f}_1\|_{0,s-1}+\|F_1\|_{0,s})\,
d\tau\tag 7.11
$$
Furthermore
$$
\|\psi\|_{r}\leq K_2^\prime\sum_{s=1}^r (h_{r-s}+g_{r-s})\big(
\int_0^t (\|{f}_1\|_{0,s-1}+\|F_1\|_{0,s})\, d\tau+\|{f}\|_{s-2}
\big)\tag 7.12
$$
and
$$
\|\psi\|_{r}\leq K_2^\prime\sum_{s=1}^r (h_{r-s}+g_{r-s})\int_0^t
(\|\dot{{f}}\|_{s-2} +\|{f}\|_{0})\, d\tau\tag 7.13
$$
 where $\|{f}\|_{s-2}$ should be interpreted as $0$
if $s=1$.
\endproclaim
\demo{Proof} Let us first prove (7.11).  Let $E_r$ be as in Lemma
{7.}2 and let
$$
\tilde{E}_r^2=E_r^2+\int_{\Omega} \psi^2\, \kappa dy
+\sum_{s=0}^{r-1} \int_{\Omega} g_{ab}(\hat{D}_t^s
F_1^a)(\hat{D}_t^s F_1^b)\, \kappa dy\tag 7.14
$$
Then it follows that $\tilde{E}_r$ satisfy the same inequality as
$E_r$ in Lemma {7.}2 does, even without the last term in (7.19) since the norm
of $\psi$ is included in $\tilde{E}_1$:
$$
\frac{d\tilde{E}_r}{dt}\leq  K_2^\prime\sum_{s=1}^{r}
(h_{r-s}+g_{r-s}) \big(\tilde{E}_s+\|F_1\|_{0,s}+\|{f}_1\|_{0,s-1}
\big)\tag 7.15
$$
The highest order energy in the right hand side $\tilde{E}_s$ with
$s=r$ can be removed by multiplying with the integrating factor
$e^{K t}$ and integrating. We get
$$
\tilde{E}_r\leq K_2^\prime\sum_1^{r} (h_{r-s}+g_{r-s}) \int_0^t
(\|{f}_1\|_{0,s-1}+\|F_1\|_{0,s})\, d\tau
+K_2^\prime\sum_{s=1}^{r-1} (h_{r-s}+g_{r-s})\int_0^t
\tilde{E}_s\, d\tau\tag 7.16
$$
We claim that if $r\geq 1$
$$
\tilde{E}_r\leq K_2^\prime\sum_{s=1}^{r} (h_{r-s}+g_{r-s})
\int_0^t (\|{f}_1\|_{0,s-1}+\|F_1\|_{0,s})\, d\tau\tag 7.17
$$
If $r=1$, we have just proven it in (7.16) and in general it
follows from (7.16) using induction an interpolation.

We claim that (7.11) follows from (7.17). In fact
$\|\hat{D}_t^{r-1} (g^{ab}\pa_b\psi) \|\leq C \tilde{E}_r$ and
 by Lemma {6.}3
 $\hat{D}_t^{r-1} (g^{ab}\pa_b \psi)=
 g^{ab}\pa_b D_t^{r-1} \psi+ \sum_{s=0}^{r-2}
 \binom{r-1}{s}
g^{ab}(\check{D}_t^{r-1-s} g_{bc})
 \hat{D}_t^{s} (g^{cb}\pa_b\psi)$. If follows that
 $\|\pa D_t^{r-1} \psi\|\leq \sum_{i=1}^r  g_{r-s} \tilde{E}_s$
 which together with (7.16) and interpolation prove (7.11).

 Finally, (7.12) follows from (7.11) and Lemma {7.}3.
 and (7.13) follows from (7.12) using that
 $\|{f}\|_{s-2}\leq \int_0^t \|\dot{{f}}\|_{s-2}\, d\tau$.
\enddemo

\proclaim{Lemma {7.}2}
 Suppose that $g^{ab}$ and $e^\prime=e^\prime(h)$ are smooth and satisfy
(7.4) for $r\geq 1$,
$$
E_r(t)=\Big(\sum_{s=0}^{r-1} \frac{1}{2}\!\int_\Omega{ \!\! e^\prime
(\!{D}_t^{s+1}\psi)^2\! + g_{ab} (\hat{D}_t^s \Psi^a )
(\hat{D}_t^s \Psi^b )\kappa d y}\Big)^{1/2} , \qquad \Psi^a=
g^{ab}\pa_b \psi+F_1^a\tag 7.18
$$
Then for $r\geq 1$
$$
\frac{d E_r}{dt} \leq K_2^\prime\sum_{s=1}^r (h_{r-s}+g_{r-s}
)\big(E_s +\| F_1\|_{0,s}+\|{f}_1\|_{0,s-1}\big)+K (g_{r-1}+h_{r-1}) \|\psi\|_0
\tag 7.19
$$
\endproclaim
\demo{Proof} We will prove that $ d E_r^2/dt$ is bounded by $E_r$
times the right hand side of (7.19), and (7.19) follows from this since $d
E_r/dt=(d E_r^2/dt)/(2E_r)$. With the notation $\hat{D}_t=D_t+\div
V$ and $\check{D}_t=D_t-\div V$ we have $\hat{D}_t (e^\prime g^2)
=(\check{D}_t e^\prime)g^2+ 2e^\prime g\hat{D}_t g$, for any functions $e^\prime$ and
$g$, and since also $D_t\kappa/\kappa=\div V$ it follows that
$$
\multline \frac{d E^2_r}{dt}=\sum_{s=0}^{r-1} \int_\Omega{ \Big( e^\prime
({D}_t^{s+1}\psi)({D}_t^{s+2}\psi)+ g_{ab}(\hat{D}_t^s \Psi^a )
(\hat{D}_t^{s+1} \Psi^b )
\Big)\,  \kappa \,d y}\\
+\frac{1}{2}\sum_{s=0}^{r-1} \int_\Omega{ \Big((\hat{D}_t e^\prime)
({D}_t^{s+1}\psi)^2+ (\check{D}_t g_{ab})(\hat{D}_t^s \Psi^a )
(\hat{D}_t^s \Psi^b ) \Big)\, \kappa \,d y}.
\endmultline\tag 7.20
$$
Here the terms on the second row are bounded by $K_2^\prime
E_r^2$. Therefore it remains to look on the terms on the first
row.
 Since $\hat{D}_t (e^\prime)=\kappa^{-1}D_t(\kappa e^\prime)$ it follows that
 $$
 e^\prime D_t^{s+2}\psi= \hat{D}_t^s D_t^2(e^\prime\psi)+ \sum_{i=0 }^{s+1} B^{s}_i D_t^i \psi,
 \qquad B_i^s=\!\!\!\!\!\sum_{j=\max(0,i-2)}^{s}\!\!\!\!\!\
  \kappa^{-1}c_{ij}^s (D_t^{s-j}\kappa)(D_t^{j+2-i} e^\prime)\tag 7.21
$$
Furthermore, using Lemma {6.}3 we get, since $\underline{\Psi}_a=
\pa_a \psi+\underline{F}_a$,
 $$
 \hat{D}_t^{s+1}\Psi^a =g^{ab}\pa_a
 D_t^{s+1}\psi+g^{ab}D_t^{s+1}\underline{F}_b-
 \sum_{i=0}^{s}\binom{s+1}{i} g^{ab} (\check{D}_t^{s+1-i} g_{bc}) \hat{D}_t^i
 \Psi^c\tag 7.22
 $$
 Since $e^\prime=e^\prime(h)$ it follows that the $L^2$ norm of the sums in (7.21)
and (7.22) are
 bounded by $K_2^\prime\sum_{s=0}^r (g_{r-s}+h_{r-s}) E_s$, which is
 included in the right hand side of (7.19) and so is
 $\|g^{ab}D_t^{s+1}\underline{F}_b\|$.
 Therefore it remains to consider
$$
\multline \sum_{s=0}^{r-1} \int_\Omega{ \!\Big(
({D}_t^{s+1}\psi)(\hat{D}_t^{s}D_t^2(e^\prime\psi)+
\big(\hat{D}_t^s\Psi^a \big) \big(\pa_a {D}_t^{s+1}\psi )\big)
\Big)\,  \kappa \,d y}\\
=\sum_{s=0}^{r-1}\int_\Omega{ \Big( ({D}_t^{s+1}\psi )
\big(\hat{D}_t^{s}D_t^2(e^\prime\psi)-\kappa^{-1} \pa_a\big(\kappa
\hat{D}_t^s \Psi^a \big)\Big)\, \kappa
dy}=\sum_{s=0}^{r-1}\int_\Omega{ \Big( ({D}_t^{s+1}\psi )
\big(\hat{D}_t^{s} {f}_1)  \kappa dy}
\endmultline\tag 7.23
$$
Here we have integrated by parts using that
$D_t^{s+1}\psi\big|_{\pa\Omega}=0$, that $D_t^2
(e^\prime\psi)-\kappa^{-1}\pa_a\big(\kappa \Psi^a\big)={f_1}$ and that
$\hat{D}_t^s \big(\kappa^{-1}\pa_a\big(\kappa
\Psi^a\big)\big)=\kappa^{-1}\pa_a\big(\kappa
\hat{D}_t^s\Psi^a\big)$.
 \qed\enddemo

One can get additional space regularity from taking time
derivatives of the equation (7.1) and solving the Dirichlet
problem for the Laplacian.

\proclaim{Lemma {7.3}} If $r\geq 1$ then
$$
 \|\psi\|_r\leq K_2^\prime\sum_{i=1}^r (g_{r-i}+h_{r-i})
 (\|\psi\|_{0,i}+\|\psi\|_{1,i-1}+\|{f}\|_{i-2})
\tag 7.24
$$
where $\|{f}\|_{i-2}$ is to be interpreted as $0$ if $i=1$.
\endproclaim
\demo{Proof} We can use the estimates for the Dirichlet problem:
$$
\triangle \psi={D}_t^2 (e^\prime\psi)-{f}, \qquad
\psi\big|_{\pa\Omega}=0\tag 7.25
$$
from Theorem {6.}1, since also $e^\prime=e^\prime(h)$.
$$
\|\psi\|_{{r}-s,s}\leq K_2^\prime\sum_{i\leq r-s-2,\, j\leq s+2}
(g_{r-i-j}+h_{r-i-j}) \big( \|\psi\|_{i,j} +\|{f}\|_{i,j-2}\big)\tag
7.26
$$
where $\|{f}\|_{i,j-2}$ is to be interpreted as $0$ if $j-2<0$.
Hence using induction and interpolation we get
 $$
 \|\psi\|_r\leq K_2^\prime\sum_{i=0}^r (g_{r-i}+h_{r-i})
 (\|\psi\|_{0,i}+\|\psi\|_{1,i-1}+\|{f}\|_{i-2})\tag 7.27
 $$
 where $\|\psi\|_{1,i-1}$ is to be interpreted as $0$ if $i=0$
 and $\|{f}\|_{i-2}$ is to be interpreted as $0$ if $i-2<0$.
\enddemo

\proclaim{Theorem {7.}4} Suppose that $\phi$ is a solution of
$$
D_t^2 \phi-\triangle (p^\prime \phi)=f,
\qquad \phi\big|_{t=0}=\dot{\phi}\big|_{t=0}=0.\tag 7.28
$$
where $p^\prime=p^\prime(h)=1/e^\prime(h)$ and $f$ vanishes to all orders as $t\to 0$;
$D_t^k f\big|_{t=0}=0$, for $k\geq 0$. Let
$$
W_1=\nabla q, \qquad \triangle q=\phi,\qquad q\big|_{\pa\Omega}=0\tag 7.29
$$
and
$$
F_1=\nabla q^\prime,\qquad \triangle q^\prime=f,\
\qquad q^\prime\big|_{\pa\Omega}=0\tag 7.30
$$
Then for $r\geq 1$
$$
\|W\|_{r+1}\leq K_2^\prime\sum_{s=1}^r
(h_{r-s}+g_{r-s}) \int_0^t (\|\dot{F}_1\|_{s-1}+\|F_1\|_0)\, d\tau
\tag 7.31
$$
\endproclaim
\demo{Proof} By Theorem {7.}1
$$
\|\phi \|_{r}\leq K_2^\prime\sum_{s=1}^r
(h_{r-s}+g_{r-s}) \int_0^t (\|\dot{F}_1\|_{s-1}+\|F_1\|_0)\, d\tau\tag 7.32
$$
which by inverting the Laplacian in (7.29)
proves that $\|\pa W_1\|_r$ is bounded by the right hand side of
(7.31) and it only remains to prove the estimate for $\|D_t^{r+1} W_1\|$.
Using (2.54) we can write (7.28) as
$$
\hat{D}_t^2 \phi -2\dot{\sigma} \hat{D}_t \phi
+k^{\prime\prime}\phi -\triangle(p^\prime \phi)=f,\qquad k^{\prime\prime}
=\dot{\sigma}^2-\ddot{\sigma}\tag 7.33
$$
and so
$$
\div \big( \ddot{W}_1-\nabla \big(p^\prime \div W_1\big)-2\dot{\sigma}\dot{W}_1
+k^{\prime\prime}  W_1 -F_1\big)=
 -2(\pa_a\dot{\sigma}) \dot{W}_1^a
+(\pa_a k^{\prime\prime}) W_1^a\tag 7.34
$$
and hence
$$
\div\big( \hat{D}_t^j \ddot{W}_1 - \hat{D}_t^j
\nabla \big(p^\prime \div W_1\big)- \hat{D}_t^j \big(
2 \dot{\sigma}\dot{W}_1-k^\prime W_1+F_1\big)\big)=
-\hat{D}_t^j \big( 2(\pa_a\dot{\sigma}) \dot{W}_1^a
-(\pa_a k^{\prime\prime}) W_1^a\big)\tag 7.35
$$
We claim that
$$
P \hat{D}_t^j W_1=P B_j(W_1,...,\hat{D}_t^{j-1} W_1),\qquad
B_j(W_1,...,\hat{D}_t^{j-1} W_1)=-\sum_{i=0}^{j-1} {\tsize{\binom{j}{i}}}
(\check{D}_t^{j-i} g_{ab})\hat{D}_t^{i} W_1 \tag 7.36
$$
In fact $0=P D_t^j \pa_a q= P D_t^j\big( g_{ab} W_1^b\big)
=P \sum_{i=0}^j \binom{j}{i}
(\check{D}_t^{j-i} g_{ab})\hat{D}_t^{i} W_1$. Furthermore, let
$$
\triangle q^j=-\hat{D}_t^j \big( 2(\pa_a\dot{\sigma}) \dot{W}_1^a
-(\pa_a k^{\prime\prime}) W_1^a\big),\qquad
q^j\big|_{\pa\Omega}=0\tag 7.37
$$
To say that $\div H=0$ is equivalent to saying that $(I-P) H=0$ so
it follows from (7.35)-(7.37) that
$$
\hat{D}_t^{j} W_1=PB_j(W_1,...,\hat{D}_t^{j-1} W_1)
+(I-P) \hat{D}_t^{j-2}\big( -
\nabla \big(p^\prime \div W_1\big)
+2 \dot{\sigma}\dot{W}_1-k^{\prime\prime} W_1+F_1\big) +\nabla q^{j-2}\tag 7.38
$$
Here
$$
\hat{D}_t^{j-2}\big(\nabla^a \big(p^\prime \div W_1\big)\big)=
 \sum_{i=0}^{j-2} \sum_{m=0}^i
{\tsize{\binom{j-2}{i}}}{\tsize{\binom{i}{m}}}
(\hat{D}_t^{j-2-i} g^{ab})
\pa_b \big( (D_t^{i-m}p^\prime(h))( D_t^m \div W_1)\big)\tag 7.39
$$
and by interpolation
$$
\|\hat{D}_t^{j-2}\big(\nabla^a \big(p^\prime \div W_1\big)\big)\|
\leq  K_2^{\prime}
\sum_{s=0}^{j-2+1} (g_{j-2-s}+h_{j-3-s})\|\div W_1\|_s\tag 7.40
$$
Since the projection $P$ maps $L^2$ to $L^2$ and since inverting (7.37)
maps $L^2$ to $H^1$ (in fact to $H^2$) it therefore follows that
$$
\|D_t^{r+1} W_1\|\leq K_2^{\prime}
\sum_{s=0}^{r}\big(  g_{r-s}\|D_t^s W_1\|+
 (g_{r-1-s}+h_{r-2-s})\|\div W_1\|_s\big)
+ K_2^{\prime}\sum_{s=0}^{r-1} g_{r-1-s} \|D_t^s F_1\|
\tag 7.41
$$
(7.31) follows from this if we use (7.32) to estimate $\phi=\div W_1$
and use that $ \|D_t^s F_1\|\leq \int_0^t \|D_t^s \dot{F}_1\|\, d\tau$.
\qed\enddemo

\proclaim{Corollary {7.}5} With assumptions as in Theorem {7.}1 we
have
$$
\||\psi\||_{r,\infty}\leq K_2^\prime\big(
(g_{r+r_0-1}+h_{r+r_0-1})\|| {f}\||_{0,\infty}
+\||{f}\||_{r+r_0-1,\infty}\big)\tag 7.42
$$
where $r_0=[n/2]+1$.
\endproclaim
\demo{Proof} By Sobolev's lemma $\|\psi\|_{r,\infty}\leq C
\|\psi\|_{r+r_0}$, which in turn can be estimated by Theorem
{7.}1. Then we estimate the $L^2$ norms of ${f}$ by $L^\infty$ norms
and use interpolation. \qed\enddemo

Let us now prove that the solution of (7.1) depends smoothly on
parameters if the metric $g$ and the inhomogeneous term ${f}$ do. We
have: \proclaim{Lemma {7.}6} Suppose that ${f}\in
C^m\big(B^k,C^\infty_{00}([0,T]\times\overline{\Omega})\big)$,
$g^{ab}\in
C^m\big(B^k,C^\infty([0,T]\times\overline{\Omega})\big)$, where
$B^k=\{r\in \bold{R}^k,\; |r|\leq \varepsilon\}$. Suppose also that $g$
satisfies the coordinate condition uniformly in $B^k$. Let $\psi$
be the solution of
$$
\square \psi =c^2 D_t^2 \psi -\triangle \psi={f},  \qquad
\psi\big|_{\pa\Omega}=0 \qquad
\psi\big|_{t=0}=D_t\psi\big|_{t=0}=0\tag 7.43
$$
where $\triangle$ is given by (7.3). Then $\psi\in
C^m\big(B^k,C^\infty_{00}([0,T]\times\overline{\Omega})\big)$.
\endproclaim
\demo{Proof}   We note that it suffices to prove the statement in
the theorem for $m=1$, since the general case follows from this by
induction. In fact, assuming that $\psi\in
C^l\big(B^k,C^\infty_{00}([0,T]\times\overline{\Omega})\big)$,
$l<m$ then for $|\alpha|\leq l$:
$$
\square D_r^\alpha \psi=D_r^\alpha {f} -\sum_{\beta+\gamma=\alpha,\,
|\beta|<|\alpha|}c^\alpha_{\beta }\square^{(\gamma)} D_r^\beta
\psi_r, \qquad
 \tag 7.44
$$
where $\gamma=(\gamma_j,...\gamma_1)$ are ordered multi-indices
with  $\gamma_i\in \{1,...,k\}$. Here
$D_r^\gamma=D_{r^{\gamma_j}}\cdot\cdot\cdot D_{r^{\gamma_1}}$,
where $D_{r^i}=\pa/\pa r^i$, and $\square_r^{(\gamma)}$ are the
repeated commutators defined inductively by
$\square^{(i,\gamma)}=[D_{r^{i}},\square^{(\gamma)}]$, where
$\square^{(\,\,)}= \square_r$. The right hand side of (7.44) is
then in
$C^1\big(B^k,C^\infty_{00}([0,T]\times\overline{\Omega})\big)$ so
it follows from the statement in the theorem for $m=1$ that
 $\psi\in
C^{l+1}\big(B^k,C^\infty_{00}([0,T]\times\overline{\Omega})\big)$.

Let us write $\psi_r$, $g_r$ ${f}_r$, and $\square_r=\square_{g_r}$,
to indicate the dependence of $r$. First we will prove that
${f}_r\in C(B^k,C_{00}^\infty)$ and $g_r\in C(B^k,C^\infty)$ implies
that $\psi_r \in C(B^k,C_{00}^\infty)$. We will only prove this
for $r=0$ since the proof in general is just a notational
difference from the proof for $r=0$. We have
$$
\square_r (\psi_r-\psi_0)={f}_r-{f}_0-(\square_r-\square_0)\psi_0\tag
7.45
$$
The right hand side is in $C(B^k,C_{00}^\infty)$ and tends to $0$
in $C(B^k,C^N)$, for any $N$, as $r\to 0$. Since in Corollary
{7.}5 we have uniform bounds for $\square_r^{-1}$, it follows that
$\psi_r-\psi_0$ tends to $0$ in $C(B^k,C^N)$, for any $N$, as
$r\to 0$, i.e. $\psi_r\in C(B^k,C_{00}^\infty)$.

Let us now assume that ${f}_r\in C^1(B^k,C_{00}^\infty)$ and $g_r\in
C^1(B^k,C^\infty)$. Let $\dot{\psi}_r=D_r\psi_r$ be defined by
$$
\square_r \dot{\psi}_r= \dot{{f}}_r-\dot{\square}_r \psi_r,
\qquad\quad  \dot{\psi}_r\big|_{\pa{\Omega}}=0,\qquad
\dot{\psi}_r\big|_{t=0}=D_t\dot{\psi}_r\big|_{t=0}=0,\tag 7.46
$$
where $\dot{\square}_r=[D_r, \square ]$, $\dot{{f}}_r=D_r {f}_r$ and
$D_r=(D_{r^1},...,D_{r^k})$. Since the right hand side of (7.46)
is in $ C(B^k,C_{00}^\infty)$ it follows as above that
$\dot{\psi}_r\in  C(B^k,C_{00}^\infty)$.
 It remains to prove that
$\psi_r$ is differentiable. We have
$$
\square_r \big(\psi_r-\psi_0-r\dot{\psi}_0\big) ={f}_r-{f}_0-r
\dot{{f}}_0 +\big(\square_r-\square_0-r\dot{\square}_0\big)\psi_0
+r(\square_r-\square_0)\dot{\psi}_0\tag 7.47
$$
Then the right hand side divided by $r$ tends to $0$ in $C^N$, for
any $N$, as $r\to 0$. In view of the uniform bounds for
$\square^{-1}_r$ in Corollary {7.}5 it follows that
$(\psi_r-\psi_0-r\dot{\psi}_0)/r$,  tends to $0$ in $C^N$ for any
$N$, as $r\to 0$. Hence $\psi\in C^1(B^k, C_{00}^\infty)$.
 \qed\enddemo

\head{8. Tame estimates for the divergence free equation.}\endhead
 Let
$$
\ddot{W}+AW+B_0 W+B_1\dot{W}=H,\qquad \div H=0,\qquad
W\big|_{t=0}=\dot{W}\big|_{t=0}=0\tag 8.1
$$
where $H$ is smooth and vanishes to order $r$ as $t\to 0$, i.e.
$D_t^k H\big|_{t=0}=0$, for $k\leq r$. Here $B_i W^a=P\big( g^{ab}
\beta_{bc}^i W^c\big)$, $i=0,1$, are projected multiplication
operators. Let
$$
g_s=\|| g\||_{s+1,\infty}, \qquad h_s=\||h\||_{s+2,\infty},\qquad
\beta_s=\||\beta^0\||_{s,\infty}\!+\||\beta^1\||_{s,\infty} ,
\qquad k_s=g_s+h_s+\beta_s \tag 8.2
$$
and let $K_2^{\prime\prime}$
 denote a continuous function of
 $$
k_0+c_0^{-1}+c_1+T^{-1}+r,  \tag 8.3
$$
which in what follows also depends on the order $r$ of
differentiation.
 We will use that by interpolation;
$ k_s k_r\leq K_2^{\prime\prime} k_{r+s}
 $, for $r,s\geq 0$.
  We will prove the following estimates

\proclaim{Theorem {8.}1} Suppose that (2.4)-(2.5) hold for $0\leq
t\leq T$ and suppose also that $x$ is smooth for $0\leq t\leq T$,
$T\leq 1$.  Then (8.1) has a smooth solution for $0\leq t\leq T$.
It satisfies the estimates
$$
\|\dot{W}\|_r+\|W\|_r\leq K_2^{\prime\prime} \sum_{s=0}^r
k_{r-s}\int_0^t\|H\|_s\, d\tau,\tag 8.4
$$
for $r\geq 0$ and for $m=0,1,2$.
$$
\sum_{j=0}^{m+1}\|\hat{D}_t^j W\|_r\leq K_2^{\prime\prime}
\sum_{s=0}^r\sum_{j=0}^m k_{r+m-j-s}\int_0^t\|\hat{D}_t^ j H\|_s\,
d\tau,\tag 8.5
$$
Furthermore, for $r\geq 1$,
$$\multline
\|\dddot{W}\|_{r-1}+\|\ddot{W}\|_{r-1}+\|\dot{W}\|_r+\|W\|_r\\
\leq K_2^{\prime\prime}\sum_{s=0}^{r-1}
\int_0^t\big(k_{r-1-s}\|\ddot{H}\|_s+k_{r-s}\|\dot{H}\|_s
+k_{r+1-s}\|H\|_s+\|\curl \underline{H}\|_s\big)\, d\tau .
\endmultline \tag 8.6
$$
\endproclaim

As pointed out before, existence of smooth solutions for (8.1) was
proven in \cite{L3} so we only need to prove the estimates. The
theorem with the norms replaced by norms with just space
derivatives was proven in \cite{L3}. The theorem here is actually
simpler to prove than the one there.

First we note that we can reduce to the case $B_0=B_1=0$ since the
general case follows from this case. Let us show this for (8.4).
Assume that (8.4) holds for the case when $B_0=B_1=0$. Then using
(8.4) applied to the equation $\ddot{W}+AW=H_1=H+B_0 W+B_1 \dot{W}
$ gives using that, by Theorem {6.}1, $\| B_i W\|_s\leq
K_2^{\prime\prime}\sum_{k=0}^s (\beta_{s-k}+g_{s-k})\|W\|_k$ and
(8.4),
$$
\|\dot{W}\|_r+\|W\|_r\leq K_2^{\prime\prime} \sum_{s=0}^r
(g_{r-s}+h_{r-s}+\beta_{r-s} )\int_0^t(\|H\|_s+
\|W\|_s+\|\dot{W}\|_{s}) , d\tau \tag 8.7
$$
We can now first remove the highest order terms from the right
hand side by a Gr\"onwall type of argument since they occur in the
left. In fact let $f(t)=\int_0^t \|\dot{W}\|_r+\|W\|_r\, d\tau$.
Then by (8.7) $f^\prime\leq K_2^{\prime\prime} f
+K_2^{\prime\prime}\|H\|_r+\Sigma_{r-1}$, where $\Sigma_{r-1}$ is
the sum of the first $r-1$ terms in the right of (8.7). Hence
multiplying by the integrating factor $e^{K_2^{\prime\prime} t}$
and integrating gives $f(t)\leq K_2^{\prime\prime}\int_0^t
\|H\|_r+\Sigma_{r-1}\, d\tau$. It follows that
$$\multline
\|\dot{W}\|_r+\|W\|_r\leq K_2^{\prime\prime} \sum_{s=0}^{r}
(g_{r-s}+h_{r-s}+\beta_{r-s} )\int_0^t\|H\|_s\, d\tau \\
+ K_2^{\prime\prime} \sum_{s=0}^{r-1}(g_{r-s}+h_{r-s}+\beta_{r-s}
) \int_0^t (\|W\|_s+\|\dot{W}\|_{s}) , d\tau \endmultline \tag 8.8
$$
This proves (8.4) for $r=0$. (8.4) for $r\geq 1$ now follows by
induction, using (8.4) for $s\leq r-1$ in the terms on the second
row of (8.8) together with the interpolation (8.4). The proof of
that we can reduce to the case $B_0=B_1=0$ also for (8.5) follows
in a the same way, using that we have already proven the estimate
for smaller $m$. The prove that (8.6) can be reduced to the case
$B_0=B_1=0$, we also estimate the lower order terms using (8.5).

We now in what follows assume that $B_0=B_1=0$. Of course we could
have included these operators in the calculations that follows but
the argument becomes more clear without them.
 Lowering the indices in (8.1):
$$
\underline{G}\ddot{W}_{}+\underline{A} W_{}=\underline{G}H\tag 8.9
$$
Let ${\Cal L}_T^I $, $I\in{\Cal T}$, stand for a product of Lie
derivatives of $|I|$ vector fields in ${\Cal T}$ and let ${W}_{{}
I}= \hat{\Cal L}_T^I W_{}$. If we apply repeatedly apply Lie
derivatives ${\Cal L}_{T}$ and the projection in between, see
section 4, we obtain
$$
c^{I_1 I_2} \big(\underline{G}_{I_1}\ddot{W}_{{} I_2}+
\underline{A}_{I_1}^{} W_{{}I_2}
-\underline{G}_{I_1}H_{I_2}\big)=0 \tag 8.10
$$
where the sum is over all combination of $I_1+I_2=I$ and $c^{I_1
I_2}=1$. Here $G_I$ and $A_I$ are the operators given by (4.39)
and (4.40).  If we raise the indices again we get
$$
\ddot{W}_{{} I}+ {A} W_{{} I} =-\tilde{c}^{I_1 I_2}{A}_{I_1} W_{{}
I_2} -\tilde{c}^{I_1 I_2}{G}_{I_1}\ddot{W}_{{} I_2}+c^{I_1
I_2}{G}_{I_1}H_{I_2}\tag 8.11
$$
where $\tilde{c}^{I_1 I_2}=1$, if $|I_2|<|I|$, and $\tilde{c}^{I_1
I_2}=0$ if $|I_2|=|I|$.

 Let us define energies
$$
E_I^{}=\langle \dot{W}_{{} I},\dot{W}_{{} I}\rangle +\langle
{W}_{{} I},(A+I){W}_{{} I}\rangle, \qquad\quad E_s^{{\Cal
T}}=\sum_{|I|\leq s,I\in{\Cal T}}\sqrt{E_I^{}} \tag 8.12
$$
Note that in the sum we also included all time derivatives
$\hat{\Cal L}_{D_t}$. The reason for this is that when calculating
commutators second order time derivatives show up in the first
term on the right in (8.10). We get by differentiating (8.11), see
the end of section 3,
$$
\dot{E}_I^{}= 2\langle \dot{W}_{{} I}, \ddot{W}_{{} I} +(A+I)W_{{}
I}\rangle +\langle \dot{W}_{{} I}, \dot{G}\dot{W}_{{} I}\rangle
+\langle {W}_{{} I}, (\dot{A}+\dot{G}){W}_{{} I}\rangle\tag 8.13
$$
We now, want to estimate the right hand side by $E_r^{\Cal T}$,
where $r=|I|$. Here, by (3.27), the last two terms can be bounded
by $(h_0 c_0^{-1}+g_0)E_I$. Therefore, (8.12) can be estimated by
the $L^2$ norm of the right hand side of (8.10). In the sums with
$c^{I_1I_2}$, in (8.9), we have $|I_2|<|I|$. Since we included
time derivatives up to highest orders in $E_r^{\Cal T}$ it follows
that $\|\ddot{W}_{I_2}\|\leq E_r^{\Cal T}$. Here $G_{I_1}$ is a
bounded operator so the term with $G_{I_1}$ can be controlled.
$A_{I_1}$ is an operator of order one so the term with $A_{I_1}$
can be controlled by a constant times $\|\pa
W_{I_2}\|+\|W_{I_2}\|$. However, at this point we only have
control of tangential derivatives of $W$. One could combine the
estimate here with the curl estimate given later to get control
over all derivatives up to order $r=|I|$. Instead we will add a
lower order term to the energy such that the its time derivative
cancels the terms with $A_{I_1}$ and replaces them with lower
order terms that may be controlled. Let
$$
D_I^{}=2\tilde{c}^{I_1 I_2} \langle W_{{} I},A_{I_1} W_{{}
I_2}\rangle\tag 8.14
$$
where the sum is over all $I_1+I_2=I$, with$|I_2|<|I|$. Then
$$
\dot{D}_I^{}=2\tilde{c}^{I_1 I_2} \big(\langle \dot{W}_{{}
I},A_{I_1} W_{{} I_2}\rangle +\langle W_{{}
I},\dot{A}_{I_1}{W}_{{} I_2}\rangle+\langle W_I,{A}_{I_1}
\dot{W}_{{} I_2}\rangle\big) . \tag 8.15
$$
Hence
 $$\multline
 \dot{E}_I+\dot{D}_I =2\tilde{c}^{I_1 I_2}\big( \langle W_{{} I},\dot{A}_{I_1}{W}_{{}
 I_2}\rangle
 +\langle W_I,{A}_{I_1} \dot{W}_{{} I_2}\rangle \big)
 +\langle W_I,\dot{A} W_I\rangle \\
 +2\langle \dot{W}_{{} I},
-\tilde{c}^{I_1 I_2}{G}_{I_1}\ddot{W}_{{} I_2}+c^{I_1
I_2}{G}_{I_1}H_{I_2}+W_{{} I}\rangle +\langle \dot{W}_{{} I},
\dot{G}\dot{W}_{{} I}\rangle +\langle {W}_{{} I}, \dot{G}{W}_{{}
I}\rangle
\endmultline \tag 8.16
$$
Here, the terms on the first row can be controlled using (3.7) and
the terms on the terms on the second row can be controlled using
(3.10). We get
 $$
| \langle U, A_J W\rangle|\leq \|| h\||_{s+1,\infty} c_0^{-1} \langle
U, A U\rangle^{1/2}\langle W, A W\rangle^{1/2},\qquad |\langle U,
G_J W\rangle|\leq  \||g\||_{s,\infty}\|U\| \|W\|, \qquad |J|=s\tag
8.17
$$
Hence, we have proven that
$$
|\dot{E}_I^{}+\dot{D}_I^{}|\leq C E_r^{{\Cal T}} \sum_{s=0}^r
(h_{r-s}c_0^{-1}+g_{r-s})\big(E_s^{{\Cal T}}+\|H\|_s^{\Cal
T}\big), \qquad\quad |{D}_I^{}|\leq CE_r^{{\Cal
T}}\sum_{s=0}^{r-1}h_{r-s}c_0^{-1} E_s^{{\Cal T}}\tag 8.18
$$
so it follows that
 $$
 E_r^{\Cal T}\leq  C(h_0c_0^{-1} +g_0) \int_0^t\!\!(E_r^{\Cal T}
 +\|H\|_r)d\tau
 +C \sum_{s=0}^{r-1}
(h_{r-s}c_0^{-1}+g_{r-s})\Big( \int_0^t \!\!\big(E_s^{{\Cal
T}}+\|H\|_s^{\Cal T}\big)\, d\tau + E_s^{{\Cal T}}\Big)\tag 8.19
 $$
  Using induction and interpolation (8.4) we get the same
  inequality without
  $E_{s}^{\Cal T}$, for $s\leq r-1$, but instead multiplied by a
  constant $K_2^{\prime\prime}$ depending on (8.3) and $r$.
  Then, by a Gr\"onwall type of argument, see the beginning of the proof,
   we can also remove
  $E_{r}^{\Cal T}$ from the right hand side of (8.17), i.e.,
  with $f(t)=\int_0^t E_r^{\Cal T}\, d\tau$, we get and
  inequality $f^\prime\leq K_2^{\prime\prime} f+K_2^{\prime\prime} \sum_{s=0}^r
  (h_{r-s}+g_{r-s}) \|H\|_s^{\Cal T}$ and multiplying by the
  integrating factor $e^{K_2^{\prime\prime} t}$ gives:
 \proclaim{Lemma {8.}2}For $r\geq 0$ we have
 $$
E_r^{{\Cal T}}\leq K_2^{\prime\prime}
\sum_{s=0}^r(h_{r-s}+g_{r-s})\int_0^t \|H\|_s^{\Cal T}\, d\tau
\tag 8.20
$$
\endproclaim
This proves the estimates for tangential derivatives. We now also
want to have estimate for the curl and then by the results in
section 5 the estimates for all derivatives will follow from this.

Let ${w}_a=g_{ab}W^b$ and $\dot{w}_a=g_{ab} \dot{W}^b$ and let
$\curl w_{ab}=\pa_a w_b-\pa_b w_a$. Then  $D_t w_a=\dot{w}_a
+\dot{g}_{ab} W^b$, where $\dot{g}_{ab}=\check{D}_t g_{ab}$ and
$\dot{W}^b =\hat{D}_t W^b$. It follows that  $D_t \curl
w_{ab}=\curl \dot{w}_{ab} +\pa_a (\dot{g}_{bc}W^c )-\pa_b
(\dot{g}_{ac} W^c)$ Since the curl of $A$ vanishes it follows that
$\curl\ddot{w}=\curl \underline{H}$, if $\ddot{w}_a=g_{ab}
\ddot{W}^b$. Since the curl commutes with the Lie derivative it
therefore follows that
$$\align
|D_t \curl \dot{w}|_{r-1}^{\Cal U} &\leq C\sum_{s=0}^r g_{r-s}
|\dot{W}|_{s}^{\Cal U}+|\curl \underline{H}|_{r-1}^{\Cal U}
\tag 8.21\\
|D_t \curl w|_{r-1}^{\Cal U}&\leq C\sum_{s=0}^r
g_{r-s}|W|_{s}^{\Cal U}\tag 8.22
\endalign
$$
and by Lemma {5.}3:
$$\align
|W|_r^{\Cal U}&\leq K_1\sum_{s=1}^r g_{r-1-s}\big( |\curl
w|_{s-1}^{\Cal U}+|\div W|_s^{\Cal U}+
|W|_s^{\Cal T}) +K_1 g_{r-1}|W|_0\tag 8.23\\
|\dot{W}|_r^{\Cal U}&\leq K_1\sum_{s=1}^r g_{r-1-s}\big( |\curl
\dot{w}|_{s-1}^{\Cal U} +|\div \dot{W}|_s^{\Cal U}+
|\dot{W}|_s^{\Cal T}) +K_1 g_{r-1} |\dot{W}|_0\tag 8.24
\endalign
$$
Let
$$
C_r^{{\Cal U}}=\|\curl \dot{w}\|_{r-1}+ \|\curl w\|_{r-1}.\tag
8.25
$$
Since $\div W=\div\dot{W}=0$ it now follows that
$$
|\dot{C}_r^{{\Cal U}}| \leq K_1 \sum_{s=0}^{r}g_{r-s}
(\|\dot{W}\|_s+\|W\|_s)+ C\|\curl \underline{H}\|_{r-1}\leq K_1
\sum_{s=0}^{r}g_{r-s} ( C_s^{{\Cal U}}+E_s^{{\Cal T}})+ C\|\curl
\underline{H}\|_{r-1}. \tag 8.26
$$
This together with Lemma {8.}2 and the argument for its proof
gives
 \proclaim{Lemma {8.}3} For $r\geq 0$ we have
$$
\|\dot{W}\|_{r} +\|W\|_{r} +E_r^{{\Cal T}} \leq
K_2^{\prime\prime}\sum_{s=0}^r (g_{r-s}+h_{r-s})\int_0^t \|H\|_s\,
d\tau . \tag 8.27
$$
\endproclaim
This proves the first part of Theorem {8.}1.

In order to prove the second part of Theorem {8.}1 we replace the
energy in (8.10) by
$$
E_{r,m}^{\Cal T}=\sum_{I=I_1+I_2, \,\, I_1\in {\Cal T},\,\,
|I_1|\leq r,\,\,\, I_2=\{D_t,...,D_t\!\},\, |I_2|=m } \sqrt{E_I}
\tag 8.28
$$
i.e. we take two additional time derivatives. Noting that, the
argument leading up to Lemma {8.}2 only requires that we have at
least as many time derivatives as tangential space derivatives so
it follows from its proof that
 \proclaim{Lemma {8.}4} For $r\geq 0$ we have
 $$
E_{r,m}^{\Cal T}\leq K_2^{\prime\prime} \sum_{s=0}^{r}\sum_{j=0}^m
k_{r+m-j-s}\int_0^t \|\hat{D}_t^j{H}\|_s^{\Cal T}\, d\tau \tag
8.29
$$
\endproclaim
We have $\ddot{w}_a=-\underline{A}_a W+\underline{H}_a$ and
$\dddot{w}_a=g_{ab}\hat{D}_t \ddot{W}^b
=D_t\ddot{w}_a-(\check{D}_t g_{ab})W^b$ so, since $\curl  A W=0$,
$\curl \ddot{w}=\curl\underline{H}$ and $\curl\dddot{w}_{ab}=
\curl\,(D_t\underline{H})_{ab}-\pa_a \big(\check{D}_t
g_{bc})W^b\big) +\pa_b \big((\check{D}_t g_{ac})W^c\big)$. Hence
$|\curl \ddot{w}|_{s-1}=|\curl H|_{s-1}$ and $|\curl
\dddot{w}|_{s-1}\leq K_1\sum_{j=0}^{s} g_{s-j}\,|W|_j + |D_t\curl
\underline{H}|_{s-1}$ so using the estimates (8.23)-(8.23) with
$(W,\dot{W})$ replaced by $(\ddot{W},\dddot{W})$ gives:
$$\align
|\ddot{W}|_{r}^{\Cal U}&\leq K_1 \sum_{s=0}^{r} g_{r-1-s}
|\ddot{W}|_{s}^{\Cal T}+
K_1\sum_{s=1}^{r} g_{r-1-s} |\curl \underline{H}|_{s-1} \tag 8.30\\
|\dddot{W}|_{r}^{\Cal U} &\leq K_1 \sum_{s=0}^{r}
g_{r-1-s}(|W|_{s}^{\Cal U}+|\dddot{W}|_{s}^{\Cal T}\big)
+K_1\sum_{s=1}^{r}g_{r-1-s} |D_t\curl \underline{H}|_{s-1} \tag
8.31
\endalign
$$
Using the estimate in Lemma {8.}3, with $r$ replaced by $r-1$,
 the estimate in Lemma {8.}4, the estimate  $\|D_t^k\curl \underline{H}\|_{s-1}
 \leq \int_0^t \|D_t^{k+1} \curl \underline{H}\|_{s-1}\, d\tau$,
  and interpolation gives:
 \proclaim{Lemma {8.}5} For $r\geq 0$ and $m=0,1,2$, we have
$$
\sum_{j=0}^{m+1}\|\hat{D}_t^{j}{W}\|_{r} +E_{r,m}^{{\Cal T}} \leq
K_2^{\prime\prime}\sum_{s=0}^{r} \sum_{j=0}^m k_{r+m-j-s}\int_0^t
\|\hat{D}_t^j {H}\|_s\, d\tau . \tag 8.32
$$
\endproclaim

We now want to use the estimate for the additional time
derivatives to get estimates for $\|W\|_r$ and $\|\dot{W}\|_r$ in
(8.6) since by Lemma {5.}5:
$$
\align c_0 \|W\|_{r}&\leq  K_2^{\prime\prime}\big(\|\curl
w\|_{r-1} +\|\div W\|_{r-1}+\| W\|_{r-1,A}^{\Cal T}
+\|\dot{W}\|_{r-1}^{\Cal
T}+\sum_{s=0}^{r-1} g_{r-1-s}\|W\|_s \big)\rightalignspace \tag 8.33\\
c_0 \|\dot{W}\|_{r}&\leq  K_2^{\prime\prime}\big(\|\curl
\dot{w}\|_{r-1} +\|\div \dot{W}\|_{r-1}+\| \dot{W}\|_{r-1,A}^{\Cal
T} +\|\ddot{W}\|_{r-1}^{\Cal T}+\sum_{s=0}^{r-1}
g_{r-1-s}\|\dot{W}\|_s \big)\rightalignspace \tag 8.34
\endalign
$$
where
$$
\|W\|_{s,A}^{\Cal T} =\sum_{|I|=s,I\in{\Cal T}}\| A\hat{\Cal
L}_T^I W\|\tag 8.35
$$
The terms of order $r-1$ or less in (8.33)-(8.34) can be estimated
by Lemma {8.}3 and Lemma {8.}5 with $r$ replaced by $r-1$.
Therefore it only remains to estimate the terms involving the curl
and the operator $A$. The estimate for $\|A W_J\|$ with $|J|\leq
r-1$ we get from (8.11);
$$
{A} W_{{} J} =-\tilde{c}^{J_1 J_2}{A}_{J_1} W_{{} J_2} -c^{I_1
J_2}\big( {G}_{J_1}\ddot{W}_{{} J_2}-{G}_{I_1}H_{J_2}\big)\tag
8.36
$$
where $\tilde{c}^{J_1 J_2}=0$ if $|J_2|=r-1$. Here the terms in
the parenthesis can be estimated by Lemma {8.}5, with $r$ replaced
by $r-1$. Note that $\|H_{J_2}\|\leq \int_0^t \|\dot{H}_{J_2}\|\,
d\tau $ since we assume that $H$ vanishes to all orders as $t\to
0$. Since by (3.8) $A_{J_1}$ is order one and $|J_2|\leq r-2$ in
the first term on the right it follows that also this term can be
estimated by Lemma {8.}3 with $r$ replaced by $r-1$. We hence get
 $$
 \|W\|_{r-1,A}^{\Cal T}\leq K_2^{\prime\prime}\sum_{s=0}^{r-1} \int_0^t
\big(k_{r-1-s}\|\ddot{H}\|_s+k_{r-s}\|\dot{H}\|_s+k_{r+1-s}\|{H}\|_s\big)\,
d\tau \tag 8.37
$$
From (8.9) we also get the estimate for $\|A\dot{W}_J\|$ with
$|J|\leq r-1$ by letting one of the derivatives in (8.9) be a time
derivative so $I=J+\{D_t\}$. If we write this out we get
$$
{A} \dot{W}_{{} J} =-\tilde{c}^{J_1 J_2}{A}_{J_1} \dot{W}_{{} J_2}
-c^{J_1 J_2}\dot{{A}}_{J_1}{W}_{{} J_2} -c^{J_1
J_2}\big({G}_{J_1}\dddot{W}_{{} J_2}+\dot{G}_{J_1}\ddot{W}_{{}
J_2}-{G}_{J_1}\dot{H}_{I_2}-\dot{G}_{J_1}H_{J_2}\big) \tag 8.38
$$
where $\tilde{c}^{J_1 J_2}=0$ if $|J_2|=r-1$. Here the terms in
the parenthesis can be estimated by Lemma {8.}5, with $r$ replaced
by $r-1$. Since $A_{J_1}$ and $\dot{A}_{J_1}$ is order $1$, see
(3.8), we can also estimate all the other by Lemma {8.}3 with $r$
replaced by $r-1$, apart from the term with $\dot{A} W_{J}$. This
term is estimated by $\|\dot{A} W_{J}\|\leq C
h_0(\|W\|_r+\|W\|_{r-1})$, where the last term can again be
estimated by Lemma {8.}3 with $r$ replaced by $r-1$. We hence get
 $$
 \|\dot{W}\|_{r-1,A}^{\Cal T}\leq  h_0 \|W\|_r+
 K_2^{\prime\prime}\sum_{s=0}^{r-1}\int_0^t
\big(k_{r-1-s}\|\ddot{H}\|_s+k_{r-s}\|\dot{H}\|_s+k_{r+1-s}\|{H}\|_s\big)\,
d\tau \tag 8.39
$$
(8.33) and (8.34) together with (8.37) and (8.39) therefore gives
$$
 c_0 \|W\|_r \leq K_2^{\prime\prime}\|\curl w\|_{r-1}+
 K_2^{\prime\prime}\sum_{s=0}^{r-1} \int_0^t
\big(k_{r-1-s}\|\ddot{H}\|_s+k_{r-s}\|\dot{H}\|_s+k_{r+1-s}\|{H}\|_s\big)\,
d\tau \tag 8.40
$$
and
$$
 c_0 \|\dot{W}\|_r \leq h_0 c_0^{-1}\|\curl w\|_{r-1}+K_2^{\prime\prime} \|\curl \dot{w}\|_{r-1}
 +K_2^{\prime\prime}\sum_{s=0}^{r-1} \int_0^t
\big(k_{r-1-s}\|\ddot{H}\|_s+k_{r-s}\|\dot{H}\|_s+k_{r+1-s}\|{H}\|_s\big)\,
d\tau \tag 8.41
$$
It therefore only remains to control the curl. With $C_r^{\Cal U}
=\|\curl\dot{w}\|_{r-1}+\|\curl w\|_{r-1}$ it hence follows from
(8.26) and (8.40)-(8.41) that
$$
|\dot{C}_{r}^{\Cal U}|\leq K_2^{\prime\prime} C_r^{\Cal U}+
K_2^{\prime\prime}\sum_{s=0}^{r-1} \int_0^t
\big(k_{r-1-s}\|\ddot{H}\|_s+k_{r-s}\|\dot{H}\|_s+k_{r+1-s}\|{H}\|_s\big)\,
d\tau +\|\curl \underline{H}\|_{r-1}  \tag 8.42
$$
where we also used that, by Lemma {8.}3,
$\|\dot{W}\|_0+\|W\|_0\leq K_2^{\prime\prime}\int_0^t\|H\|_0\,
d\tau$. Integrating this equation gives a bound for $C_r^{\Cal U}$
in terms of the integral in (8.38). Using this bound in
(8.40)-(8.41) then gives
 \proclaim{Lemma {8.}6} For $r\geq 1$ we have
 $$
 \|\dot{W}\|_r+\|W\|_r\leq
 K_2^{\prime\prime}\sum_{s=0}^{r-1}\int_0^t
\big(k_{r-1-s}\|\ddot{H}\|_s+k_{r-s}\|\dot{H}\|_s+k_{r+1-s}\|{H}\|_s+\|\curl
\underline{H}\|_{s}\big)\, d\tau \tag 8.43
$$
\endproclaim
This concludes the proof of Theorem {8.}1.

\head 9. Existence and tame estimates for the inverse of the
modified linearized operator.\endhead We want to show existence
and tame estimates for the inverse of the linearized operator.
However, first we will show existence and estimates for the
modified linearized operator $L_1$ given by (2.55):
$$
L_1 W=F,\qquad W\big|_{t=0}=\dot{W}\big|_{t=0}=0\tag 9.1
$$
where $F$ is smooth and vanishes to all orders as $t\to 0$.

Let
 $$
 n_s=\||g\||_{s+2,\infty}+\||\omega\||_{s+1,\infty}
 +\||h\||_{s+3,\infty},\tag 9.2
 $$
 and let $K_3^{\prime\prime}$ be denote a continuous function of
 $$
 n_0+c_0^{-1}+T^{-1}+c_1+r\tag 9.3
 $$
 which in what follows also depends on the order of differentiation
 $r$.
 In the proof that follows we will use the interpolation
 $$
 n_s n_r\leq K_3^{\prime\prime} n_{s+r}\tag 9.4
 $$

 \proclaim{Theorem {9.}1} Suppose that (2.4)-(2.5) hold
for $0\leq t\leq T$ and suppose also that $x$ is smooth for $0\leq
t\leq T$ an that $T\leq 1$. Then (9.1), where $F$ is smooth and
vanishing to all orders as $t\to 0$, has a smooth solution for
$0\leq t\leq T$. It satisfies the estimates
$$
\|\dot{W}\|_r+\|W\|_r\leq K_3^{\prime\prime} \sum_{s=1}^r
n_{r-s}\int_0^t\|F\|_s\, d\tau,\qquad r\geq 1\tag 9.5
$$
\endproclaim
\demo{Proof of Theorem {9.}1} The proof of existence and the
estimate (9.5) for the modified linearized operator uses the
orthogonal decomposition of the vector field into its divergence
free part and a gradient of a function vanishing on the boundary.
The solution will be further divided into four parts and for each
of them we use either Theorem {6.}1, Theorem {7.}1 or Theorem
{8.}1. This gives us the estimates in Theorem {9.}2. These
estimates imply the estimate (9.5). Furthermore, the estimates
holds for iterates, i.e given an iterate for $W$ we define $\delta
h=\Psi^\prime(x)W$ by (9.10) and then define $W_{ij}$ by
(9.8)-(9.20) and (9.12)-(9.13). This gives us a new iterate for
$W$. Theorem {9.}2 then gives us uniform bounds for the iterates
and applied to the equations for differences of iterates gives us
convergence, see \cite{L2,L3}.

 Now, the solution $W$ of
$$
L_1 W=F\tag 9.6
$$
can be obtained as the sum of four terms, see section 3,
$$
W=W_0+W_1,\qquad W_1=W_{10}+W_{11},\qquad W_{0}=W_{00}+W_{01} \tag
9.7
$$
where $W_0$ is divergence free and
$$
W_{1i}^a=g^{ab}\pa_b q_{1i},\qquad i=0,1,\qquad
q_{1i}\big|_{\pa\Omega}=0,\qquad \triangle q_{11}=\varphi,\qquad
\triangle q_{10}=-e^\prime(h)\delta h \tag 9.8
$$
and $\varphi$ satisfies an ordinary differential equation:
$$
D_t^2\varphi +\div\Phi\, \varphi=\div F+\div\Phi\, e^\prime(h)\delta
h,\qquad \qquad \varphi\big|_{t=0}=\dot{\varphi}\big|_{t=0}=0 \tag 9.9
$$
where $\div\Phi=D_t^2 e(h)+ D_t\div V$ and $\delta h$ satisfies
the wave equation
$$\align\leftalignspace
D_t^2 (e^\prime(h)\delta h) -\triangle\delta h &=\triangle\big(
(\pa_c h) W^c\big) -\div B_1\dot{W}+2\dot{\sigma}\div \dot{W}
-\div B_0 W-(D_t^2 e(h)+\dot{\sigma}^2)\div W,\rightalignspace \tag 9.10\\
 \delta h\big|_{\pa\Omega}&=0,\qquad\qquad
\delta h\big|_{t=0}=\dot{\delta h}\big|_{t=0}=0.\tag 9.11\endalign
$$

Here, the divergence free parts satisfy the evolution equation for
the normal operator, (3.16)
$$
\ddot{W}_{00}+A W_{00}-B_{10} \dot{W}_{00}-B_{00} W_{00}
=-AW_{10}+B_{11} \dot{W}_1+B_{01} W_1 +PF, \qquad
W_{00}\big|_{t=0}\!=\dot{W}_{00}\big|_{t=0}\!=0\tag 9.12
$$
and
$$
\ddot{W}_{01}+A W_{01} -B_{10} \dot{W}_{01}-B_{00} W_{01}
=-AW_{11},\qquad\qquad
W_{01}\big|_{t=0}\!=\dot{W}_{01}\big|_{t=0}\!=0\tag 9.13
$$

If
$$
E_r=\sum_{i,j=0}^1 \|\dot{W}_{ij}\|_r+\|W_{ij}\|_r+\|\delta h
\|_{r}\tag 9.14
$$
then we get from Theorem {9.}2 that for $r\geq 1$
$$
E_r\leq K_3^{\prime\prime}\sum_{s=1}^r n_{r-s} \int_0^t
(\|F\|_s+E_s)\, d\tau\tag 9.15
$$
Using a Gr\"onwall type of argument and induction as in section 8
it follows that for $r\geq 1$ we have
$$
E_r\leq K_3^{\prime\prime}\sum_{s=1}^r n_{r-s}\int_0^t \|F\|_s\,
d\tau\tag 9.16
$$
which proves Theorem {9.}1.
 \qed\enddemo

 \proclaim{Theorem {9.2} } If $r\geq 1$ we have
 $$\align
 \|\dot{W}_{11}\|_{r} +\|W_{11}\|_r &\leq K_3^{\prime\prime}\sum_{s=1}^r n_{r-s}
 \int_0^r \|F\|_s+\|W_{10}\|_{s} \, d\tau\tag 9.17\\
\|\ddot{W}_{11}\|_{r} &\leq K_3^{\prime\prime}\sum_{s=1}^r n_{r-s}
 \Big( \|F\|_s+\|W_{10}\|_{s}+\int_0^r \|F\|_s+\|W_{10}\|_{s} \,
 d\tau\Big)\rightalignspace \tag 9.18\\
\|\dot{W}_{01}\|_{r} +\|W_{01}\|_r &\leq
K_3^{\prime\prime}\sum_{s=1}^r n_{r-s}
 \int_0^r\big( \|F\|_s+ \|W_{10}\|_{s} \big)\, d\tau\tag 9.19\\
 \|\dot{W}_{00}\|_{r} +\|W_{00}\|_r &\leq K_3^{\prime\prime}\sum_{s=1}^r n_{r-s}
 \int_0^r \!\!\big(\|F\|_s+ \|\pa W_{10}\|_{s}\big)\,
 d\tau\rightalignspace \tag 9.20\\
 \leftalignspace
 \|\dot{W}_{10}\|_{r} +\|W_{10}\|_{r+1}&\leq K_3^{\prime\prime}\sum_{s=1}^r n_{r-s}
 \int_0^t (\|\dot{W}\|_s +\|W\|_s) \, d\tau\tag 9.21
 \endalign
$$
where $W=W_0+W_1$, $W_1=W_{10}+W_{11}$ and $W_0=W_{00}+W_{01}$.
 \endproclaim
 \demo{Proof of Theorem {9.}2} (9.17) follows directly from applying
Lemma {9.}3 below to (9.9).
Here we write $\div F_1=\div F +\div\Phi\,e^\prime(h)\delta h=
\div \big( F-\div\Phi W_{10}\big)+ (\pa_a \div\Phi)W_{10}^a $
and use Theorem {6.}1.
This also gives the additional estimate (9.18).

 To prove (9.19) we use the second part of Theorem {8.}1 applied
 to (9.13). Since $\curl AW_{11}=0$ and we have an estimate
 for and additional time derivative in (9.18), (9.19) follows.

 By the first part of Theorem {8.}1:
$$
\|\dot{W}_{00}\|_{r} +\|W_{00}\|_r \leq
K_3^{\prime\prime}\sum_{s=0}^r n_{r-1-s}
 \int_0^r \!\!\big(\|F\|_s+ \|\pa W_{10}\|_{s} +\|\dot{W}_1\|_{s}+\|W_1\|_s\big)\,
 d\tau\tag 9.22
 $$
Using (9.21) and (9.17) gives (9.20). Note that the sum above
starts at $s=0$ but since the lower norms are included in the
higher norms and since we replaced $n_{r-s-1}$ by $n_{r-s}$ we can
start the sum in (9.20) from $s=1$.

 The estimate (9.21) follows
from Theorem {7.}1 applied to (9.10). If ${f}$ denotes the right
hand side of (9.10) then
$$
\|\dot{{f}}\|_{s-2}\leq K_3^{\prime\prime}\sum_{k=0}^s n_{s-1-k}
\big(\|\dot{W}\|_{k} +\|W\|_k\big)\tag 9.23
$$
Since
 $\|{f}\|_0\leq K_3^{\prime\prime}(\|W\|_2+ \|\dot{W}\|_1)$ we obtain from (7.13)
 for $r\geq 2 $
$$
\|\delta h\|_r\leq K_3^{\prime\prime}\sum_{s=0}^r n_{r-1-s}
\int_0^t (\|\dot{W}\|_s+\|W\|_s)\, d\tau\tag 9.24
$$
For $r=1$ we use (7.12), which gives (9.24) also for $r=1$ if we
write ${f}={f}_1+\triangle (h_c W^c)$. This proves that $\|\delta h\|_r$
is bounded by the right hand side of (9.21).
To get the estimate also for $W_{10}$ we use Theorem {6.}1 to estimate
the solution of the Dirichlet problem in (9.8) using that
$\|e^\prime(h)\delta h\|_r\leq
K_3^{\prime\prime} \sum_{s=0}^r \|h\|_{r-s,\infty} \|\delta h\|_s $.
This gives that $\|\pa W_{10}\|_r$ and  $\|W_{10}\|_r$ are bounded by
the right hand side of (9.21) but it remains to show that
$\| D_t^{r+1} W_{10}\|$ is bounded by the right hand side of (9.21)
which follows from Theorem {7.}4.
 \qed\enddemo

\proclaim{Lemma {9.}3} Let $\varphi$ be the solution of
 $$
 D_t^2 \varphi+ k \varphi=f,\qquad
 \dot{\varphi}\big|_{t=0}=\varphi\big|_{t=0}=0, \tag 9.25
 $$
where $D_t^j f_1\big|_{t=0}=0$
for $j\geq 0$ and  $k=\div\Phi=D_t^2e(h)+D_t \div V$.
 Then for $r\geq 1$ we have
 $$
 \align
 \|\dot{\varphi}\|_{r-1}+\|\varphi\|_{r-1}&\leq K_3^{\prime\prime}\sum_{s=0}^{r-1}
  n_{r-1-s}\int_0^t
 \!\!\! \|f\|_s\, d\tau,\tag 9.26\\
 \|\ddot{\varphi}\|_{r-1}&\leq K_3^{\prime\prime}\sum_{s=0}^{r-1}  n_{r-1-s}\Big(
 \|f\|_s+\! \int_0^t\!\!\! \|f\|_s\, d\tau\Big)\tag 9.27
 \endalign
 $$
Furthermore let $W_1$ and $F_1$ be defined by
 $$\align
 W_1&=\nabla q,\qquad \triangle q=\varphi,\qquad
 q\big|_{\pa\Omega}=0\tag 9.28\\
F_1&=\nabla q^\prime \qquad \triangle q^\prime=f,\qquad
 q^\prime\big|_{\pa\Omega}=0\tag 9.29
\endalign
 $$
Then for $r\geq 1$ we have
 $$
 \align
 \|\dot{W}_1\|_{r}+\|W_1\|_{r}&\leq K_3^{\prime\prime}\sum_{s=0}^r  n_{r-s}\int_0^t
 \!\!\! \|F_1\|_s\, d\tau,\tag 9.30\\
 \|\ddot{W}_1\|_{r}&\leq K_3^{\prime\prime}\sum_{s=0}^r  n_{r+1-s}\Big(
 \|F_1\|_s+\! \int_0^t\!\!\! \|F_1\|_s\, d\tau\Big)\tag 9.31
 \endalign
 $$
\endproclaim
\demo{Proof of Lemma {9.}3}
(9.25) is an ordinary differential equation and (9.26)-(9.27) follows
as in the proof of Proposition {10.}1 in \cite{L3}. (9.25) can be
written
$$
\hat{D}_t^2 \varphi -2\dot{\sigma} \hat{D}_t \varphi
+k^\prime \varphi=f,\qquad k^\prime=\dot{\sigma}^2 +D_t^2 e(h)\tag 9.32
$$
and so
$$
\div \big( \ddot{W}_1-2\dot{\sigma}\dot{W}_1+k^\prime W_1 -F_1\big)=
 -2(\pa_a\dot{\sigma}) \dot{W}_1^a
+(\pa_a k^\prime) W_1^a\tag 9.33
$$
and hence
$$
\div\big( \hat{D}_t^j \ddot{W}_1 - \hat{D}_t^j \big(
2 \dot{\sigma}\dot{W}_1-k^\prime W_1+F_1\big)\big)=
-\hat{D}_t^j \big( 2(\pa_a\dot{\sigma}) \dot{W}_1^a
-(\pa_a k^\prime) W_1^a\big)\tag 9.34
$$
We claim that
$$
P \hat{D}_t^j W_1=P B_j(W_1,...,\hat{D}_t^{j-1} W_1),\qquad
B_j(W_1,...,\hat{D}_t^{j-1} W_1)=\sum_{i=0}^{j-1} {\tsize{\binom{j}{i}}}
(\check{D}_t^{j-i} g_{ab})\hat{D}_t^{i} W_1 \tag 9.35
$$
In fact $0=P D_t^j \pa_a q= P D_t^j\big( g_{ab} W_1^b\big)
=P \sum_{i=0}^j \binom{j}{i}
(\check{D}_t^{j-i} g_{ab})\hat{D}_t^{i} W_1$. Furthermore, let
$$
\triangle q^j=-\hat{D}_t^j \big( 2(\pa_a\dot{\sigma}) \dot{W}_1^a
-(\pa_a k^\prime) W_1^a\big),\qquad
q^j\big|_{\pa\Omega}=0\tag 9.36
$$
To say that $\div H=0$ is equivalent to saying that $(I-P) H=0$ so
it follows from (9.34)-(9.36) that
$$
\hat{D}_t^{j} W_1=PB_j(W_1,...,\hat{D}_t^{j-1} W_1)
+(I-P) \hat{D}_t^{j-2}\big(
2 \dot{\sigma}\dot{W}_1-k^\prime W_1+F_1\big) +\nabla q^{j-2}\tag 9.37
$$
Since the projection $P$ maps $L^2$ to $L^2$ and since inverting (9.36)
maps $L^2$ to $H^1$ (in fact to $H^2$) it therefore follows that:
$$
\|D_t^{r+1} W_1\|\leq K_3^{\prime\prime}
\sum_{s=0}^{r} n_{r-s}\|D_t^s W_1\|+\sum_{s=0}^{r-1} n_{r-2-s} \|D_t^s F_1\|\tag 9.38
$$
Sine it also follows from (9.26) that
$$
 \align
 \|\pa\dot{W}_1\|_{r-1}+\|\pa W_1\|_{r-1}&\leq K_3^{\prime\prime}\sum_{s=0}^r  n_{r-s}\int_0^t
 \!\!\! \|F_1\|_s\, d\tau,\tag 9.39\\
 \|\pa \ddot{W}_1\|_{r-1}&\leq K_3^{\prime\prime}\sum_{s=0}^r  n_{r-s}\Big(
 \|F_1\|_s+\! \int_0^t\!\!\! \|F_1\|_s\, d\tau\Big)\tag 9.40
 \endalign
 $$
the lemma follows from also estimating
 $\|D_t^s F_1(t,\cdot)\|\leq \int_0^t \|D_t^{s+1} F_1(\tau ,\cdot)\|\, d\tau$.
\qed\enddemo

 \head 10. Estimates for the enthalpy in terms of the
coordinate.\endhead
 We have now proved that the linearized operator is invertible.
 However, since we think of $h=\Phi(x)$ as a functional of $x$ we
 must also estimate the $L^\infty$ norms of $h$ in terms of the
 $L^\infty$ norms of $x$. For the corresponding problem for the
 incompressible case in \cite{L3} we could take advantage of the
 Schauder
 estimates. However for the wave equation there are no $L^\infty$
 estimates that do not loose regularity.
 For wave equations it is best to get the
  $L^{\infty}$ norms from the
 $L^2$ norms using Sobolev's lemma.
 These estimates where obtained in Corollary {7.}5.
 However, in the estimates there the $L^\infty$ norm also occurred
 in the right hand side, due to that we assumed that $e^\prime =e^\prime (h)$.
 This estimate can easily be improved by estimating the $L^2$
 norm of the solution of the nonlinear wave equation instead.
 We will however for simplicity assume that $e(h)=c h$, so that
 $e^\prime (h)=c$, where $0<c<\infty$ is a constant. In that case there
 are no $h$ terms in the results in section 7. We have
$$
D_t^2 c h -\triangle h=(\pa_i V^j)(\pa_j V^i), \qquad
h\big|_{\pa\Omega}=0\tag 10.1
$$
The estimates in Corollary {7.}5 were however formulated for
vanishing initial data.  Therefore let $\tilde{h}=h-h_0$ and where
$h_0$ satisfies the equation (10.1) to all orders as $t\to 0$, if
$x-x_0$ vanishes to all orders as $t\to 0$, see section 2.
 It follows that
$$
D_t^2 (c\tilde{h})-\triangle \tilde{h}=-D_t^2 (c h_0)+\triangle
h_0-(\pa_i V^j)(\pa_j V^i)\tag 10.2
$$
vanishes to all orders as $t\to 0$. We can therefore apply
Corollary {7.}5, which gives:

\proclaim{Lemma {10.}1 } We have
 $$
 \||g\||_{s+1,\infty}+\||\omega\||_{s,\infty}
 \leq K_1 \||x\||_{s+2,\infty}\tag 10.3
 $$
 Suppose that $e(h)=ch$, where
$0<c<\infty$. Then for $r\geq 2$
$$
\||h\||_{r,\infty}\leq K_2 \||x\||_{r+r_0+1,\infty}.\tag 10.4
$$
Here $K_i$ are as in Definition {5.}2 and $K_2$ also depends on a
bound for $\||h_0\||_{r+r_0+1,\infty}$.
\endproclaim
\demo{Proof} The first inequality follows directly from the
definitions and interpolation. If ${f}$ denotes the right hand side
of (10.2) then
$$
\||{f}\||_{s,\infty}\leq K_2 \big( \||
h_0\||_{s+2,\infty}+\||x\||_{s+2,\infty} \big)\tag 10.5
$$
As pointed out above if $e(h)=ch$ then $h_r$ in section 7 vanishes
and  by Corollary {7.}5 and interpolation we have
 $$
 \||\tilde{h}\||_{r,\infty}\leq K_2\big( \||x\||_{r+r_0+1,\infty}
 +\||h_0\||_{r+r_0+1,\infty}\big)\tag 10.6
 $$
 However since we are just looking on fixed initial data we can
also include the norms of $h_0$ in the constants and since in fact
we also have a lower bound for $\|| x\||_{1,\infty}\geq C c_1>0$,
by the coordinate condition, the lemma follows.
  \qed\enddemo

\demo{Remark} Note that $K_2$ in (10.4) depends on $h_0$. However,
$h_0$ is a function which is fixed once we fixed initial data so
this just leads to an $r$ dependence of the constant.
\enddemo

It now follows that with $K_i$ as in Definition {5.}2 and
$K_i^\prime$ as in sections 7, 8, and 9 we have
$$
K_i^\prime\leq K_{i+r_0+1}\tag 10.7
$$

\head 11. Estimates for the physical condition and coordinate
condition. \endhead

We assume that the physical condition and the coordinate condition
initially at time $0$ for some constants $c_0>0$ and $c_1<\infty$
and we need to show that this implies that they will hold with $
c_0$ replaced by $c_0/2$ and $c_1$ replaced by $2c_1$,  for
 $0\leq t\leq T$, if $T$ is sufficiently small.

 Now, for the coordinate condition this is easy since we can just
 estimate the physical condition by the time derivative of $g$,
 which can be estimated by $\||x\||_{2,\infty}$:
 \proclaim{Lemma {11.}1} Let $M(t)=\sup_{y\in\Omega}\sqrt{|\pa x/\pa y|^2+|\pa y/\pa
 x|^2}$. Then
 $$
 M(t)\leq 2M(0),\qquad \text{for}\qquad t\leq T,\qquad
 \text{if}\qquad T \||x\||_{2,\infty}M(0)\leq 1/8\tag 11.1
 $$
 Let $N(t)=\sup_{y\in\pa\Omega} |\na_N h|^{-1}$. Then assuming
 that $T$ is so small that (11.1) hold we have
 $$
 N(t)\leq 2N(0)\qquad \text{for}\qquad t\leq T,\qquad
 \text{if}\qquad T \||h\||_{2,\infty}M(0)N(0)\leq 1/8\tag 11.2
 $$
 \endproclaim
 \demo{Proof} We have $|D_t \pa x/\pa y|\leq \||x\||_{2,\infty} $
 and $|D_t \pa y/\pa x|\leq |\pa y/\pa x|^2 |D_t \pa x/\pa y|$ so
 $ M^\prime(t)\leq (1+M^2) \||x\||_{2,\infty}\leq 2M^2
 \||x\||_{2,\infty}$, since also $M(t)\geq 1$. Hence
 $$
 M(t)\leq M(0) \big(1-2\||x\||_{2,\infty} M(0) t\big)^{-1}, \qquad
 \text{when}\qquad 2\||x\||_{2,\infty} M(0) t< 1.\tag 11.3
 $$
 Now, $\na_N h=N^a \pa_a h$, where $N$ is the unit normal,
 so $D_t \na_N h=\na_N D_t h+ (D_t N^a) \pa_a h=\na_N D_t h+(D_t
 N^a) g_{ab} N^b \na_N h$, since $h\big|_{\pa\Omega}=0$.
 Furthermore $0=D_t (g_{ab} N^a N^b)= 2g_{ab} (D_t N^a) N^b+(D_t
 g_{ab}) N^a N^b$ and $N^a=(\pa y^a/\pa x^i) N^i$, where
 $\delta_{ij} N^i N^j=1$. Hence
 $|D_t\na_N h|\leq M \big(|\pa D_t h|+ |\pa D_t x| |\na_N h|\big)$
 Therefore if $N(t)=\sup_{y\in\pa\Omega}|\na_N h|^{-1}$, we have
 $N^\prime\leq M \||h\||_{2,\infty}N^2 + M \||x\||_{2,\infty}N/2$
 and if we use (11.1) and multiply with the
 integrating factor, $\tilde{N}(t)=N(t)e^{-tM(0) \||x\||_{2,\infty}}$
 we get $\tilde{N}^\prime\leq 2e^{1/8}M(0) \||h\||_{2,\infty}
 \tilde{N}^2$. Hence
 $$
 {N}(t)\leq {N}(0) e^{1/8} \big(1-N(0)2e^{1/8} M(0) \||h\||_{2,\infty}
 t\big)^{-1}\!\!\!\!,\quad \text{when}\quad
 N(0)2 e^{1/8}M(0) \||h\||_{2,\infty}t<1\tag 11.4
 $$
 This proves the lemma. \qed
 \enddemo

 To satisfy the condition in (11.1) we just need to choose $T$ so small that
 $T\||x\||_{2,\infty}c_1  \leq 1/8$.
 Remarking that $x=u+x_0$, where $x_0$ is a fixed and that in the
 Nash-Moser iteration we will only apply our estimates to functions
 satisfying $\||u\||_{r_0+4,\infty}\leq 1$.

However for the physical condition this is a bit more difficult.
One has to control the $\||h\||_{2,\infty}$ and the estimate for
$h$ in terms of $x$ used interpolation so they are in terms of a
constant $K$ that is a continuous function of $T^{-1}$ and tends
to infinity as $T\to 0$. In the compressible case we never used
interpolation in time so the corresponding estimate there was
easier. We therefore have to redo the estimates for the wave
equation for the lowest norms without using interpolation in time.
This is however, follows from standard estimates for the wave
equation. Those we have here also work if we do not use
interpolation. If we do not use interpolation, then the estimates
in sections 5,6 and 7 still hold, with constants independent of
$T\leq 1$, but instead of depending linearly on the highest norms
of $x$ they are polynomials in the highest norms of $x$ occurring
in the estimates.
 \proclaim{Lemma {11.}2} There are continuous increasing functions $C_r$ such
 that for $T\leq 1$ we have
 $$
 \||h\||_{r,\infty}\leq
 C_r\big(c_1,\||x\||_{r+r_0+1,\infty},\||h_0\||_{r+r_0+1,\infty} \big)
 \tag 11.5
$$
 \endproclaim
\demo{Proof} The proof is the same as the proof of Lemma {10.}1
using that the estimates in sections 5,6 and hold with constant of
the form above. \qed\enddemo

Summing up, we have hence proven that
 \proclaim{Lemma {11.}3} Let $C_2$ be as in Lemma {11.}2
 and let $c_1$ and $c_0$ be constants such that the coordinate
 condition (2.13) and the physical condition (2.14) hold when $t=0$.
 Suppose $0<T\leq 1$ is fixed such that
$$
Tc_1 (1+\||x_0\||_{2,\infty})\leq 1/8,\qquad  T c_1
C_2\big(2c_1,1+\||x_0\||_{r_0+3,\infty},\||
h_0\||_{r_0+3,\infty}\big) \leq c_0/8,\tag 11.6
$$
where $C_2$ is as in Lemma {11.}2.
Then for $0\leq t\leq T$, the coordinate condition hold with $c_1$
replaced by $2c_1$ and $c_0$ replaced by $c_0/2$ if
$$
\|| u\||_{r_0+4,\infty}\leq 1\tag 11.7
$$
where $r_0=[n/2]+1$ is the Sobolev exponent.
Here $u=x-x_0$ and $x_0$ is the approximate solution.
 \endproclaim
\demo{Proof} It follows from Lemma {11.}2 that
$$
\|| h\||_{2,\infty}\leq C_2\big(c_1,\|| u+x_0\||_{r_0+3,\infty},
\|| h_0\||_{r_0+3,\infty}\big) \tag 11.8
$$
In view of (11.2), the physical condition with $c_0$ replaced by
$c_0/2$ hold if $T$ is so small that (11.5) hold and
$$
T c_1 C_2\big(2c_1,\|| u+x_0\||_{r_0+3,\infty},\||
h_0\||_{r_0+3,\infty}\big) \leq c_0/8\qed\tag 11.9
$$
\enddemo

 Recalling again that in the Nash-Moser iteration we
will only consider $u$ for which (11.7) hold. From now on we will
therefore assume that $0<T\leq 1$ is fixed and so small that
(11.6) hold.

 \head 12. Tame estimates for the inverse of the
linearized operator in terms of the coordinate.\endhead

 We have
\proclaim{Theorem {12.}1} Suppose that $T>0$ is so small that the
conditions in Lemma {11.}3 hold. Suppose also that $x=u+x_0$ and
$\delta \Phi$ are smooth in $[0,T]\times\overline{\Omega}$ and
that $\delta \Phi $ and $u$ vanish to infinite order as $t\to 0$.
 Then there are constants $K$, depending on the approximate solution
 $(x_0,h_0)$, on $(c_0,c_1)$ and on $r$, such that there is a
 smooth solution $\delta x$ of
 $$
 \Phi^\prime(x)\delta x=\delta \Phi,\qquad\text{in}\qquad [0,T]\times\overline{\Omega}\qquad
 \qquad \delta x\big|_{t=0}=\dot{\delta x}\big|_{t=0}=0\tag 12.1
 $$
satisfying
 $$
 \|\delta\dot{ x}\|_r+\|\delta x\|_r\leq K\sum_{s=1}^r \||x\||_{r+r_0+4-s,\infty}
 \int_0^t\|\delta \Phi\|_s\, d\tau\tag 12.2
 $$
 for $r\geq 1$ if
$$
 \||u\||_{r_0+4,\infty}\leq 1\tag 12.3
 $$
 Moreover
 $$
 \||x\||_{r+r_0+4,\infty}\leq K+\||u\||_{r+r_0+4,\infty}\tag 12.4
 $$
 \endproclaim
 \demo{Proof} First show existence the equation
 where the vector field is expressed in the Lagrangian frame,
 $W^a=(\pa y^a/\pa x^i) \delta x^i$ and $F^a=(\pa y^a/\pa
 x^i)\delta \Phi^i$:
 $$
 L_0 W=F,\qquad\text{in}\qquad [0,T]\times\overline{\Omega}\qquad
 \qquad W\big|_{t=0}=\dot{W}\big|_{t=0}=0\tag 12.5
 $$
and that it satisfies the estimate
 $$
 \|\dot{W}\|_r+\|W\|_r\leq K\sum_{s=1}^r \||x\||_{r+r_0+4-s,\infty}
 \int_0^t\|F\|_s\, d\tau\tag 12.6
 $$
 for $r\geq 1$ if
 We have already proved this for $L_1 W=F$ in
 Theorem {9.}1, using Lemma {10.}1 and Lemma {11.}3.
 Therefore it remains to prove the result for $L_0 W=L_1 W-B_3
 W=F$, where $B_3$ is given by (2.63).
We have $\|B_3 W\|_s\leq
K\sum_{k=0}^s\||x\||_{s+r_0+3-k,\infty}\|W\|_k\leq
K\sum_{k=1}^s\||x\||_{s+r_0+4-k,\infty}\|W\|_k$, if $s\geq 1$.
Applying the theorem to the equation $L_1 W=F+B_3 W$ and using
interpolation we get that for $r\geq 1$
$$
\|\dot{W}\|_r+\|W\|_r\leq K\sum_{s=1}^r
\||x\||_{r+r_0+4-s,\infty}\int_0^t (\|F\|_s+\|W\|_s)\, d\tau\tag 12.7
$$
If we put up an iteration $L_1 W^{k+1} =F-B_3 W^{k}$, for $k\geq
0$ and $W^{0}=0$ then (12.6) is going to be true with $W$ in the
right hand side replaced by $W^{(k)}$ and in the left by
$W^{k+1}$. It is easiest to first show convergence to a smooth
solution and then afterwards prove the estimate (12.2). To show
convergence we can just consider the estimate (12.6) where we
include the norms of $x$ in the constants and estimate the lower
order norms by higher order norms. Let
$\tilde{W}^{k+1}=W^{k+1}-W^{k}$. Then $L_1 \tilde{W}^{1}=F$ and
$L_1\tilde{W}^{k+1}=-B_3\tilde{W}^{k}$, for $k\geq 1$. Hence if
$E^k_r=\sum_{j=1}^k
\|\hat{D}_t\tilde{W}^{j}\|_r+\|\tilde{W}^{j}\|_r$ we have
$E_{r}^{k+1}\leq C_r\int_0^t (\|F\|_r+E_r^{k})\, d\tau$. Using a
Gr\"onwall type of argument one therefore get uniform bounds
 $E_r^k\leq C_r^\prime \int_0^t \|F\|_r\, d\tau$, for $0\leq t\leq T\leq 1$.
 This proves convergence of to a smooth solution. $W$.
 Once we have a smooth solution it will satisfy the estimate
 (12.7). By a Gr\"onwall type of argument and induction as in
 section 8 it follows that the solution also satisfies the
 estimate (12.6), for some other constant $K$.

 Finally, we want to deduce the estimate for $\delta x$ and
 $\delta\Phi$.
The estimate (12.6) is in terms of $W^a=\delta x_i \pa y^a/\pa
x_i$ and $F^a=\delta \Phi_i \pa y^a/\pa x^i$,
 turning them into an estimate for $\delta x$  and $\delta \Phi$
 just produces lower
order terms of the same form:
$$
\|\delta \dot{x}\|_r+\|\delta x\|_r\leq K(\|\dot{W}\|_r
+\|{W}\|_r+\|| x\||_{r+2,\infty}\|W\|_{0}\big), \qquad \|F\|_r\leq
K(\|\delta \Phi\|_r +\|| x\||_{r+1,\infty}\|\delta
\Phi\|_{0}\big), \tag 12.8
$$
where we used the interpolation inequality in Lemma {5.}7. By
(12.6)  $\|W\|_0\leq \|W\|_1\leq K \int_0^t\|F\|_1\, d\tau\leq
K\int_0^t\|\delta \Phi\|_1\, d\tau$ and by interpolation
$$
\||x\||_{r+r_0+4-s,\infty} \||x\||_{s+1,\infty} \leq
 K\||x\||_{r+r_0+3,\infty},\qquad 1\leq s\leq r \tag 12.9
$$
Hence (12.3) follows.
 \qed\enddemo

  If we also use Sobolev's lemma, $r_0=[n/2]+1$,
and estimate the integrals by the $L^\infty$ norms and use
interpolation and (12.2) we get
$$
\||\delta x\||_{r,\infty}\leq K
 \big( \|\delta \Phi\|_{r+r_0,\infty}+ \|| x\||_{r+2r_0+4,\infty}
 \|\delta \Phi\||_{0,\infty}\big),\qquad r\geq 0\tag 12.10
 $$
 Furthermore we want to turn it into estimates for
 $$
 \tilde{\Phi}(u)=\Phi(u+x_0)-\Phi(x_0)\tag 12.11
 $$
 see (2.23). Then $\tilde{\Phi}^\prime(u)=\Phi^\prime(u+x_0)$.
Let $\psi(u)$ denote the right inverse of
$\tilde{\Phi}^\prime(u)$. Then since $\||x\||_{r,\infty}\leq \||
u\||_{r,\infty}+\||x_0\||_{r,\infty}$ we can again include
 $\||x_0\||_{r,\infty}$ in the constants. Hence we have proven
 that

 \proclaim{Theorem {12.}2} Suppose that $T>0$ is so small that the
conditions in Lemma {11.}3 hold. Suppose also that $x=u+x_0$ and
$\tilde{g}$ are smooth in $[0,T]\times\overline{\Omega}$ and that
$\tilde{g} $ and $u$ vanish to infinite order as $t\to 0$.
 Then there are constants $K$, depending on the approximate solution
 $(x_0,h_0)$, on $(c_0,c_1)$ and on $r$, such that
 the linearized operator $\tilde{\Phi}^\prime(u)$ has a right
inverse $\psi(u)$ satisfying
$$
\||\psi(u) \tilde{g}\||_{r,\infty}\leq K
 \big( \||\tilde{g}\||_{r+r_0,\infty}+ \|| u\||_{r+2r_0+4,\infty}
 \|\tilde{g}\||_{0,\infty}\big),\qquad r\geq 0\tag 12.12
 $$
 if
 $$
 \|| u\||_{r_0+4,\infty}\leq 1 \tag 12.13
 $$
\endproclaim

 \head 13. Tame estimate for the second variational
derivative.\endhead
 We now first want to show that the Euler map $\Phi(x)$ given
 by (2.20)-(2.21) is $C^k$, i.e. that $\Phi(x)$ depends
 smoothly on $k$ parameters if $x$ does.
 To be more precise, with $B^k=\{r\in \bold{R}^k;\, |r|\leq 1\}$, we want to show that
 $\Phi(x)-\Phi(x_0)\in C^\infty(B^k,C_{00}^\infty)$ if
 $x-x_0\in C^\infty(B^k, C_{00}^\infty )$, where $x_0$ is the
 approximate solution satisfying (2.18) and $C_{00}^\infty$ is given by
 (2.19). That $D_t^2$ and $\pa_i$ in (2.20) depends smoothly on
 parameters is obvious so we only need to prove that $h=\Psi(x)$,
 given by (2.21) does. Subtracting off the approximate solution
 $h_0$
 of (2.21) we get (10.2), where the right hand side is in
 $C^k(B^k, C_{00}^\infty)$ if $x-x_0\in C^k(B^k, C_{00}^\infty)$.
 Hence it follows from Lemma {7.}6 that
 $h-h_0\in C^k(B^k, C_{00}^\infty)$.

We must now also obtain tame estimates for the second variational
derivative. If $x$ depends smoothly on the parameter $r$ then the
variational derivative of $x$ is $\delta x(t,y)=\pa x(r,t,y)/\pa
r\big|_{r=0}$, we can e.g. take $x=x(t,y)+r\, \delta x(t,y)$. The
first variational derivative $\Phi^\prime(x)$ of the Euler map is
given by
$$
\Phi^\prime(x)\delta x_i=\delta\Phi(x)_i=\pa \Phi(x)_i/\pa
r\big|_{r=0}=D_t^2 \delta x_i -\pa_k h\,\, \pa_i\delta x^k
+\pa_i\, h^{\,\prime}(\delta x), \tag 13.1
$$
where $\delta h=h^{\,\prime}(\delta x)=\Psi^\prime(x)\delta x$
satisfies
$$
D_t^2(c\delta h)-\triangle \delta h=-\delta\triangle p -\pa_k p\,
\triangle \delta x^k\! -2(\pa_i\pa_k p)\pa^i\delta x^k,\qquad
\delta \triangle h=2\pa_k V^i\, \pa_i\delta x^l\,\pa_l V^k
\!-2\pa_k V^i\,\pa_i \delta v^k\tag 13.2
$$
and $\delta p\big|_{\pa\Omega}=0$, where $\delta v=D_t\delta x$.

Now, let $x$ depend smoothly on two parameters $r$ and $s$, such
that $\pa^2 x/\pa r\pa s=0$,  and also set $\epsilon x=\pa x/\pa
s\big|_{s=0} $, e.g. $x=x(t,y)+r\delta x(t,y)+s\,\epsilon
\,x(t,y)$. Then the second variational derivative is given by
$$
 \Phi^{\prime\prime}(x)(\delta x,\epsilon x)_i=
 \epsilon \delta \Phi(x)_i= \pa\big( \pa \Phi_i(x)/\pa r\big|_{r=0}\big)/\pa s\big|_{s=0}
 .\tag 13.3
$$

We have \proclaim{Lemma {13.}1}
$$
\Phi^{\prime\prime}(\epsilon x,\delta x)_i= \pa_k h\,\, \big(\pa_i
\epsilon x^l\,\pa_l \delta x^k +\pa_i \delta x^l\,\pa_l \epsilon
x^k\big) -\pa_k h^\prime(\epsilon x)\,\, \pa_i \delta x^k -\pa_k
h^\prime(\delta x)\,\, \pa_i \epsilon x^k +\pa_i\,
h^{\prime\prime}(\epsilon x,\delta x)\tag 13.4
$$
where $h^\prime=\Psi^\prime(x)$ and
$h^{\prime\prime}=\Psi^{\prime\prime}(x)$.
\endproclaim
The estimates for $h^\prime=\Psi^\prime(x)$ and
$h^{\prime\prime}=\Psi^{\prime\prime}(x)$ must also be obtained:

\proclaim{Lemma {13.}2} Let $h=\Psi(x)$ and let $\delta
h=h^{\,\prime}(\delta x)=\Psi^\prime(x)\delta x$ be the
variational derivative. We have
$$
 \||\delta h\||_{r,\infty}\leq K\big( \||\delta
x\||_{r+r_0+1,\infty} +\||x\||_{r+2r_0+3,\infty}\||\delta
x\||_{0,\infty}\big)\tag 13.5
$$
and with $h^{\prime\prime}(\delta x,\epsilon
x)=\Psi^{\prime\prime}(x)(\delta x,\varepsilon x)$ the second
variational derivative, we have
$$\multline
\|h^{\prime\prime}(\delta x,\epsilon x)\|_{r,\infty}\\ \leq
 K\big( \||\delta x\||_{r+3r_0+6,\infty} \||\epsilon x\||_{0,\infty}
 +\||\epsilon x\||_{r+3r_0+6,\infty}\||\delta x\||_{0,\infty}
 +\||x\||_{r+3r_0+6,\infty} \||\delta x\||_{0,\infty}
 \||\epsilon x\||_{0,\infty}\big)
 \endmultline
 \tag 13.6
$$
\endproclaim
\demo{Proof of Lemma {13.}2} The proof is similar to the proof of
Lemma {10.}1. We have
 $$
 D_t^2 (ch)-\triangle h=(\pa_i V^j)(\pa_j V^i), \qquad
 \triangle h=\kappa^{-1}\pa_a\big( \kappa g^{ab}\pa_b h\big)\tag
 13.7
 $$
It follows that
$$
 D_t^2 (c\delta h)-\triangle \delta h=
 {\delta}\big((\pa_i V^j)(\pa_j V^i)\big)
 +\kappa^{-1}\pa_a\big(\delta( \kappa g^{ab})\pa_b h\big)
 -\div\delta x\, \kappa^{-1}\pa_a\big(\kappa g^{ab}\pa_b
 h\big),\tag 13.8
 $$
since $\delta \kappa=\kappa \div\delta x$, see \cite{L1}. Using
the estimate for $h$ in terms of $x$ in Lemma
 {10.}1 and Corollary {7.}5 gives if ${f}$ denotes the right hand
 side of (13.8)
 $$
\||{f}\||_{r,\infty}\leq K\big( \||\delta x\||_{r+2,\infty}
+\||x\||_{r+r_0+3,\infty}\||\delta x\||_{0,\infty}\big)\tag 13.9
 $$
 and hence by Corollary {7.}5
$$
\||\delta h\||_{r,\infty}\leq K\big( \||\delta
x\||_{r+r_0+1,\infty} +\||x\||_{r+2r_0+3,\infty}\||\delta
x\||_{0,\infty}\big)\tag 13.10
$$
In the proof we use the interpolation inequalities in Lemma {5.}7
and the fact that $\||x\||_{r_0+4,\infty}\leq K$.

 To calculate the second variation we apply  ${\epsilon}$ to this,
 where we assumed that $\epsilon\delta x=0$. Note that
 $\epsilon\div\delta x=-(\pa_i \epsilon x^k)(\pa_k \delta x^i)$.
$$\multline
 D_t^2 (c\epsilon\delta h)-\triangle\epsilon \delta h=
 \epsilon{\delta}\big((\pa_i V^j)(\pa_j V^i)\big)\\
 -\div\epsilon x\, \kappa^{-1}\pa_a\big(\delta( \kappa g^{ab})\pa_b h\big)
 -\div\delta x\, \kappa^{-1}\pa_a\big( \epsilon(\kappa g^{ab})\pa_b
 h\big)\\
 +\big(\div\epsilon x\, \div\delta x\,
 +(\pa_i \epsilon x^k)\pa_k \delta x^i)\big)
 \kappa^{-1}\pa_a\big(\kappa g^{ab}\pa_b h\big)
 -\div\delta x\, \kappa^{-1} \pa_a\big(\epsilon (\kappa g^{ab}) \pa_b h\big)
 +\kappa^{-1}\pa_a\big( \epsilon\delta (\kappa g^{ab})\, \pa_b h\big) \\
 +\kappa^{-1}\pa_a\big(\delta( \kappa g^{ab})\pa_b \epsilon h\big)
 -\div\delta x\, \kappa^{-1}\pa_a\big(\kappa g^{ab}\pa_b \epsilon h\big)
 +\kappa^{-1}\pa_a\big(\epsilon( \kappa g^{ab})\pa_b \delta h\big)
 -\div\epsilon x\, \kappa^{-1}\pa_a\big(\kappa g^{ab}\pa_b \delta h\big)
 \endmultline\tag 13.11
 $$
 Here $\delta g^{ab}=-g^{ac} g^{bd}\delta
 g_{cd}$, $\delta g_{ab} =\delta_{ij}(\pa_a \delta x^i)(\pa_b x^j) +
 \delta_{ij}(\pa_a x^i)(\pa_b \delta x^j)$
 so since $\epsilon\delta x=0$, we have
 $\epsilon\delta g_{ab}=2\delta_{ij}(\pa_a \delta x^i)(\pa_b \epsilon
 x^j)$.
 The first terms on the right of (13.11) gives raise to term of
 the form
 $$
 (\pa \epsilon v)(\pa\delta v),\quad
 (\pa \epsilon v )(\pa \delta x) (\pa v),\quad
 (\pa \delta v )(\pa \epsilon x) (\pa v),\quad
 (\pa \epsilon x)(\pa \delta x)(\pa v)(\pa v)\tag 13.12
 $$
 multiplied by powers of $\pa y/\pa x$ or $\pa x/\pa y$.
 If ${f}_1$ denotes any of these terms then
 $$
 \||{f}_1\||_{r,\infty}\leq K\big( \||\delta x\||_{r+4,\infty}
 \||\epsilon x\||_{0,\infty} +\||\epsilon x\||_{r+4,\infty}
 \||\delta x\||_{0,\infty} +\||x\||_{r+5,\infty}\||\delta
 x\||_{0,\infty}\||\epsilon x\||_{0,\infty} \big)\tag 13.13
 $$
The terms on the second and third row in (13.11) gives raise to
terms of the form
$$
(\pa^2 \delta x)(\pa \epsilon x)(\pa h), \quad
 (\pa \delta x)(\pa^2 \epsilon x)(\pa h),
 \quad (\pa \delta x)(\pa \epsilon x) (\pa^2 h) ,
 \quad (\pa \delta x)(\pa \epsilon x) (\pa h)(\pa^2 x)\tag 13.14
$$
multiplied by powers of $\pa y/\pa x$ or $\pa x/\pa y$. These can
be estimated by (13.13) plus
$$
K\||h\||_{r+4,\infty}\||\delta x\||_{0,\infty}\||\epsilon
x\||_{0,\infty} \leq K\||x\||_{r+r_0+5,\infty}\||\delta
x\||_{0,\infty}\||\epsilon x\||_{0,\infty}\tag 13.15
$$
by Lemma {10.}1.
 The terms on the last row in
(13.11) gives raise to terms of the form
$$
(\pa^2 \delta x) (\pa \epsilon h),\quad
 (\pa^2 \epsilon x)(\pa \delta h),\quad
 (\pa \delta x)(\pa^2 \epsilon h),\quad
 (\pa \epsilon x)(\pa^2 \delta h),\quad
 (\pa \epsilon x)(\pa \delta h)(\pa^2 x),\quad
 (\pa\delta x)(\pa \epsilon h) (\pa^2 x) \tag 13.14
$$
multiplied by powers of $\pa y/\pa x$ or $\pa x/\pa y$. These can
be estimated by
$$\multline
 K\big( \||\delta x\||_{r+3,\infty}
 \||\epsilon h\||_{0,\infty} +\||\epsilon x\||_{r+3,\infty}
 \||\delta h\||_{0,\infty}
 +\||\delta h\||_{r+3,\infty}
 \||\epsilon x\||_{0,\infty} +\||\epsilon h\||_{r+3,\infty}
 \||\delta x\||_{0,\infty}\big)\\
 +K\||x\||_{r+4,\infty}(\||\delta
 x\||_{0,\infty}\||\epsilon h\||_{0,\infty}
 +\||\delta h\||_{0,\infty}\||\epsilon x\||_{0,\infty} \big)
 \endmultline \tag 13.15
 $$
 If we also use (13.5) we see that this can be estimated by
$$\multline
 K\big( \||\delta x\||_{r+3,\infty}
 \||\epsilon x\||_{r_0+1,\infty} +\||\epsilon x\||_{r+3,\infty}
 \||\delta x\||_{r_0+1,\infty}
 +\||\delta x\||_{r+r_0+4,\infty}
 \||\epsilon x\||_{0,\infty} +\||\epsilon x\||_{r+r_0+4,\infty}
 \||\delta x\||_{0,\infty}\big)\\
 +K\||x\||_{r+4,\infty}(\||\delta
 x\||_{0,\infty}\||\epsilon x\||_{r_0+1,\infty}
 +\||\delta x\||_{r_0+1,\infty}\||\epsilon x\||_{0,\infty} \big)\\
 +K\|| x\||_{2r_0+3,\infty}\big( \||\delta x\||_{r+3,\infty}
 \||\epsilon x\||_{0,\infty} +\||\epsilon x\||_{r+3,\infty}
 \||\delta x\||_{0,\infty}\big)\\
 +K\||x\||_{r+2r_0+6,\infty} \big(\||\delta x\||_{0,\infty}
 \||\epsilon x\||_{0,\infty} +\||\epsilon x\||_{0,\infty}
 \||\delta x\||_{0,\infty}\big)\\
 +K\||x\||_{r+4,\infty}\||x\||_{2r_0+3,\infty}(\||\delta
 x\||_{0,\infty}\||\epsilon x\||_{0,\infty}
 +\||\delta x\||_{0,\infty}\||\epsilon x\||_{0,\infty} \big)
 \endmultline \tag 13.16
 $$
Using interpolation again this is bounded by
$$
K\big(\|| \delta x\||_{r+r_0+4,\infty}\||\epsilon x\||_{0,\infty}
+\|| \epsilon x\||_{r+r_0+4,\infty}\||\delta x\||_{0,\infty}
+\||x\||_{r+2r_0+6,\infty}\||\delta x\||_{0,\infty}\||\epsilon
x\||_{0,\infty} \big)\tag 13.17
$$
 It follows that
 $$\multline
 \||D_t^2 (c\epsilon\delta h)-\triangle \epsilon\delta h
 \||_{r,\infty}\\
 \leq K\big(\|| \delta x\||_{r+r_0+4,\infty}\||\epsilon x\||_{0,\infty}
+\|| \epsilon x\||_{r+r_0+4,\infty}\||\delta x\||_{0,\infty}
+\||x\||_{r+2r_0+6,\infty}\||\delta x\||_{0,\infty}\||\epsilon
x\||_{0,\infty} \big)
\endmultline \tag 13.18
 $$
Hence by Corollary {7.}5 :
$$\multline
\||\epsilon\delta h\||_{r,\infty}\\
\leq K\big(\|| \delta x\||_{r+2r_0+3,\infty}\||\epsilon
x\||_{0,\infty} +\|| \epsilon x\||_{r+2r_0+3,\infty}\||\delta
x\||_{0,\infty} +\||x\||_{r+3r_0+5,\infty}\||\delta
x\||_{0,\infty}\||\epsilon x\||_{0,\infty} \big)
\endmultline \tag 13.19
$$
 \enddemo

\proclaim{Theorem {13.}3} Suppose that $T>0$ is so small that the
assumptions in Lemma {11.}3 hold. Then there are constants $K$
depending on the approximate solution $(x_0,h_0)$, on $(c_0,c_1)$
and on $r$ such that
$$
\multline
 \||\Phi^{\prime\prime}(u+x_0)(\epsilon x,\delta
x)\||_{r,\infty}\\
 \leq K\big( \||\delta x\||_{r+2r_0+4,\infty} \||\epsilon x\||_{0,\infty}
 +\||\epsilon x\||_{r+2r_0+4,\infty}\||\delta x\||_{0,\infty}
 +\||u\||_{r+3r_0+6,\infty} \||\delta x\||_{0,\infty}
 \||\epsilon x\||_{0,\infty}\big)
 \endmultline
 \tag 13.20
$$
if
$$
 \|| u\||_{r_0+4,\infty}\leq 1. \tag 13.21
 $$
\endproclaim

\head{14. The smoothing operators.}\endhead We will work in
H\"older spaces since the standard proofs of the Nash-Moser
theorem uses H\"older spaces. The H\"older norms for functions
defined on a compact convex set $B\in\bold{R}^{1+n}$ are given by,
if $k<a\leq k+1$, where $k\geq 0$ is an integer,
$$
\||u\||_{a,\infty}=\|| u\||_{H^a}=\sup_{(t,y),(s,z)\in B}\sum_{|\alpha|=k}
\frac{ |\pa^\alpha u(t,y)-\pa^\alpha
u(s,z)|}{|(t,y)-(s,z)|^{a-k}}+ \sup_{(t,y)\in B}|u(t,y)|\tag 14.1
$$
and $\||u\||_{H^0}=\sup_{(t,y)\in B} |u(t,y)|$.  Since we use the
same notation for the $C^k$ norms, $\||u\||_{k,\infty}$
we will indicate the difference by simply using letters $a,b,c,d,e,f$ etc for
the H\"older norms and $i,j,k,l, ..$ for the $C^k$ norms. Since a
Lipschitz continuous function is differentiable almost everywhere
and the norm of the derivative at these points is bounded by the
Lipschitz constant, we conclude that for integer values this is
the same if the $L^\infty$ norm of $\pa^\alpha u$ for
$|\alpha|\leq k$, and furthermore, since all our functions are
smooth it is the same as the supremum norm. Our tame estimates for
the inverse of the linearized operator and the second variational
derivative are only for $C^k$ norms with integer exponents.
However, since $\||u\||_{k,\infty}\leq C\||u\||_{a,\infty}\leq C\||u\||_{{k+1,\infty}}$,
if $k\leq a\leq k+1$, see (14.2),  they also hold for non integer
values with a loss of of one more derivative.

In \cite{L3} we used smoothing only in the space directions but
here we will use smoothing also in the time direction. Therefore
we define the H\"older space time norms as above.

 They satisfy
$$
\||u\||_{a,\infty}\leq C\||u\||_{b,\infty} ,\qquad a\leq b\tag 14.2
$$
and they also satisfy the interpolation inequality
$$
\|| u\||_{c,\infty}\leq C\||u\||_{a,\infty}^\lambda \||u\||_{b,\infty}^{1-\lambda}\tag 14.3
$$
where $0\leq a\leq c\leq b$, $0\leq \lambda\leq 1$ and $c=\lambda
a+(1-\lambda)b$. Furthermore, the H\"older spaces are rings:
$$
\|| uv\||_{a,\infty}\leq C(\||u\||_{a,\infty}\|| v\||_{0,\infty}+\|| u\||_{0,\infty}
\||v\||_{a,\infty})\tag
14.4
$$

 For the Nash-Moser technique, apart from
tame estimates one also needs smoothing operator $S_\theta$
 that satisfy the following properties with respect to the H\"older norms:
 \proclaim{Lemma {14.}1} Let $\||u\||_a$ denote the H\"older norms
 in (14.1) with $B=[0,T]\times\overline{\Omega}$, where $T\leq 1$. Let
 $C_{00}^\infty=C_{00}^\infty\big([0,T]\times\overline{\Omega}\big)$ be as in (14.11).
 Then there is a family of
 smoothing operators $S_\theta: C_{00}^\infty\to C_{00}^\infty$,
 $1\leq \theta<\infty$ such that
$$\align
\||S_\theta u\||_{a,\infty}&\leq C\|| u\||_{b,\infty},\qquad a\leq b\tag 14.5\\
\||S_\theta u\||_{a,\infty}&\leq C \theta^{a-b} \||u\||_{b,\infty},\qquad a\geq b \tag 14.6\\
\||(I-S_\theta) u\||_{a} &\leq C \theta^{a-b}\||u\||_{b,\infty},\quad a\leq b\tag 14.7\\
\||(S_{2\theta}-S_\theta)u\||_{a,\infty} &\leq C\theta^{a-b}\||u\||_{b,\infty},\qquad a\geq 0  \tag 14.8
\endalign
$$
where the constants $C$ only depend on the dimension and an upper
bound for $a$ and $b$.
\endproclaim
 The last property, (14.8) follows from (14.6) for $a\geq b$ and from (14.7) for $a\leq b$.
Alternatively, it follows from the following stronger property
$$
\||\tfrac{d}{d\theta} S_\theta u\||_{a,\infty}\leq C\theta^{a-b-1}\||u\||_{b,\infty},\qquad a\geq 0.\tag 14.9
$$

For functions supported in the interior of a compact set $K$ there
there are smoothing operators, see \cite{H1}, that satisfy the above properties
(14.5)-(14.9), with respect to the H\"older norms. These are constructed as
follows. Let the Fourier transform $\hat{\phi}\in C_0^\infty$ be $1$ in a neighborhood
of the origin and set $\phi_\theta(z)=\theta^{1+n}\phi(\theta z)$, and
$S_\theta u=\chi \phi_\theta *u$, where $\chi\in C_0^\infty$ is $1$ on a neighborhood
of $K$. However we have
functions defined on the compact set
$[0,T]\times\overline{\Omega}$ that do not have compact support in
$\Omega$. Therefore we need to extend these functions to have
compact support in some larger set, without increasing the
H\"older norms with more than with a multiplicative constant.
There is a standard extension operator in \cite{S} that turns to
have these properties, see Lemma {14.}2 below.

First however, we note that we will only apply the smoothing
operators to functions that vanish to all orders as $t\to 0$.
Hence we can extend these functions to be $0$ for $t\leq 0$
without changing the H\"older norm. Then we extend the functions
defined for $t\leq T$ to functions supported in $[0,2]$, using the
extension in Lemma {14.}2, for $y\in \overline{\Omega}$ fixed,
noting that H\"older continuity in $(t,y)$ follows from
differentiability and H\"older continuity in each direction using
the triangle inequality and the linearity of the extension
operator in Lemma {14.}2. Then we want to extend the functions
defined in $\overline{\Omega}=\{y; |y|\leq 1\}$ to functions
supported in $\{y; |y|\leq 2\}$. In order to do this we first
remove a region around the origin and  introduce polar coordinates
$r\geq 0$ and $\omega\in S^{n-1}$, which is a nonsingular change
of variables away from the origin. Then we use the extension
operator in Lemma {14.}2 to, for fixed $t$ and $\omega$ extend in
the radial direction from function defined for $r\leq 1$ to
functions supported in $r\leq 2$. By the remark above, H\"older
continuity in $(t,r,\omega)$ follows from H\"older continuity in
each direction.

Doing the extensions above we hence obtain an extension
$\tilde{u}$ of $u$ defined in $[0,T]\times \overline{\Omega}$ such
that
$$
\|| \tilde{u}\||_{a,\infty}\leq C \|| u\||_{a,\infty} , \qquad a\geq 0\qquad  \text{supp}
(\tilde{u})\in \{ (t,y);\, 0\leq t\leq 2,\, |y|\leq 2\}\tag
14.10
$$
for $u$ in
$$
C_{00}^\infty=C_{00}^\infty\big([0,T]\times\overline{\Omega}\big)
=\{ u\in C^\infty([0,T]\times \overline{\Omega}), \, D_t^k
u\big|_{t=0}=0,\,  k\geq 0\}\tag 14.11
$$
We note that, in fact the constant in (14.10) is independent of
$T$.

Once we have the extension operator we can use the smoothing
operators in \cite{H1,H2}, defined for compactly supported
functions, applied to the extension of our function. Let us call
the smoothing operators defined in \cite{H1,H2}
$\tilde{S}_\theta$. These satisfy the properties (14.5)-(14.9). By
(14.10) the smoothing operators $\hat{S}_\theta u$ given by the
restriction of $\tilde{S}_\theta \tilde{u}$ to
$[0,T]\times\overline{\Omega}$ then also satisfy the properties
(14.5)-(14.9) if $u$ is in (14.11). However, $\hat{S}_\theta u$ is
not in (14.11) anymore. In our estimates we will
 only apply the smoothing operators to functions that vanish to
 all orders as $t\to 0$ and in our estimates we need also
 $S_\theta u$ to vanish to all orders as $t\to 0$.
 We therefore have to modify our smoothing operators
 so that this is true. Let $\chi(t)\in C^\infty$ be a function
 such that $\chi(t)=0$, when $t\leq 0$ and $\chi(t)=1$, when
 $t\geq 1$, and let $\chi_\theta(t)=\chi(\theta t)$. Then
 $$
 S_\theta u=\chi_\theta \tilde{S}_\theta
 (\tilde{u})\big|_{[0,T]\times \overline{\Omega}} \tag 14.12
 $$
 is in (14.11) and we claim that for functions $u$ in (14.11),
 (14.5)-(14.9) hold. This follows from Lemma {14.}3 below, since
 the smoothing operators defined in \cite{H1,H2} are convolution
 operators of the form in Lemma {14.}3.

 \proclaim{Lemma {14.}2} There is a linear extension
operator ${\Cal Ext }: H^a((-\infty,0])\to  H^a((-\infty,+\infty))$,
where $H^a$ are the H\"older spaces,
 such that ${\Cal Ext}(f)=f$ when $r\leq 0$, and
$$
\|{\Cal Ext }( f)\|_{a,\infty}\leq C\|f\|_{a,\infty}\tag 14.13
$$
Here $C$ is bounded when $a$ is bounded.
Furthermore, if $r\geq -c$ in the support of $f$ , where $c>0$, then
$r\leq c$  in the support of ${\Cal Ext} (f)$.
\endproclaim
\demo{Proof} Let ${\Cal Ext} (f)(r)=\tilde{f}(r)$, where $\tilde{f}(r)=f(r)$,
 when $r\leq 0$, and
$$
\tilde{f}(r)=\int_1^\infty f(r-2\lambda r)\,
\psi_1(\lambda)\, d\lambda, \qquad r>0\tag 14.14
$$
where $\psi_1$ is a continuous function on $[1,\infty)$, such that
$$
\int_1^\infty \psi_1(\lambda)\, d\lambda =1, \qquad \int_1^\infty
\lambda^k \psi_1(\lambda)\, d\lambda =0,\quad k>0, \qquad
|\psi_1(\lambda)|\leq C_N (1+\lambda)^{-N}, \quad N\geq 0\tag
14.15
$$
The existence of such a function was proved in \cite{S} where the
extension operator was also introduced. In \cite{S} it was proven
that this operator is continuous on the Sobolev spaces but it was
not proven there that it is continuous on the H\"older spaces so
we must prove this.

First we note that if $f\in C^k$ then the extension is in $C^k$.
In fact
$$
\tilde{f}^{(j)}(r)=\int_1^\infty f^{(j)} (r-2\lambda
r)(1-2\lambda)^j \, \psi_1(\lambda)\, d\lambda, \qquad r>0
1\tag 14.16
$$
From the continuity of $\pa_r^j f$ and (14.14)-(14.15) it follows
that $\lim_{r\to +0}\tilde{f}^{(j)}(r)=f^{(j)}(0)$,
that $\tilde{f}$ is in $C^k$, and that $\|\tilde{f}\|_{k,\infty}
\leq C_k \|f\|_{k,\infty}$, if $k$ is an  integer.

Suppose now that $k<a\leq k+1$ where $k$ is an integer. We must
now estimate
$$
\sup_{r,\rho} \frac{|\tilde{f}^{(k)}(r)-\tilde{f}^{(k)}(\rho)|}
{|r-\rho|^{a-k}} \tag 14.17
$$
by $C_k \|\pa_r^k f\|_{a-k,\infty}$. If $r\leq 0$ and $\rho\leq 0$ there
is nothing to prove. Also if $r<0<\rho$ or $\rho<0<r$, then
$|r-\rho|\geq |\rho| $ and $|r-\rho|\geq |r|$ so in this case, we
can reduce it to two estimates with either $r=0$ or $\rho=0$. Also
it is symmetric in $r$ and $\rho$ so it only remains to prove the
assertion when $r>\rho\geq  0$.
It follows from the H\"older continuity of $f^{(k)}$ and the last
estimate in (14.15) that for $r,\rho\geq 0$,
$$
\Big|\int_1^\infty\big( f^{(k)} (r-2\lambda r)- f^{(k)}
(\rho-2\lambda \rho)\big) (1-2\lambda)^k \, \psi_1(\lambda)\,
d\lambda \Big|\leq C_k  \| f^{(k)}\|_{a-k,\infty} |r-\rho|^{a-k} \tag
14.18
$$
which proves the lemma.\qed\enddemo

\proclaim{Lemma {14.}3} Let $\chi(t)\in C^\infty$ be a function
 such that $\chi(t)=0$, when $t\leq 0$ and $\chi(t)=1$, when
 $t\geq 1$, and let $\chi_\theta(t)=\chi(\theta t)$. Let
the Fourier transform $\hat{\phi}\in C_0^\infty$ be $1$ in a neighborhood
of the origin and let $\chi_1\in C_0^\infty$ be $1$ on a neighborhood of
$\{(t,y);\, 0\leq t\leq 2,\, |y|\leq 2\}$. Set
 $S_\theta =\chi_1 \phi_\theta* u$, where
$\phi_\theta((t,y))=\phi(\theta (t,y))/\theta^{1+n}$.
 Then
 $$
\|| (1-\chi_\theta)S_\theta u\||_{a,\infty}\leq C|\theta|^{a-b}
\||u\||_{b,\infty},\tag 14.19
$$
if $u$ is smooth and vanishes for $t\leq 0$ and for $t\geq 2$.
\endproclaim
\demo{Proof} First we note that by interpolation it suffices to prove the estimate
for $a=k$ an integer.
Since $u$ vanishes to infinite order as $t\to 0$ we have if $k<b\leq k+1$
$$
|u(t,y)|=\Big| \int_0^t (\pa_t^k u)(s,y)\, \frac{(t-s)^{k-1}}{(k-1)!} \, ds\Big|
\leq  \int_0^t s^{b-k} \||u\||_{b,\infty} \, \frac{(t-s)^{k-1}}{(k-1)!} \, ds
\leq C_{b} t^b \||u\||_{b,\infty} .\tag 14.20
$$
Since $\phi$ is fast decaying we have if $|\alpha|=k$,
$$
\big| D^\alpha \big(   (1-\chi_\theta)\chi_1 \phi_\theta * u\big)\big|\leq
C_N \theta^{k-b}\||u\||_{b,\infty} \int\int_0^\infty \frac{ |s|^b  \, \, ds\, dy }
{(1+|\theta t- s|+|\theta x- y|)^N}\tag 14.21
$$
for any $N$.
Here the integral is uniformly bounded when $\theta t\leq C$ so the lemma follows.
 \qed\enddemo

\head 15. The Nash Moser Iteration.\endhead At this point, given
the results stated in sections 11-14, the problem is now reduced
to a completely standard application of the Nash-Moser technique.
One can just follow the steps of the proof of \cite{AG,H1,H2,K1}
replacing their norms with our norms. The main difference is that
we have a boundary, but we have constructed smoothing operators
that satisfy the required properties for the case with a boundary.
The proof of the Nash-Moser that we outline below is similar to the
one in \cite{L3}. The only difference is that now we also smooth
in time. We will follow the formulation from \cite{AG} which
however is similar to \cite{H1,H2}. The theorem in \cite{AG} is
stated in terms of H\"older norms, with a slightly different
definition of the H\"older norms for integer values. However, the
only properties that are used of the norms are the smoothing
properties and the interpolation property in section 14, which we
proved with the usual definition, i.e. the one used in \cite{H1}.

Let us also change notation and call $\tilde{\Phi}(u)$ from last
section $\Phi(u)$. For $k<a\leq k+1$, where $k\geq 0$ is an
integer, let
$$
\||u\||_{a,\infty}=\|| u\||_{H^a}=\sup_{(t,y),(s,z)\in
[0,T]\times\overline{\Omega}}\sum_{|\alpha|=k} \frac{ |\pa^\alpha
u(t,y)-\pa^\alpha u(s,z)|}{|(t,y)-(s,z)|^{a-k}}+ \sup_{(t,y)\in
[0,T]\times\overline{\Omega}}|u(t,y)|\tag 15.1
$$
and $\||u\||_{H^0}=\sup_{(t,y)\in [0,T]\times\overline{\Omega}}
|u(t,y)|$. The estimates we proved for the inverse of the
linearized operator and the second derivative of the operator
where in terms of $L^\infty$ norms, i.e. H\"older norms for
integer values. However, since $\||u\||_{k,\infty}\leq \||u\||_{a,\infty}\leq
\||u\||_{k+1}$, if $k\leq a\leq k+1$, it follows that they also
hold for non integer values with loss of an additional derivative.

(${\Cal H}_1$):  $\Phi$, is twice differentiable and satisfies
$$\multline
\||\Phi^{\prime\prime}(u)(v_1,v_2)\||_{a,\infty}\\
\leq C_a \Big( \||
v_1\||_{a+\mu,\infty} \|| v_2\||_{\mu^\prime,\infty}+
 \|| v_1\||_{\mu^\prime,\infty}
\||v_2\||_{a+\mu,\infty} +\|| u\||_{a+\mu,\infty} \||v_1\||_{\mu^\prime,\infty} \||v_2
\||_{\mu^\prime,\infty}\Big),
\endmultline \tag 15.2
$$
where $\mu=3r_0+8$, for $u,v_1,v_2\in C^\infty_{00}$, if
$$
\|| u\||_{\mu^\prime,\infty}\leq 1,\qquad \mu^\prime=r_0+4\tag 15.3
$$
where $K$ is the constant in Lemma {12.}1

(${\Cal H}_2$): If  $u\in C^\infty_{00}$ satisfies (15.3) then
there is a linear map
 $\psi(u)$ from  $C^\infty_{00}$
to $ C^\infty_{00}$ such that $\Phi^\prime(u)\psi(u) =Id$. It
satisfies
$$
\||\psi(u) g\||_{a,\infty}\leq C_a\big( \|| g\||_{a+\lambda,\infty} +\||g
\||_{0,\infty}\, \||u\||_{a+\lambda,\infty} \big), \tag 15.4
$$
where $\lambda =2r_0+6$.

\proclaim{Theorem {15.}1} Suppose that $\Phi$ satisfies (${\Cal
H}_1$),
 (${\Cal H}_2$) and $\Phi(0)=0$. Suppose that $\mu\geq \mu^\prime$
 and let
$\alpha>\lambda+\mu+\mu^\prime$, $\alpha\notin \Bbb{N}$. Then

i) There is neighborhood $W_\delta=\{ f\in C_{00}^\infty ; \, \||
f\||_{\alpha+\lambda,\infty}\leq \delta^2\}$, $\delta>0$, such that, for
$f\in W_\delta $, the equation
$$
\Phi(u)=f\tag 15.5
$$
has a solution $u=u(f)\in C_{00}^\infty$. Furthermore,
$$
\|| u(f)\||_{a,\infty}\leq C\||f\||_{\alpha+\lambda,\infty},\qquad a<\alpha\tag
15.6
$$
\endproclaim

In the proof, we construct a sequence $u_j\in C_{00}^\infty$
converging to $u$, that satisfy
 $\|| u_j\||_{\mu^\prime}\leq 1$ and $\||S_i u_i\||_{\mu^\prime}\leq
 1$, for all $j$, where $S_i$ is the smoothing operator in (15.7).
 The estimates (15.2) and (15.4) will only be used for convex
combinations of these and hence within the domain (15.3) for which
these estimates hold.

Following \cite{H1,H2,AG,K1,K2} we set
$$
u_{i+1}=u_i+\delta u_i,\qquad \delta u_i=\psi(S_i u_i) g_i,  \quad
u_0=0, \qquad S_i=S_{\theta_i}, \quad \theta_i=\theta_0 2^i,\quad
\theta_0\geq 1\tag 15.7
$$
and $g_i$ are to be defined so that $u_i$ formally converges to a
solution. We have
$$\multline
\Phi(u_{i+1})-\Phi(u_i)=\Phi^\prime(u_i)(u_{i+1}-u_i)+e^{\prime\prime}_i
=\Phi^\prime(u_i)\psi(S_i u_i) g_i +e^{\prime\prime}_i\\
=(\Phi^\prime(u_i)-\Phi^\prime(S_i u_i))\psi(S_i u_i)
g_i+g_i+e^{\prime\prime}_i =e^\prime_i+e^{\prime\prime}_i+g_i
\endmultline \tag 15.8
$$
where
$$\align
e_i^\prime&=(\Phi^\prime(u_i)-\Phi^\prime(S_i u_i))\delta u_i\tag 15.9 \\
e_i^{\prime\prime}&=\Phi(u_{i+1})-\Phi(u_i)-\Phi^\prime(u_i)\delta u_i \tag 15.10\\
e_i&=e^\prime_i+e^{\prime\prime}_i\tag 15.11
\endalign
$$
Therefore
$$
\Phi(u_{i+1})-\Phi(u_i)=e_i+g_i\tag 15.12
$$
and adding, we get
$$
\Phi(u_i)=\sum_{j=0}^i g_j+S_i E_i +e_i +(I-S_i)E_i, \qquad
E_i=\sum_{j=1}^{i-1} e_j\tag 15.13
$$
To ensure that $\Phi(u_i)\to f$ we must have
$$
\sum_{j=0}^i g_j+S_i E_i=S_i f\tag 15.14
$$
Thus
$$
g_0=S_0 f, \qquad g_i=(S_{i}-S_{i-1})(f-E_{i-1})-S_i e_{i-1}\tag
15.15
$$
and
$$
\Phi(u_i)=S_i f+e_i +(I-S_i) E_i\tag 15.16
$$
Given $u_0,u_1,...,u_i$ these determine $\delta u_0,\delta
u_1,...,\delta u_i$ which by (15.9)-(15.10) determine
$e_1,...,e_{i-1}$, which by (15.15) determine $g_i$. The new term
$u_{i+1}$ is the determined by (15.7).

\proclaim{Lemma {15.}2} Assume that $\|| u_i\||_{\mu^\prime,\infty}\leq
1$, $\||u_{i+1}\||_{\mu^\prime,\infty} \leq 1$, and $\|| S_i
u_i\||_{\mu^\prime,\infty}\leq 1$. Then
$$\multline
\|| e_i^\prime\||_{r,\infty}\leq C_r \big( \||(I-S_i)
u_i\||_{r+\mu,\infty}\||\delta u_i\||_{\mu^\prime,\infty}
+\||(I-S_i) u_i\||_{\mu^\prime,\infty}\||\delta u_i\||_{r+\mu}\big)\\
+C_r\||S_i u_i\||_{r+\mu,\infty} \||(I-S_i) u_i\||_{\mu^\prime,\infty} \||\delta
u_i\||_{\mu^\prime,\infty}\big)
\endmultline\tag 15.17
$$
and
$$
\|| e_i^{\prime\prime}\||_{r,\infty}\leq C_r \big( \||\delta
u_i\||_{r+\mu,\infty}\||\delta
u_i\||_{\mu^\prime,\infty}+\||u_i\||_{r+\mu,\infty}\||\delta
u_i\||_{\mu^\prime,\infty}^2\big) \tag 15.18
$$
\endproclaim
\demo{Proof} The proof of (15.17) makes use of
$$
(\Phi^\prime(u_{i})-\Phi^\prime(S_i u_i))\delta u_i = \int_0^1
\Phi^{\prime\prime}(S_i u_i+s(I-S_i)u_i)(u_i-S_i u_i ,\delta
u_i)\, ds \tag 15.19
$$
together with (15.2). The proof of (15.18) makes use of
$$
\Phi(u_{i+1})-\Phi(u_i)-\Phi^\prime(u_i)\delta u_i = \int_0^1
(1-s) \Phi^{\prime\prime}(u_i+s\delta u_i)(\delta u_i,\delta
u_i)\, ds \tag 15.20
$$
together with (15.2). \qed\enddemo

Let $\tilde{\alpha}>\alpha$ and
$\tilde{\alpha}-\mu>2\alpha-\mu-\mu^\prime$.
 Throughout the proof $C_a$ will stand for constants that
depend on $a$ but are independent of $i$ and $n$ in (15.21).

Our inductive assumption $(H_n)$ is,
$$
\|| \delta u_i\||_{a,\infty}\leq \delta \theta_i^{a-\alpha} , \qquad
0\leq a\leq \tilde{\alpha},\qquad i\leq n\tag 15.21
$$
If $n=0$ then $\delta u_0 =\psi(0) S_0f$, and if $a\leq
\tilde{\alpha}$, we have
 $\|| \delta u_0\||_{a,\infty}\leq C_{\tilde{\alpha}}
 \|| f\||_{\alpha+\lambda,\infty}\leq C_{\tilde{\alpha}} \delta^2$,
so it follows that (15.21) hold for $n=0$ if we choose $\delta$ so
small that $C_{\tilde{a}} \delta \leq
\theta_0^{\tilde{\alpha}-\alpha}$. We must now prove that $(H_n)$
implies $(H_{n+1})$ if $C_{\tilde{\alpha}}^\prime\delta \leq 1$,
where $C_{\tilde{\alpha}}^\prime$ is some
 constant that only depends on $\tilde{\alpha}$ but is
independent of $n$.

\proclaim{Lemma {15.}3} If (15.21) hold then with a constant
$C_{a}$ independent of $i\leq n$
 $$
 \sum_{j=0}^i\||\delta u_j\||_{a,\infty}\leq
C_a \delta
\big(\min(i,1/|\alpha-a|)+1)(\theta_i^{a-\alpha}+1\big), \qquad
0\leq a\leq \tilde{\alpha}\tag 15.22
$$
\endproclaim
\demo{Proof} Using (15.21) we get $\sum_{j=0}^i\||\delta
u_j\||_{a,\infty}\leq C_a\delta \sum_{j=0}^i 2^{j(a-\alpha)}$ and noting
that $\sum_{j=0}^i 2^{-sj}\leq C(\min{(1+1/s,i)}+1)$, if $s>0$,
(15.22) follows. \qed\enddemo

  It follows from
(15.22):
 \proclaim{Lemma {15.}4} If $(H_n)$, i.e. (15.21),  hold,
and $\tilde{\alpha}>\alpha$, then for $i\leq n+1$ we have
$$\align
\|| u_i\||_{a,\infty}&\leq C_a \delta
(\min(i,1/|\alpha-a|)+1)(\theta_i^{a-\alpha}+1), \qquad 0\leq
a\leq \tilde{\alpha}
\tag 15.23\\
\|| S_i u_i\||_{a,\infty}&\leq  C_a\delta
 (\min(i,1/|\alpha-a|)+1)(\theta_i^{a-\alpha}+1) ,
\qquad a\geq 0 \tag 15.24\\
\|| (I-S_i) u_i\||_{a,\infty}&\leq  C_a\delta \theta_i^{a-\alpha}, \qquad
0\leq a\leq \tilde{\alpha} \tag 15.25
\endalign
$$
where the constants are independent of $n$.
\endproclaim
\demo{Proof} The proof of (15.23) is just summing up the series
$u_{i+1}=\sum_{j=0}^i \delta u_j$, using Lemma {15.}3.  (15.24)
follows from (15.23) using (14.5) for $a\leq \tilde{\alpha}$ and
(14.6), with $b=\tilde{\alpha}$ for $a\geq \tilde{\alpha}$.
(15.25) follows from (14.7) with $b=\tilde{\alpha}$ and (15.23)
with $a=\tilde{\alpha}$. \qed\enddemo
 Since we have assumed that $\alpha>\mu^\prime$, we note
  that in particular, it follows that
 $$
 \||u_i\||_{\mu^\prime,\infty}\leq 1\quad \text{and}\quad
 \||S_i u_i\||_{\mu^\prime,\infty}\leq 1,\qquad
 \text{for}\quad i\leq n+1\quad\text{if}\quad
  C_{\mu^\prime} \delta\leq 1.
 \tag 15.26
 $$
  As a consequence of Lemma {15.}4 and Lemma {15.}2 we
get
 \proclaim{Lemma {15.}5} If $(H_n)$ is satisfied and
$\alpha>\mu^\prime$, then for $i\leq n$,
$$
\align \|| e_i^\prime\||_{a,\infty}\leq C_a\delta^2
\theta_i^{a-(2\alpha-\mu-\mu^\prime)},
\qquad 0\leq a\leq \tilde{\alpha}-\mu \tag 15.27\\
\|| e_i^{\prime\prime}\||_{a,\infty}\leq C_a\delta^2
\theta_i^{a-(2\alpha-\mu-\mu^\prime)}, \qquad 0\leq a\leq
\tilde{\alpha} -\mu\tag 15.28
\endalign
$$
where the constants are independent of $n$.
\endproclaim

As a consequence of Lemma {15.}5 and (14.8) we get
 \proclaim{Lemma
{15.}6} If $(H_n)$ is satisfied, then for $i\leq n+1$,
$$\align
\||\ S_i e_{i-1}\||_{a,\infty} &\leq  C_a\delta^2
\theta_i^{a-(2\alpha-\mu-\mu^\prime)},
\qquad a\geq 0\tag 15.29\\
\|| (S_{i}-S_{i-1}) f\||_{a,\infty} &\leq C_a\theta_i^{a-\beta} \||
f\||_{\beta,\infty},
\qquad a\geq 0 \tag 15.30\\
\||(I-S_i) f\||_{a,\infty} &\leq C_a\theta_i^{a-\beta} \|| f\||_{\beta,\infty},
\qquad 0\leq a\leq \beta  \tag 15.31
\endalign
$$
Furthermore, if $\tilde{\alpha}-\mu>(2\alpha-\mu-\mu^\prime)$:
$$\align
\|| (S_{i}-S_{i-1}) E_{i-1}\||_{a,\infty} &\leq  C_a\delta^2
\theta_i^{a-(2\alpha-\mu-\mu^\prime)},
\qquad a\geq 0\tag 15.32\\
\||(I-S_i) E_{i}\||_{a,\infty} &\leq  C_a\delta^2
\theta_i^{a-(2\alpha-\mu-\mu^\prime)}, \qquad  0\leq a\leq
\tilde{\alpha}-\mu \tag 15.33
\endalign
$$
Here the constants $C_a$ are independent of $n$.
\endproclaim
\demo{Proof} (15.29) follows from (15.27).  For $a\leq
\tilde{\alpha}-\mu$ we use (14.5) with $b=a$ and for $a\geq
\tilde{\alpha}-\mu$, we use (14.6) with $b=\tilde{\alpha}-\mu$.
(15.30) follows from (14.8) and (15.31) follows from (14.7).
 Now, $E_i=\sum_{j=0}^{i-1} e_j$ so by Lemma {15.}5
 $\|| E_i\||_{\tilde{\alpha}-\mu,\infty} \leq C_a\delta^2 \sum_{j=0}^{i-1}
 \theta_j^{\tilde{\alpha}-\mu-(2\alpha-\mu-\mu^\prime)}\leq
 C_a^\prime \delta^2
 \theta_i^{\tilde{\alpha}-\mu-(2\alpha-\mu-\mu^\prime)}$,
 since we assumed that the exponent is positive.
(15.32) follows from this and (14.8) with $b=\tilde{\alpha}-\mu$
and similarly (15.33) follows from (14.7) with
$b=\tilde{\alpha}-\mu$.\qed\enddemo

It follows that: \proclaim{Lemma {15.}7} If $(H_n)$ is satisfied,
$\tilde{\alpha}-\mu>(2\alpha-\mu-\mu^\prime)$, and
$\alpha>\mu^\prime$ then for $i\leq n+1$,
$$
\|| g_i\||_{a,\infty}\leq C_a\delta^2\theta_i^{a-(2\alpha-\mu-\mu^\prime)}
+C_a \theta_i^{a-\beta}\||f\||_{\beta,\infty} ,\qquad a\geq 0. \tag 15.34
$$
\endproclaim

Using this lemma and (15.4) we get
 \proclaim{Lemma {15.}8} If
$(H_n)$ holds, $\tilde{\alpha}-\mu>(2\alpha-\mu-\mu^\prime)$,
$\alpha>\mu^\prime$, $\alpha>\lambda$ then, for  $i\leq n+1$, we
have
$$
\|| \delta u_i\||_{a,\infty}\leq C_a\delta^2
\theta_i^{a+\lambda-(2\alpha-\mu-\mu^\prime)} +C_a \|| f\||_{\beta,\infty}
\theta_i^{a+\lambda-\beta} , \qquad a\geq 0. \tag 15.35
$$
\endproclaim
\demo{Proof} Using Lemma {15.}7, (15.24), and (15.4) we get
$$\multline
\|| \delta u_i\||_{a,\infty}\leq C_a \big(\||g_i\||_{a+\lambda,\infty}+
\||g\||_{\lambda,\infty} \||S_i u_i\||_{a+\lambda,\infty}\big)\\
 \leq C_a\big(\delta^2
\theta_i^{a+\lambda-(2\alpha-\mu-\mu^\prime)} +\|| f\||_{\beta,\infty}\,
\theta_i^{a+\lambda-\beta}\big)\\
+ C_a\big(\delta^2 \theta_i^{\lambda-(2\alpha-\mu-\mu^\prime)}
+\|| f\||_{\beta,\infty} \,\theta_i^{\lambda-\beta}\big)
 \delta
 (\min(i,1/|\alpha-a-\lambda|)+1)(\theta_i^{a+\lambda-\alpha}+1)
\endmultline \tag 15.36
$$
Using that $\alpha>d$ we get (15.35). \qed\enddemo

If, we now pick $\beta=\alpha+\lambda$, and use the assumptions
that $\lambda+\alpha<2\alpha-\mu-\mu^\prime$, and $\||
f\||_{\alpha+\lambda,\infty}\leq \delta^2$,
 we get that for $i\leq n+1$,
$$
\|| \delta u_i\||_{a,\infty}\leq C_a \delta^2 \theta_i^{a-\alpha},\qquad
a\geq 0, \tag 15.37
$$
If we pick $\delta>0$ so small that
$$
C_{\tilde{\alpha}}\delta \leq 1,\tag 15.38
$$
the assumption $(H_{n+1})$ is proven.

The convergence of the $u_i$ is an immediate consequence of Lemma
{15.}3:
$$
\sum_{i=0}^\infty \||u_{i+1}-u_i\||_{a,\infty}\leq C_{a}\delta ,\qquad
a<\alpha \tag 15.39
$$
It follows from Lemma {15.}6 that
$$
\|| \Phi(u_i)-f\||_{a,\infty}\leq C_a \delta^2 \theta_i^{a-\alpha-\lambda}
\tag 15.40
$$
which tends to $0$, as $i\to \infty$, if $a<\alpha+\lambda$. In
particular it follows from (15.39) that $\||
u_j\||_{\mu^\prime,\infty}\leq 1$, if $\delta>0$ is chosen small enough.

It remains to prove $u\in C_{00}^\infty$. Note that in Lemma
{15.}8 we proved a better estimate than $(H_n)$. In fact if we let
$\gamma=2(\alpha-\mu)-(\alpha+\lambda)>0$ and
$\alpha^\prime=\alpha+\gamma$, then $\||
f\||_{\alpha^\prime+\lambda,\infty}\leq C$ implies that
$$
\|| \delta u_i\||_{a,\infty}\leq C_a \theta_i^{a-\alpha^\prime},\qquad
a\geq 0\tag 15.41
$$
Using this new estimate, in place of $(H_n)$, we can go back to
Lemma {15.}4-Lemma {15.}8 and replace $\alpha$ by $\alpha^\prime$.
Then it follows from Lemma {15.}8 that
$$
\|| \delta u_i\||_{a,\infty}\leq
C_a\theta_i^{a+\lambda-2(\alpha^\prime-\mu)} +C_a
\theta_i^{a+\lambda-\beta}\|| f\||_{\beta,\infty}\tag 15.42
$$
and if we now pick
$\gamma^\prime=2(\alpha^\prime-\mu)-(\lambda-\alpha^\prime)=2\gamma$
and
$\alpha^{\prime\prime}=\alpha^\prime+\gamma^\prime=\alpha+2\gamma$,
and use that $\|| f\||_{\alpha^\prime+\gamma^\prime,\infty}\leq C$ we see
that
$$
\|| \delta u_i\||_{a,\infty}\leq C_a
\theta_i^{a-\alpha^{\prime\prime}},\qquad a\geq 0\tag 15.43
$$
Since the gain $\gamma>0$ is constant, repeating this process
yields that (15.41) holds for any $\alpha^\prime$ and hence that
(15.39)-(15.40) hold for any $a\geq 0$, (without $\delta$ ). It
follows that $u_j$ is a Cauchy sequence in $C^k\big([0,T]\times
\overline{\Omega})$, for any $k$, and since $u_j\in C_{00}^\infty$
it follows that $u_j\to u\in C^\infty_{00}$ and $\Phi(u_j)\to f\in
C_{00}^\infty$, and since $\Phi$ is continuous it follows that
$\Phi(u)=f$. (15.6) follows from (15.37) with
$\delta^2=\||f\||_{\alpha+\lambda,\infty}$. This concludes the proof of
Theorem {15.}1.

 \head 16. Existence
 of initial data satisfying the compatibility conditions. \endhead
 In this section we show that there are initial data satisfying
 the compatibility conditions. We do not attempt to find the
 most general class of initial data that do so, the purpose
 is simply to show that our local existence theorem is not
 about the empty set. The set of initial data we construct
 can then easily be extended to a much larger set using
 essentially the same proof.

 Let us therefore start by making some simplifying assumptions.
 First we assume that $e(h)=h$.
We now want to find a formal power series solution in $t$ of the
system
$$
D_t^2 x_i=-\pa_i h,\qquad D_t^2 h-\triangle h=(\pa_i V^j)(\pa_j
V^i),\qquad v_i=D_t x_i \tag 16.1
$$
with initial data
$$
x\big|_{t=0}=f_0,\qquad D_t x\big|_{t=0}=v_0,\qquad h\big|_{t=0}
=h_0,\qquad D_t h\big|_{t=0}=h_1=-\div V_0\tag 16.2
$$

\proclaim{Lemma {16.}1} Let $h_l=D_t^l h\big|_{t=0}$, for $l\geq
0$ and let $V_0=v\big|_{t=0}$. If (16.1) and (16.2) holds then
 $$
h_{l+2}= \triangle h_l+F_{l+2}(h_{-1},...,h_{l-1}),\qquad l\geq
0\tag 16.3
$$
where $F_2=(\pa_i V_0^j)(\pa_j V_0^i)$ and in general
$F_{l+2}=F_{l+2}(h_{-1},...,h_{l-1})$ is a sum of the form
$$
F_{l+2}= C_{\alpha_1...\alpha_n}^{ l_1...l_n}
h_{l_1}^{\alpha_1}\cdot\cdot\cdot h_{l_n}^{\alpha_n} ,\qquad
h_l^\alpha =\pa_x^\alpha h_l ,\quad l\geq 0, \quad
h^{\alpha}_{-1}=\pa_x^{\alpha^\prime} V_0^{i_k}, \quad
\alpha=(\alpha^\prime,i_k)\tag 16.4
$$
with
$$
|\alpha_1|+l_1+...+|\alpha_n|+l_n=l+2, \quad n\geq 2, \quad -1\leq
l_i\leq l-1,\quad 1\leq |\alpha_i|+l_i\leq l+1\tag 16.5
$$
\endproclaim
\demo{Proof}
 If we us that $[D_t,\pa_i]=-(\pa_i V^k)\pa_k $ we obtain from differentiating the
 first equation in (16.1)
$$
\pa_i D_t^l h=-D_t^{l+1} v_i + a^{l\, k_1...k_n}(\pa_i
V_{k_1})\cdot\cdot\cdot(\pa V_{k_{n-1}}) V_{k_n}\tag 16.6
$$
where $V_{k}=D_t^k V$. Here the sum is over $k_1+...+k_n=l+2-n$,
$n\geq 2$, $k_n\geq 1$,  and  terms in the sum consists of
contractions over $n-1$ pairs of indices. From differentiating
$\div V=\pa_i V^i$,
$$
D_t^{l+1}\div V=\div D_t^{l+1} V+ d^{\,l\, k_1...k_n}(\pa
V_{k_1})\cdot\cdot\cdot\pa V_{k_{n}}\tag 16.7
$$
where the sums are over $k_1+...+k_n=l+2-n$, $n\geq 2$ and terms
in the sum consists of contractions over $n$ pairs of indices.

It follows from (16.1) that $D_t(D_t h+\div V)=0$. Hence
$$
\triangle D_t^{l} h-D_t^{l+2}h= e^{\,l\,k_1...k_n}(\pa
V_{k_1})\cdot\cdot\cdot\pa V_{k_{n}} + d^{\,l\,k_1...k_n}(\pa
V_{k_1})\cdot\cdot\cdot(\pa^2 V_{k_{n-1}}) V_{k_n}\tag 16.8
$$
where $k_n\geq 1$ in the last sum.
 Finally we note that we can
turn $V_{k+1}$, for $k\geq 1$ into $h_k$:
 $$
V_{k+1}=-\pa h_k+a^{k k_1...k_n}(\pa V_{k_1})\cdot\cdot\cdot(\pa
V_{k_{n-1}}) V_{k_n}\tag 16.9
 $$
 where $k_n\geq 1$ and $V_1=-\pa h$. This proves the general form
 (16.4) and it remains to prove the range of the indices in
 (16.5). To prove this we note that the terms in (16.8) are
 contractions over $n$ pairs of indices and this is still true
 for the terms we obtain by using (16.9) to replace the factors
 of $V$ by factors of $h$. This proves that
 $|\alpha_1|+...+|\alpha_n|=2n$. On the other hand when we replace
 factors of $V_{k+1}$ by $h_{k}$ the number of time derivatives
 go down by one for each factor, so we conclude that
 $l_1+...+l_n=l+2-2n$. This proves (16.5) apart from the last
 statement. That $|\alpha_i|\geq 1$ is clear and if $l_i=-1$ then $|\alpha_i|\geq
 2$, in view of (16.4), so in general $|\alpha_i|+l_i\geq 1$, and
 since $n\geq 2$, (16.5) follows.
\enddemo

We will now obtain a formal power series solution in the distance
to the boundary of the system for $h_l$ in Lemma {16}.1. In order
to do this, we will first choose simpler initial data for (16.2),
$$
f_0=y, \qquad v_0= \pa \phi,\qquad k_0=-\triangle \phi \tag 16.10
$$
where $\phi$ is to be determined. Let $h_{-1}=-\phi$, then (16.3)
hold also for $l=-1$ with $F_1=0$.

\proclaim{Lemma {16.}2} Suppose that $g_{ab}=\delta_{ab}$. Suppose
also that $h_{0,\,l}$ and $h_{1,\,l}$ are smooth for $l\geq -1$,
and let $F_l$ be as in Lemma {16.}1, and $F_{1}=0$. Then the
system
$$
\triangle h_{l}=h_{l+2}+F_{l+2}(h_{-1},...,h_{l-1}) ,\qquad l\geq
-1\tag 16.11
$$
with boundary conditions
$$
{h}_{l}\big|_{\pa\Omega}=h_{0,\,l}, \qquad \na_N h_{l}
\big|_{\pa\Omega}=h_{1,l} \qquad l\geq -1\tag 16.12
$$
has a formal power series solution in the distance to the
boundary:
$$
\overline{h}_{l}(r,\omega)\sim \sum h_{n,\,l}(\omega)
(1-r)^n/n!\tag 16.13
$$
Let $\chi$ be smooth such that $\chi(d)=1$, when $|d|\leq 1$,
$\chi(d)=0$, when $|d|\geq 2$ and $\chi\geq 0$. Then there are
$\varepsilon_{ln}>0$ such that
$$
\overline{h}_{l}(r,\omega)= \sum \chi((1-r)/\varepsilon_{ln})
h_{n,\,l}(\omega) (1-r)^n/n!\tag 16.14
$$
are smooth functions and such that (16.11) hold to infinite order
at the boundary.
\endproclaim
\demo{Proof} We have
$$
\triangle = \na_N^{\,\, 2} +\tr\theta \,\na_N +
\overline{\triangle}\tag 16.15
$$
where $\theta $ is the second fundamental form of the boundary and
$\overline{\triangle}$ is the tangential Laplacian on the
boundary. In the case $g_{ab}=\delta_{ab}$, $\na_N =\pa_r$ and
$\tr\theta=(n-1)/r$, where $r$ is the radial derivative.
Furthermore $\overline{\triangle}=r^{-2}\triangle_{\omega}$, where
$\triangle_\omega$ is the angular Laplacian on $S^{n-1}$. Hence we
have the system
$$
\pa_r^2 h_{l}\!=-\tfrac{1}{r^2}\triangle_\omega
h_{l}+h_{l+2}-\tfrac{n-1}{r}\pa_r h_{l}+
F_{l+2}(h_{-1},...,h_{l-1}),\qquad l\geq -1 \tag 16.16
$$
We want this system to be satisfied to all orders at the boundary,
so if
 $$
 h_{k,\,l}=\pa_r^k h_l\big|_{\pa\Omega},\qquad
F_{m,l} =\pa_r^m F_l(h_{-1},...,h_{l-3})\big|_{\pa\Omega}\tag
16.17
 $$
 we want
$$
h_{k,\,l}=\sum_{2i+2j+m\leq k,\, m\leq k-1} c_{ij\,km}
\triangle_\omega^{i} h_{m,\,l+2j}+\sum_{2i+2j+m\leq k}d_{ij\,k
m}\triangle_\omega^i F_{m,\,l+2j}\tag 16.18
$$
to hold for all $k\geq 2$ and $l\geq -1$, where $h_{0,\, l}$ and
$h_{1,\, l}$ are the given boundary conditions.

We now want to use induction. Note that the first term in the
right of (16.18) contains $h_{k^\prime,l^\prime}$ for
$k^\prime+l^\prime \leq k+l$, and $k^\prime<k$, and the second
term contains $h_{k^\prime,l^\prime}$ for $k^\prime+l^\prime<
k+l$, by the last inequality in (16.5). Assume that we found
$$
h_{k,\,l} ,\qquad \text{for}\qquad k+l\leq N,\qquad l\geq -1
,\quad k\geq 0\tag 16.19
$$
such that (16.18) hold for $k+l\leq N$ and $k\geq 2$.  Note that
if $N=0$ then $k\leq 1$ so there is nothing to prove. We know want
to find $h_{k,\,l}$ for $k+l=N+1$ such that (16.18) hold also for
$k+l\leq N+1$. This is again proven by induction. Assume that in
addition to (16.19) we found
 $$
 h_{k,\,l},\qquad\text{for}\quad k+l= N+1,\quad\text{and}\quad 0\leq k\leq M,\quad
 l\geq -1\tag 16.20
 $$
 such that (16.18) hold for $k\geq 2$. Note that for $M\leq 1$
 there is nothing to prove.
 Since $m$ in the first sum on the right of (16.18) is less than $k$ in
the left it follows that we can find $h_{k+1,\,l-1}$ such that
(16.18) hold.

In order to prove (16.14) we note that $h_{m,l}$ are smooth
functions on $S^{n-1}$. Hence we can use the usual trick of
choosing $\varepsilon_{ml}$ so small that $(\| h_{
m,l}\|_{m+1}+1)\,\varepsilon_{ml}\leq 1/2$, in which case the sum
converges in $H^{m}$ for any $m$ and $r$.
 \qed\enddemo

 Now, we want to find a formal power series solution in $t$ of the
 system (16.1) with initial data in (16.2) of the form
 $$
 f_0(y)=\tilde{f}_0(y) +y,\qquad
 v_0=\tilde{v}_{0}-\pa\, \overline{h}_{-1}, \qquad
 h_0=\tilde{h}_0+\overline{h}_0\qquad h_1= -\div \tilde{V}_0+\triangle\,
 \overline{h}_{-1} \tag 16.21
 $$
 where $\overline{h}_0$, $\overline{h}_{-1}$ are given by Lemma {16.}2
 and $\tilde{f}_0$, $\tilde{v}_0$ and $\tilde{h}_0$ vanish to
 infinite order at the boundary.
 Let $h_l$, for $l\geq 2$ be defined by (16.3).
Then it inductively follows that
 $$
 h_l=\tilde{h}_l +\overline{h}_l,\tag 16.22
 $$
 where $\overline{h}_l$ are as in Lemma {16.}2 and $\tilde{h}_l$ vanish
 to infinite order at the boundary. Therefore if we choose
 boundary data in (16.12) such that $h_{0,\,l}=0$ for $l\geq 0$,
 $h_{1,0}\leq c_0<0$. Then it follows that we can choose
 $\tilde{h}_0$ so that $h_0>0$ in $\Omega$ and $\na_N
 h_0\big|_{\pa\Omega}\leq -c_0<0$.
 Moreover it follows that
 $$
 {h}_l\big|_{\pa\Omega}=0,\qquad l\geq 0\tag 16.23
 $$
and hence the compatibility conditions are satisfied to all
orders.

We can now construct smooth functions in
$[0,T]\times\overline{\Omega}$, satisfying the initial conditions
(16.2) and the equations (16.2) to infinite order as $t\to 0$:
 \proclaim{Lemma {16.}3} Suppose that initial data $f_0$, $v_0$ and
 $h_0$ and $h_1=-\div V_0$ satisfy the compatibility conditions
 for all orders, i.e. if $h_l$ for $l\geq 2$ are defined by (16.3)
 then
 $$
 h_l\big|_{\pa\Omega}=0,\qquad \text{for}\qquad l\geq 0.\tag 16.24
$$
 Then there are smooth functions $(x,h)$ in $[0,T]\times\overline{\Omega}$,
 such that (16.2) hold, (16.1)
 is satisfied to infinite order as $t\to 0$ and
 $h\big|_{\pa\Omega}=0$.
 \endproclaim
 \demo{Proof} Let $\chi$ be as in Lemma {16.}2 and set
 $$
  h(t,y)=\sum_{l=0}^\infty \chi(t/\varepsilon_l) h_l(y) t^l/l!\tag
  16.25
  $$
  where $\varepsilon_l>0$ are chosen so that
  $(\|h_l\|_{m}+1)\varepsilon_m\leq 1/2$. Then it follows that the
  sum converges in $H^m$ for any $m$ so $h$ is smooth and
  satisfies $h\big|_{\pa\Omega}=0$ and $D_t^l h\big|_{t=0}=h_l$.
  Furthermore let $x(t,y)$ be defined by
  $$
  D_t^2 x_i=-\pa_i h,\qquad x_i\big|_{t=0}=f_0,\qquad D_t
  x\big|_{t=0}=v_0.\qed\tag 16.26
  $$
  \enddemo

\head 17. The general case, when the enthalpy is a strictly
increasing function of the density. \endhead
 We will now outline how to generalize the existence result
 obtained for $e(h)=ch$ to the case when $e(h)$ is a smooth
 strictly increasing function satisfying $1/c_1^\prime \leq
 e^\prime(h)\leq c_1^\prime$. First we will show that the functional
 $h=\Psi(x)$, i.e. the solution of (2.21), exist for $x$ in a bounded set,
 $\||u\||_{r_0+r_1,\infty}\leq 1$, and for $T$ sufficiently
 small. Since in the Nash-Moser iteration we need a bound
 for $\||h\||_{3,\infty}$, we will first show that we can obtain
 such a bound as well as a bound for $\|h\|_{r_0+3}$ independent
 of $T$. Local existence for the nonlinear wave equation follows
 from a standard argument using essentially the same estimates,
 so it is just a question of showing that
 we have {\it a priori} bounds up to some time $T>0$ that only
 depends on the approximate solution $(x_0,h_0)$ and is
 independent of $x=u+x_0$ as long as $\||u\||_{r_0+r_1,\infty}\leq 1$.
Once we have a bound for $\||h\||_{4,\infty}$, the bound for
higher derivatives follows from this.

The equation we study in this section is
$$
D_t\big( e^\prime (h) D_t h\big) -\triangle h ={f},\qquad h\big|_{t=0}=0,\qquad
f=(\pa_i V^j)(\pa_j V^i) \tag 17.1
$$
where
$$
e^\prime+1/e^\prime\leq c_1^\prime,\qquad
\sum_{a,b} |g^{ab}|+|g_{ab}|\leq c_1^\prime ,\qquad
|\pa x/\pa y|^2+|\pa y /\pa x|^2\leq c_1^\prime\tag 17.2
$$
for some constant $0<c_1^\prime<\infty$. Let $K_{10}^\prime$ denote a continuous
function of $c_1^\prime$ that also depends on the order of differentiation $r$
but is independent of a lower bound for $T$.

We prove the following Theorem:
\proclaim{Theorem {17.}1} Let $r_0=[n/2]+1$ be the Sobolev exponent and let
$k\geq 1$. There are continuous function $C_k$ and $D_k$ such that if
$T>0$ is so small that
$$
T C_1(\|| x_0\||_{r_0+2+1,\infty},\||h_0\||_{r_0+2+1,\infty})\leq 1,\tag 17.3
$$
and $u=x-x_0$ is small that
$$
\||u\||_{r_0+2+k,\infty}\leq 1\tag 17.4
$$
then (17.1) has a smooth solution for $0\leq t\leq T$ satisfying
$$
\||h\||_{k,\infty}\leq
D_k(\|| x_0\||_{r_0+2+k,\infty},\||h_0\||_{r_0+2+k,\infty}).\tag 17.5
$$
Furthermore, there is a continuous function $K_2^\prime$ of
$$
\||x\||_{2,\infty}+\||h\||_{2,\infty}+\||h_0\||_{2,\infty}+1/T+c_1^\prime\tag 17.6
$$
depending also on $r$ such that
$$
\||h\||_{r,\infty}\leq K_2^\prime\big( \|| x\||_{r+r_0+2,\infty}
+\||h_0\||_{r+r_0+2,\infty}\big) \tag 17.7
$$
\endproclaim
\demo{Proof} Since by Sobolev's Lemma
$$
\||h\||_{1,\infty}\leq C(\|h\|_{1,r_0}+\|h\|_{0,r_0+1})\tag 17.8
$$
where $C$ is independent of $T$. It follows from Lemma {17.}2 that
$$
d \tilde{E}_r/dt\leq M_r(\tilde{E}_{r+1}+1)^{r+2},\qquad
r\geq r_0\tag 17.9
$$
where $M_r= K_{10}^\prime\||x\||_{r+3,\infty}^{(r+2)^2}$.
Integrating this inequality gives
$- d(\tilde{E}_{r+1}+1)^{-(r+1)}/dt\leq M_r(r+1) $
and hence $(\tilde{E}_{r+1}(t)+1)^{-(r+1)}\geq (\tilde{E}_{r+1}(0)+1)^{-(r+1)}
-M_r(r+1) t$. If $M_r(r+1) t \leq (\tilde{E}_{r+1}(0)+1)^{-(r+1)}/2$ it follows
that $\tilde{E}_{r+1}(t)+1\leq 2(\tilde{E}_{r+1}(0)+1)$.
Since $\tilde{E}_{r+1}(0)\leq K_{20}^\prime\sum_{s=0}^r
\big[\big[\big[\pa x\big]\big]\big]_{r-s,\infty} (\|h_0\|_{1,s}+\|h_0\|_{0,s+1})$
this proves the first part of the theorem for $k=1$ and for $k\geq 2$ it follows from
also using Theorem {17.}3.

It follows from  Lemma {17}.2 and interpolation in space time that
$$
\frac{d \hat{E}_{r+1}}{dt}\leq K_2^\prime \big( \hat{E}_{r+1}+\||x\||_{r+3,\infty}
\big),\qquad\quad\text{if}\qquad
 \hat{E}_{r+1}=\sum_{s=0}^{r} \|| x\||_{r+2-s} \tilde{E}_{s+1}\tag 17.10
$$
Multiplying by the integrating factor $e^{K_2^\prime t}$ and integrating, we get
$$
\hat{E}_{r+1}(t)\leq K_2^\prime \big(  \hat{E}_{r+1}(0)+\||x\||_{r+3,\infty}\big)
\tag 17.11
$$
By Theorem {7.}3
$
\|h\|_{r+1}\leq K_2^\prime \hat{E}_{r+1}
$
and hence
$$
\|h(t,\cdot)\|_{r+1}\leq K_2^\prime
\big(\sum_{s=0}^{r} \|| x\||_{r+2-s}  \|h_0(0,\cdot)\|_{s+1}+ \|| x\||_{r+3} \big)
\tag 17.12
$$
Using Sobolev's lemma and interpolation the estimate for $\||h\||_{r,\infty}$
follows.
\qed\enddemo

\proclaim{Lemma {17.}2}
 Suppose that $g^{ab}$ and $e^\prime =e^\prime (h)$ are smooth and satisfy
(7.4). For $s\geq 0$ let
$$
E_{s+1}(t)=\Big(\frac{1}{2}\!\int_\Omega{ \!\! e^\prime
(\!{D}_t^{s+1}h)^2\! + g_{ab} (\hat{D}_t^s H^a ) (\hat{D}_t^s H^b
)\kappa d y}\Big)^{1/2} , \qquad H^a= g^{ab}\pa_b h, \tag 17.13
$$
$E_{0}(t)=\int_{\Omega} h^2\, \kappa dy$ and
$\tilde{E}_{r+1}=\sum_{s=0}^{r+1} E_s$.
Then for $r\geq 0$
$$\align
\|h\|_{1,r}+\|h\|_{0,r+1}&\leq C\sum_{s=0}^{r}
\big[\big[\big[\pa x\big]\big]\big]_{r-s,\infty} \tilde{E}_{1+s} \tag 17.14\\
 \tilde{E}_{1+r}&\leq C\sum_{s=0}^{r}
\big[\big[\big[\pa x\big]\big]\big]_{r-s,\infty}
\big(\|h\|_{1,s}+\|h\|_{0,s+1}\big)\tag 17.15
\endalign
$$
and
$$
\frac{d \tilde{E}_{r+1}}{dt}\leq K_{10}^\prime\sum_{k=0}^{r+1}
\|h\|_{1,\infty}^k\Big( \sum_{s=0}^{\min(r+1-k,r)}
\big[\big[\big[ \pa x\big]\big]\big]_{r+1-k-s,\infty} \tilde{E}_{s+1} +
 \big[\big[\big[ \pa x\big]\big]\big]_{r+2-k,\infty} \Big)\tag 17.16
$$
\endproclaim
\demo{Proof} The proof is similar to that of Lemma {7.}2. We have
$$
\multline \frac{d E^2_{s+1}}{dt}= \int_\Omega{ \Big( e^\prime
({D}_t^{s+1}h)({D}_t^{s+2}h)+ g_{ab}(\hat{D}_t^s H^a )
(\hat{D}_t^{s+1} H^b )
\Big)\,  \kappa \,d y}\\
+\frac{1}{2}\int_\Omega{ \Big((\hat{D}_t e^\prime ) ({D}_t^{s+1}h)^2+
(\check{D}_t g_{ab})(\hat{D}_t^s H^a ) (\hat{D}_t^s H^b ) \Big)\,
\kappa \,d y}.
\endmultline\tag 17.17
$$
Here the terms on the second row are bounded by
$K_{10}^\prime (\|h\|_{1,\infty}+\|g\|_{1,\infty})E_{s+1}^2$.
The first row of the right hand side is up to lower order
$$
\multline  \int_\Omega{ \! \left( (
{D}_t^{s+1}h)\big(\hat{D}_t^{s}D_t(e^\prime  D_t h)\big)+ (\hat{D}_t^s
H^a) (\pa_a {D}_t^{s+1}h )\right)\,  \kappa \,d y}\\
=\int_\Omega{  ({D}_t^{s+1}h ) \left(\hat{D}_t^{s}D_t( e^\prime  D_t
h)-\kappa^{-1} \pa_a\big(\kappa \hat{D}_t^s H^a \big)\right)\,
\kappa dy}=\int_\Omega{ ({D}_t^{s+1}h ) (\hat{D}_t^{s} {f}) \kappa
dy},
\endmultline\tag 17.18
$$
where we have integrated by parts using that
$D_t^{s+1}h\big|_{\pa\Omega}=0$, that $D_t\big( e^\prime  D_t h\big)
-\kappa^{-1}\pa_a\big(\kappa H^a\big)={f}$ and that $\hat{D}_t^s
\big(\kappa^{-1}\pa_a\big(\kappa
H^a\big)\big)=\kappa^{-1}\pa_a\big(\kappa \hat{D}_t^sH^a\big)$.
In fact, using Lemma {6.}3 we get, since $\underline{H}_a=\pa_a h$,
 $$
 \hat{D}_t^{s+1}H^a -g^{ab}\pa_a
 D_t^{s+1}h=
 -\sum_{i=0}^{s} \binom{s+1}{i}g^{ab}(\check{D}_t^{s+1-i} g_{bc}) \hat{D}_t^i
 H^c\tag 17.19
 $$
 and, since $\hat{D}_t D_t^2 e(h)=\kappa^{-1}D_t^s\big(\kappa D_t^2 e(h))\big)$,
 $$
 \hat{D}_t^s D_t^2 e(h) -e^\prime(h) D_t^{s+2}h= \sum_{i=0}^{s-1} \binom{s}{i}
 \kappa^{-1}( D_t^{s-i}\kappa)(D_t^{2+i} e(h) )
+D_t^{2+s} e(h)-e^\prime(h)D_t^{s+2} h.\tag 17.20
$$
Here the right hand side is bounded by a constant times
$$
\sum_{i=0}^{s-1}\kappa^{-1}|D_t^{s-i}\kappa|\!\!\!\!\!\!\!
\sum_{\Sb s_1+...+s_k=i+2,\\ \,s_i\geq 1,\, k\geq 1\endSb}\!\!\!\!\!
|e^{(k)}(h)|\, |D_t^{s_1}h|\!\cdot\cdot\cdot|D_t^{s_k} h|
+\sum_{\Sb s_1+...+s_k=s+2,\\ \,s_i\geq 1,\, k\geq 2\endSb}\!\!\!\!\!
|e^{(k)}(h)|\, |D_t^{s_1}h|\!\cdot\cdot\cdot|D_t^{s_k} h|\tag 17.21
$$
The $L^2$ of this can be estimated by Theorem {7.}3:
$$
 K_{1,0}^\prime\sum_{k=0}^{s+1}
\|h\|_{1,\infty}^k\Big(\!\!\!\!\! \sum_{i=1}^{\min(s+2-k,s+1)}\!\!\!\!\!
\big[\big[\big[\,\pa x\big]\big]\big]_{s+2-k-i,\infty}
\big( \|h\|_{0,i}+\|h\|_{1,i-1}\big)+
\big[\big[\big[\pa x\big]\big]\big]_{s+2-k,\infty}\Big)\tag 17.22
$$
By Lemma {6.}3
 $\hat{D}_t^{s} (g^{ab}\pa_b h)=
 g^{ab}\pa_b D_t^{s} h+ \sum_{i=0}^{s-1}
 \binom{s}{i}
g^{ab}(\check{D}_t^{s-i} g_{bc})
 \hat{D}_t^{i} (g^{cb}\pa_b h)$ so it follows that
 $\|\pa D_t^{s} h\|\leq \sum_{i=0}^s
 \big[\big[\big[\pa x\big]\big]\big]_{r-i}\tilde{E}_{1+i}$
Since also $\|\hat{D}_t^s f\|\leq K_{10}\big[\big[\big[\pa x\big]\big]\big]_{s+2}$
the lemma follows.
\qed\enddemo

\proclaim{Theorem {17.}3}  Suppose that
$D_t^2 e(h)-\triangle h={f}$, where $f=(\pa_i V^j)(\pa_j V^i)$.
 Suppose also that $|e^{(k)}(h)|\leq C_k$. Then we have
$$
\|\, \big[ h\big]_r^j \,\|\leq K_{1,0}^\prime \sum_{k=j-1}^r
\|h\|_{1,\infty}^k\Big(\sum_{i=1}^{r-k}
\big[\big[\big[\,\pa x\big]\big]\big]_{r-k-i,\infty}
\big( \|h\|_{i-1,1}+\|h\|_{i,0}\big)+
\big[\big[\big[\,\pa x\big]\big]\big]_{r-k,\infty} \Big)\tag 17.23
$$
and
$$
\|\, \big[ h\big]_r^j \,\|\leq K_{1,0}^\prime\sum_{k=j-1}^r
\|h\|_{1,\infty}^k\Big( \sum_{i=1}^{r-k}
\big[\big[\big[\,\pa x\big]\big]\big]_{r-k-i,\infty}
\big( \|h\|_{0,i}+\|h\|_{1,i-1}\big)+
\big[\big[\big[\,\pa x\big]\big]\big]_{r-k,\infty}\Big)\tag 17.24
$$
where
$$
[h]_r^j=\!\!\!\sum_{\Sb r_1+...+r_k=r,\\ r_i\geq 1,\,
k\geq j\endSb}\!\!\!\!\!\!\!\! |h|_{r_1}\!\cdot\cdot\cdot |h|_{r_k}
,\qquad\quad
\big[\big[\big[\,\pa x\big]\big]\big]_{r,\infty}=\!\!\!
\sum_{\Sb r_1+...+r_k=r,\\ r_i\geq 1\endSb}
\!\!\!\!\!\!\!\! \||\pa x\||_{r_1,\infty}
\!\cdot\cdot\cdot \||\pa x\||_{r_k,\infty}\tag 17.25
$$
\endproclaim
\demo{Proof} The first inequality follows from Lemma {17.}4 below and interpolation
and the second follows from the first and Lemma {17.}6 below.
\qed\enddemo

\proclaim{Lemma {17.}4} Suppose that
$D_t^2 e(h)-\triangle h={f}$, where $f=(\pa_i V^j)(\pa_j V^i)$.
 Suppose also that $|e^{(k)}(h)|\leq C_k$. With notation as in
Definition {5.}1 we have
$$
[h]_{r,s}\leq K_{10}^\prime \sum_{i+j\leq r+s}
\big[\pa x\big]_{r+s-i-j} \, \big[h\big]_{i,(j,1)}\,\tag 17.26
$$
where $\big[\pa x\big]_l$ is as in Definition {5.}1,
$$
\big[h\big]_{r,s}=\!\!\!\!
\sum_{\Sb r_1+...+r_k=r,\\ s_i+...+s_k=s\\ s_i+ r_i\geq 1\endSb}
\!\!\!\! |h|_{r_1,s_1}\cdot\cdot\cdot |h|_{r_k,s_k},
\qquad \big[h\big]_{r,(s,m)}=
\!\!\!\!\sum_{\Sb r_1+...+r_k=r,\\ s_i+...+s_k=s\\ s_i+ r_i\geq 1,\, s_i\leq m,
\endSb}
\!\!\!\! |h|_{r_1,s_1}\cdot\cdot\cdot |h|_{r_k,s_k},\tag 17.27
$$
and $\big[h\big]_{0,0}=1$, $\big[h\big]_{0,(0,1)}=1$.
\endproclaim
\demo{Proof} If $h=\eta(e)$ is the inverse of $e(h)$, then
$$
|h|_{r,s}\leq C\sum_{r_1+...+r_k=r,\,\, s_1+...+s_k=s,\, r_i+s_i\geq 1}
|\eta^{(k)}(e)|\, |e|_{r_1,s_1}\cdot\cdot\cdot |e|_{r_k,s_k}. \tag 17.28
$$
Since,
$$
|e|_{r,s}\leq |\triangle h|_{r,s-2}+|f|_{r,s-2} ,
\qquad\qquad\qquad\qquad\qquad\quad s\geq 2\tag 17.29
$$
where
$$
|\triangle h|_{r,s-2}\leq K_{10} \sum_{1\leq i\leq r+2,\,\,j\leq s-2}
\big[ \pa x\big]_{r+s-i-j}\, |h|_{i,j} ,\qquad s\geq 2\tag 17.30
$$
and
$$
|f|_{r,s-2}\leq K_{10} \big[ \pa x\big]_{r+s} ,
\qquad\qquad\qquad\qquad\qquad\qquad\qquad s\geq 2\tag 17.31
$$
we obtain
$$
|e|_{r,s} \leq K_{10} \sum_{ i\leq r+2,\,\,j\leq s-2}
\big[ \pa x\big]_{r+s-i-j}\, [h]_{i,j} ,\qquad\qquad\qquad s\geq 2\tag 17.32
$$
It follows that
$$
\big[ h\big]_{r,(s,m)}\leq
K_{10}^\prime\sum_{ i+j\leq r+s}
\big[ \pa x\big]_{r+s-i-j}\, [h]_{i,(j,m-1)} ,\qquad\qquad\qquad m\geq 2\tag 17.33
$$
and the lemma follows by induction.
\qed\enddemo

\proclaim{Lemma {17.}5} For $r\geq 1$ we have with notation as in Definition {5.}3,
$$
\|h\|_r\leq K_{10}^\prime\sum_{j=0}^r \|h\|_{1,\infty}^j\Big(
\sum_{i=1}^{r-j}
\big[\big[\big[\pa x\big]\big]\big]_{r-i-j,\infty} \,
\big( \|h\|_{0,i}+\|h\|_{1,i-1} \big)+
\big[\big[\big[\pa x\big]\big]\big]_{r-j,\infty}\Big)\tag 17.33
$$
\endproclaim
\demo{Proof} We have
$$
|\triangle h|_{r,s}\leq |e(h)|_{r,s+2}+|{f}|_{r,s}\tag 17.34
$$
where
$$
|e(h)|_{r,s}\leq e^\prime(h)|h|_{r,s}+ C\sum_{k=2}^{r+s}|e^{(k)}(h)|\,
[h]_{r,s}^{\, k}
,\qquad \text{where} \quad
 [h]_{r,s}^{\, k}=\!\!\!\! \sum_{\Sb r_1+...+r_k\,=r,\\ s_1+...+s_k\,=s\\
 r_i+s_i\geq 1\endSb }\!\!\!\!
|h|_{r_1,s_1}\!\cdot\cdot\cdot |h|_{r_k,s_k}\tag 17.35
$$
and by the previous lemma and interpolation in space only
$$
\| [h]_{r,s}^{\, k} \|\leq K_{10}^\prime \sum_{i+j\leq r+s,\, j\geq k-1,i\geq 1}
\big[\big[\big[\pa x\big]\big]\big]_{r+s-i-j,\infty}  \| h\|_{1,\infty}^{j}
\big(\|h\|_{i,0}+\|h\|_{i-1,1}\big)\tag 17.36
$$
Using Theorem {6}.1 we get for $r\geq 2$
$$\multline
\|h\|_{r+2,s}\leq K_{10}\!\!\!\!\sum_{i\leq r,\, j\leq s}\!\!\!\! \big[\big[\big[
\pa x\big]\big]\big]_{r+s-i-j,\infty}
\|h\|_{i,j+2}
+K_{10}\sum_{j=0}^s  \big[\big[\big[
\pa x\big]\big]\big]_{r+s+1-j,\infty}\|h\|_{1,j,\infty}\\
+K_{10}^\prime \big[\big[\big[\pa x\big]\big]\big]_{r+s+2,\infty}
+K_{10}^\prime \sum_{i+j\leq r+s+2,\, i,j\geq 1}
\big[\big[\big[\pa x\big]\big]\big]_{r+s+2-i-j,\infty}  \| h\|_{1,\infty}^{j}
\big(\|h\|_{i,0}+\|h\|_{i-1,1}\big)\endmultline \tag 17.37
$$
Hence using induction, we get
$$\multline
\|h\|_{r}\leq K_{10}^\prime\sum_{1\leq i\leq r} \big[\big[\big[
\pa x\big]\big]\big]_{r-i,\infty}
\big( \|h\|_{0,i}+\|h\|_{1,i-1}\big)\\
+K_{10}^\prime \big[\big[\big[\pa x\big]\big]\big]_{r,\infty}
+K_{10}^\prime \sum_{i+j\leq r,\, i,j\geq 1}
\big[\big[\big[\pa x\big]\big]\big]_{r-i-j,\infty}  \| h\|_{1,\infty}^{j}\|h\|_{i}
\endmultline \tag 17.38
$$
Since $\|h\|_1\leq \|h\|_{0,1}+\|h\|_{1,0}$ the lemma now follows from induction.
\qed\enddemo

Alternatively, one can use interpolation interpolation in space-time.
\proclaim{Lemma {17.}6} We have
$$\align
\|| D_t^l h\||_{L^{2/(m-l)}}&\leq C_T\||h\||_{L^\infty}^{l/m} \||
D_t^m h\||_{L^2}^{1-l/m}\tag 17.39\\
\|| \, |D_t^{l_1} h|\cdot\cdot\cdot |D_t^{l_k} h|\, \||_{L^2}
&\leq C_T\||h\||_{L^\infty}^{k-1} \|| D_t^m h\||_{L^2}\tag 17.40
\endalign
$$
where $m=l_1+...+l_k$. Furthermore, if $h=\tilde{h}+h_0$ where
$\tilde{h}$ vanish to infinity order as $t\to 0$ then, with a
constant independent of $T$ but depending on $h_0$ we have
$$
\|| \, |D_t^{l_1} h|\cdot\cdot\cdot |D_t^{l_k} h|\, \||_{L^2} \leq
C_m(\||h\||_{\infty}, \||h_0\||_{\infty, m})\big(\sum_{j=0}^m \||
D_t^j h\||_{L^2}+1)\tag 17.41
$$
\endproclaim
\demo{Proof} The interpolation inequalities are standard and it is
also standard that the constant is independent of $T$ if $h$
vanish to infinite order as $t\to 0$. Hence in the product we can
estimate $\||D_t^{l_i} h\||_{L^{2/(m-l_i)}} \leq
\||D_t^{l_i}\tilde{h}\||_{L^{2/(m-l_i)}} +\||D_t^{l_i}
h_0\||_{L^\infty} $. \qed\enddemo

 \subheading{Acknowledgments} I would
like to thank Demetrios Christodoulou, Richard Hamilton and Kate
Okikiolu for many long and helpful discussions.

\end

 \proclaim{Lemma {16.}3 } The system
$$
h^{(N+2)}=f,\quad h^{(N+1)}=g,\qquad \triangle h^{(l)}=
h^{(l+2)}+F_{l+2},\quad h^{(l)}\big|_{\pa\Omega}=0,\quad -1\leq
l\leq N\tag 16.11
$$
where $F_l$ is as in Lemma {16.}1 can be solved if $f$ and $g$ are
sufficiently small. and $N$ sufficiently large.
\endproclaim
\demo{Proof} Inverting the Laplacian we hence get an estimate
$$
\| h^{k}\|_{N+2-k}\leq \| h^{k+2}\|_{N-k} +\|F_{k+2}\|_{N-k}\tag
16.12
$$
and hence
$$
\sum_{k=0}^N\| h^{k}\|_{N+2-k}\leq
C\sum_{k=0}^N\|F_{k+2}\|_{N-k}+\|f\|_0+\|g\|_1\tag 16.13
$$
Furthermore, using the expression for $F_l$ in Lemma {20}.1 and
Sobolev's lemma to estimate $L^\infty$ norms by $L^2$ norms we get
an estimate
$$
C\sum_{k=0}^N\|F_{k+2}(h)\|_{N-k}\leq  G\Big(\sum_{k=0}^N\|
h^{k}\|_{N+2-k}\Big),\qquad G(0)=G^\prime(0)=0\tag 16.14
$$
Hence if $\|f\|_0$ and $\|g\|_1$ are sufficiently small we get
existence for the system (16.11) by iteration.
\enddemo

If we now let
$$
h_0(t,y)=h^{(0)}+...+h^{(N+2)}(y) t^{N+2}/(N+2)!, \quad D_t
x_0\big|_{t=0}=-\pa h^{(-1)}(y),\qquad D_t^2 x_0=-\pa_i h_0\tag
16.15
$$

\comment
 We also remark that since $V$ is irrotational at time $0$
it follows that $\triangle h^{(0)}=-|\pa V|^2+\varepsilon^2
h^{(2)}$ is negative and so if ... Now, let $N+2=4k$ and let
$f\geq 0$. Then $h^{(l)}$ is positive if $l$ is a multiple of $4$
plus something of order $O(f^2)$. In particular $h^{(2)}$ is
negative plus something of order $O(f^2)$.
\endcomment

Furthermore, we can use this to construct a power series around
the boundary at time $0$. In fact,
$$
\triangle h^{(l)}=\varepsilon^2 h^{(l+2)}+F_l\tag 16.12
$$
so
$$
\triangle^{k} h^{(l)}=\varepsilon^{2k} h^{(l+2k)} +
\sum_{j=0}^{k-1} \varepsilon^{2j} \triangle^j F_{l+2j}\tag 16.13
$$
At the boundary, we require that
$$
h^{(l)}=0, \qquad \na_N h^{(l)}=h_l,\tag 16.14
 $$
Hence at the boundary we obtain the equations
$$
\triangle^{k} h^{(l)}= \sum_{j=0}^{k-1} \varepsilon^{2j}
\triangle^j F_{l+2j}\tag 16.15
$$
and
 $$
\na_N \triangle^{k} h^{(l)}= \varepsilon^{2k}\na_N  h_{l+2k}+
\sum_{j=0}^{k-1} \varepsilon^{2j} \na_N \triangle^j F_{l+2j}\tag
16.16
 $$
 Here $F_k$ is a multi-linear form in
 $$
 \pa^\alpha h^{(j)},\qquad \text{for}\qquad
 |\alpha|+j\leq k+1,\qquad
 j\leq k-1\tag 16.17
 $$
  with coefficients depending
 space derivatives of $x$ and $v$.
  We now want to make a power series expansion in $(t,d)$ where
$d$ is the distance to the boundary. Assume therefore that we know
(2.20) for $k\leq r$ and we want to prove (2.20) for $k=r+1$.
Since the tangential derivatives are determined are in the
coefficient we only need to determine normal derivatives.

We give initial conditions for $h^{(l)}\big|_{\pa\Omega}=0$ and
$\na_N h^{(l)}\big|_{\pa\Omega}=h_l$ are given. Then we get
$\triangle h^{(0)}=F_0$ where $F_0$ can be expressed in terms of
$V$ and $x$ at time $0$. This gives us the second normal
derivative of $h^{(0)}$, $\na_N^2 h^{(0)}$. After that we get
$\triangle h^{(1)}=F_1$, where $F_1$ depends on $\pa^2 h^{(0)}$,
which gives us the second normal derivative of $h^{(1)}$, $\na_N^2
h^{(1)}$.
 Then we use (16.16) to get $\na_N\triangle  h^{(l)}=\na_N h^{(l+2)} +
 \na_N F_l$, for $l=0,1$.

Then we get $\triangle^2 h^{(0)}=\triangle F_0+F_2$ and $\triangle
h^{(2)}=F_2$ where $F_2$ depends on $\pa^2 h^{(1)}$ and on $\pa^3
h^{(0)}$.

Well, there seems to be a problem here since we also need that
$$
D_t e(h)=-\div V\tag 16.18
$$
in order that the initial conditions for $h$ be such that it
corresponds to a  solution of the original equation. We must first
assume that $\div V_0\big|_{\pa\Omega}=0$.
 Let us therefore see what happens if instead we try to find a
 power series for (16.18) and
 $$
 D_t v_i=-\pa_i h\tag 16.19
 $$
 Then if $e(h)=c^2 h$
 $$
 D_t^2 v_i=\pa_i\div V +(\pa_i V^k)\pa_k h\tag  16.20
 $$
 and since $D_t\div V=\div D_t V-(\pa_k V^j)\pa_j V^k$
 $$
 D_t^3 v_i=\pa_i \big(-\triangle h -(\pa_j V^k)(\pa_k V^j)\big)
 -2(\pa_i V^k)\pa_k \div V -2(\pa_i V^j)(\pa_j V^k)\pa_k h
 -(\pa_i \pa_k h)\pa_k h \tag 16.21
 $$
 In any case, it is clear that the compatibility conditions are
 exactly that $D_t^k \div V\big|_{\pa\Omega}=0$.

Let us now also make the assumption that $x(0,y)=y$ so the
equations become simpler. We could further assume that $V$ is
irrotational and perhaps even that $\div V_0$ vanishes to infinite
order at the boundary at time $0$.  Let us therefore assume that
the vector field $V$ is irrotational at the boundary as well as
divergence free. Perhaps we should also assume that it is
irrotational to infinite order at the boundary at time $0$, i.e.
that $V_0=\pa \phi$ to infinite order at the boundary. Assume then
that we are given some boundary data for
 $\phi$ and for $\na_N \phi$ as well as for $h^{(l)}$ and $\na_N
 h^{(l)}$. Since $D_t h=-\triangle \phi$, this then gives an
 equation $\triangle\phi=\na_N^2 \phi+\tr\theta \na_N\phi
 +\overline{\triangle}\phi$ so this first give us
 $\na_N^2\phi$ then we get $\triangle h=(\pa_i V^j)(\pa_j V^i)$
 so we get $\na_N^2 h$. Similarly $\triangle D_t h=
 (\pa_k h)\triangle V^k+(\pa V)\pa^2 h$.
 At this point we don't quite know $(\pa_k h)\triangle V^k$
 more than that it is something like $\na_N\div V$ if indeed the
 vector field is irrotational. Therefore, it seems like we first
 have to get some more information about $V$ first.
 Since we know $\na_N D_t h$ at the boundary we get
 $\na_N D_t h=\na_N^3 \phi+\na_N^2 \phi+\triangle \na_N\phi
 +\text{L.O.}$. Therefore, this gives us $\na_N^3 \phi$ and hence
 we get $\na_N\div  V$, actually, that one we also get directly.
 Okey so now we have $\na_N^2 D_t h$, $\na_N^3 \phi$ and $\na_N^2
 h$. But at the same time we also get $\na_N^3 h
 =\na_N^3 \phi$ etc. Hence it all works out fine.

 Once we have found the power series expansion for $h$ and $V$
 we truncate the series and we can hopefully modify it so
 $D_t h_0=-\div V_0$. We would like to have this identity exactly
 in order that the initial conditions are of the form we want in
 (2.11). However,

 Now, let us try to find a formal solution at the boundary of the
 form $V=\pa \phi$, i.e. we will assume that $V$ is this plus something
 that vanishes to infinite order at the boundary.
 Then Euler's equations become
 $$
 \pa_i\big( D_t\phi-|\pa \phi|^2/2+h\big)=0
 $$
 Since this is supposed to hold at the boundary in particular for
 tangential derivatives it follows that we must have
 $$
D_t\phi-|\pa \phi|^2/2+h=c(t)
$$
Since also $D_t h=-\div V=-\triangle \phi$ we get if we also
assume that $c(t)=0$.
$$
D_t^2 \phi-\triangle \phi =-g^{ab}(\pa_a \phi)\pa_b h
$$

 \head
3. The physical conditions and initial conditions.\endhead. (The
difficulty comes when we want to do the same estimate for the wave
equation using the Nash Moser stuff for the compressible case.
Well, hopefully on can get away with just using tangential
derivative sin that case as well. Although I don't see how yet. In
that case we need to estimate this for the variation as well as
for the pressure (or enthalpy ) itself. Well, I don't know how to
do this jet. How about the space-time divergence? Well I don't
know that it works. The space-time divergence is a mess but
perhaps its okey. It should work fine to do this actually if we
have estimates fro both the time derivative, the tangential
derivatives and the curl and the divergence. Okey, so it all looks
fine then. This really looks great. Actually,its probably
sufficient with two time derivatives and the rest space.

Then one can construct an approximate solution ${x}_0$ for
$-T<t<T$ for some $T>0$, satisfying the initial conditions and
such that
$$
D_t^k \Phi({x}_0)\big|_{t=0}=0,\qquad\text{for all}\quad k\geq
0,\qquad {x}_0\big|_{t=0}=f_0,\quad D_t{x}_0\big|_{t=0}=v_0
$$
Let $x_\delta(t,y)=x(t-\delta,y)$ and let
${F}_\delta=\Phi({x}_\delta)$ and let Let
$\overline{F}_\delta=F_\delta$ when $t>\delta $ and $0$ when
$t<\delta$. Now, let
$$
\Psi(x)=\Phi(x+x_0)-\Phi(x_0), \qquad x\big|_{t=-T}=0,\quad D_t
x\big|_{t=-T}=0
$$
I.e. we think of $\Psi$ as the mapping such that we have these
initial conditions for a given right hand side that we want to
solve $\Psi(x)=F$. We note that $\Psi(0)=0$ and we want to use
that one can find a solution in a neighborhood of $0$ to find a
solution for the initial value problem. Now if we can solve
$$
\Psi(x)=\overline{F}_\delta -\overline{F}_0
$$
for any $\delta>0$ we get a solution to the problem since for
$0\leq t<\leq \delta$ this is equivalent to
$$
\Phi(x+x_0)=0, \qquad 0\leq t\leq \delta, \qquad x(t,y)=0,\qquad
\text{for}\quad t<0
$$
Hence $x+x_0$ defined a solution of our problem for a short time.
But
$$
\overline{F}_\delta-\overline{F}_0\to 0, \qquad
\overline{F}_\delta-\overline{F}_0=0,\qquad t<0
$$
so we only have to show that the equation
$$
\Psi(x)=F
$$
has a solution for $F$ sufficiently small and supported in $t\geq
0$. This last condition is also important because we are going to
use smoothing in time in the interval $-t_0<t<t_0$ and we don't
want it to mess up with the initial conditions.

 If $G$ denote the right hand side then
$$
\|\dot{G}\|_{s}
$$

  Applying Theorem {7.}1 therefore gives
$$
\|h\|_r\leq K
$$

 $$
D_t^2 (c \dot{h})-\triangle\dot{h} =-\pa_i \big( ( \pa^i V^k)\pa_k
h\big)-(\pa^i V^k)\pa_i \pa_k h +2(\pa_i V^k)\pa_k D_t V^i ,\qquad
\dot{h}\Big|_{\pa \Omega}\!\!\!=0\tag 10.1
$$
If $G$ denotes the right hand side above then
 $$
 \|\dot{G}\|_{s-2}\leq K \sum_{k=0}^s \||x\||_{s+2-k,\infty}
 \big(\|\dot{h}\|_k+\|h\|_k +\|V\|_{k}\big)    \tag 10.2
 $$
Applying the estimate in Theorem {7.}1 therefore gives
$$
\|\dot{h}\|_r\leq K\sum_{s=0}^r \||x\||_{r+s+3-k}  \int_0^t
\big(\|\dot{h}\|_s+\|h\|_s +\|V\|_{s}\big) \, d\tau \tag 10.3
$$

 Let
$h^{(k)}=D_t^k h$.
 Differentiating further in time, using that
 $$
 [D_t,\triangle ] f
= \pa_i\big( \pa^i V^k \, \pa_k  f \big)+ (\pa^i V^k)\pa_k\pa_i f,
\qquad [D_t,\pa_i] =-(\pa_i V^k)\pa_k\tag 16.3
 $$
and using (16.2) for the time derivatives of $V$ it follows that
let $V^{(k)}=D_t^k V$. Then
$$
D_t^k\big( \pa V\big)=(\pa V^{(k_1)})\cdot\cdot\cdot (\pa
V^{(k_n)}) \qquad k_1+...+k_n\leq k-n\tag 16.4
$$
and
$$
D_t^k h^{(k^\prime)}= (\pa V)\cdot\cdot
$$
Actually, it is probably better to do it slightly differently.
 Since
 $$
 \pa_i h=-D_t v_i
 $$
 we obtain
 $$
 \pa_i D_t^k h=-D_t^{k+1} v_i+ (\pa V^{(k_1)} )\cdot\cdot\cdot
 (\pa V^{(k_{n+1})} ) V^{(k_{n+2})},\qquad k_1+...+k_{n+2}
 =k-n,\quad n\geq 0
 $$
Furthermore
 $$
D_t^{k+1} \div V =\div D_t^{k+1} V+(\pa V^{(k_1)} )\cdot\cdot\cdot
 (\pa V^{(k_{n+1})} ) \pa V^{(k_{n+2})},\qquad k_1+...+k_{n+2}
 =k-n,\quad n\geq 0
$$
So, since $D_t^{k+1} \div V\big|_{\pa\Omega}=0$ (recall that $D_t
e(h)+\div V=0$),  it follows that on $\pa\Omega$
$$
\triangle D_t^k h =(\pa V^{(k_1)} )\cdot\cdot\cdot
 (\pa V^{(k_{n+1})} ) \pa V^{(k_{n+2})}+(\pa V^{(k_1)} )\cdot\cdot\cdot
 (\pa^2 V^{(k_{n+1})} ) V^{(k_{n+2})}
 ,\qquad k_1+...+k_{n+2}=k-n,
$$
Then time derivatives of $V$ can be solved in term of time
derivatives of $h$, since $D_t v=-\pa h$. However, it looks
simpler with just $V$. in the right hand side of the above formula
we have at most $k$ time derivatives of $v$ which means at most
$k-1$ time derivatives of $h$.

From all this it seems like one can not reduce all the time
derivatives. One therefore has to control these as well. So we can
have an induction giving high space derivatives in terms of lower
space derivatives and lower time derivatives. So it looks like in
the induction one would have to include all the time derivatives
as well. But perhaps it is okey to do it that way.

We have
$$
\triangle h= -(\pa_i V^j)(\pa_j V^i)+D_t^2 e(h),\tag 16.1
$$
and we also have

$$
D_t^{2+k} e(h)= \triangle D_t^k h + B_k\big(\pa h,\pa^2 h, ...\pa
,...,\pa D_t^{k-2} h,\pa^2 D_t^{k-2} h ,\pa\triangle D_t^{k-2} h,
\pa D_t^{k-1} h, \pa^2 D_t^{k-1} h\big)
$$
\comment
The $L^2$ norm of the highest order terms in the sums in (7.21)
and (7.22) can also be bounded buy
$K_{10}^\prime (\|h\|_{1,\infty}+\|g\|_{1,\infty})E_{s+1}$.  We hence have
$$
\frac{d E_{s+1}^2}{dt}\leq K E_{s+1}^2+\|\hat{D}_t^{s} {f}\|^2
+\sum_{i=0}^{s-1} \|| \check{D}_t^{s+1-i} g\||_{\infty}^2
E_{i+1}^2 + K \sum_{i=0}^{s-1}
  \sum_{j=0}^{s+1-i} \| \, |D_t^{s+1-i-j} \kappa ||D_t^j e^\prime |
|D_t^{i+1} h|\, \|^2\tag 7.24
$$
The terms of highest order $KE_{s+1}^2$ can be controlled by
multiplying by the integrating factor $e^{Kt}$. Similarly, the
term with $i=0$ and $j=s+1$ above can also be controlled by $K
E_{s+1}^2$. We therefore, also assume that $T$ is so small that
$KT\leq 1$. Also we note that $e^\prime $ is bounded by assumption so we
can include the terms with $j=0$ in the first sum.
$$\multline
E_{s+1}(t)^2\leq 4E_{s+1}(0)^2 +K_{10} \int_0^t
\big(\|\hat{D}_t^{s} {f}\|^2+\sum_{i=0}^{s-1}\big( \||
\check{D}_t^{s+1-i} g\||_{\infty}^2 +\|| D_t^{s+1-i} \kappa
\||_{\infty}^2\big) E_{i+1}^2
\big)\, d\tau \\
  +K_{10}
   \sum_{i+j\leq s-1}\int_0^t \| \, |D_t^{s-i-j} \kappa ||D_t^{j+1} e^\prime |
|D_t^{i+1} h|\, \|^2 \, d\tau
\endmultline
$$
For the second row we now use interpolation in space and time to
estimate them by lower order energies. Now, we use induction. For
$s=0$ the sums above disappear and the estimate is of the form
that we want. Using interpolation in space and time we see that
$$
\sum_{i+j\leq s-1}\int_0^T \| \, |D_t^{s-i-j} \kappa ||D_t^{j+1}
e^\prime | |D_t^{i+1} h|\, \|^2 \, d\tau \leq C_T
\sum_{k=0}^{s-1}\sum_{l=1}^{s-k} \||g\||_{\infty, s-l-k}^2 \|| D_t
h\||_{\infty}^{2(k+1)} \|| D_t^{l} h\||_{L^2}^2
$$
We now want to use induction to get estimates for $\int_0^T
E_{s+1}^2\, d\tau $.
\endcomment
\Refs
\ref \no[AG]\by S. Alinhac and P. Gerard\book
    Operateurs pseudo-differentiels et theorem de Nash-Msoer
    \publ Inter Editions and CNRS\yr 1991 \endref
\ref \no [BG] \by M.S. Baouendi and C. Gouaouic \paper
    Remarks on the abstract form of nonlinear Cauchy-Kovalevsky theorems
   \jour Comm. Part. Diff. Eq. \vol 2
   \pages  1151-1162 \yr 1977\endref
\ref \no  [C1] \by D. Christodoulou
      \paper Self-Gravitating Relativistic Fluids:
      A Two-Phase Model\jour Arch. Rational Mech. Anal. \vol 130\yr 1995
      \pages 343-400 \endref
\ref \no  [C2] \by D. Christodoulou
      \paper Oral Communication\yr August 95 \endref
\ref \no [CK]\by  D. Christodoulou and S. Klainerman
      \book The Nonlinear Stability of the Minkowski space-time
      \publ Princeton Univ. Press\yr 1993\endref
\ref \no [CL]\by  D. Christodoulou and H. Lindblad
      \paper On the motion of the free surface of a liquid.
      \jour Comm. Pure Appl. Math. \vol 53\pages 1536-1602\yr 2000\endref
\ref \no [CF] \by R. Courant and K. O. Friedrichs
     \book Supersonic flow and shock waves
     \publ Springer-Verlag \yr 1977\endref
\ref\no [Cr] \by W. Craig \paper An existence theory for water waves and the Boussinesq
      and Korteweg-deVries scaling limits\jour Comm. in P. D. E.
      \vol 10 \yr 1985 \pages 787-1003\endref
\ref\no [DM]\by B. Dacorogna and J. Moser
      \paper On a partial differential equation involving the Jacobian determinant.
       \jour Ann. Inst. H. Poincare Anal. Non. Lineaire\vol 7\yr 1990\pages 1-26\endref
\ref\no [DN] \by S. Dain and G. Nagy \paper
       Initial data for fluid bodies in general relativity
       \jour Phys. Rev. D \vol 65\yr 2002\endref
\ref\no [E1] \by D. Ebin \paper The equations of motion of a perfect fluid
      with free boundary are not well posed. \jour Comm. Part. Diff. Eq. \vol 10
      \pages 1175--1201\yr 1987\endref
\ref\no [E2] \by D. Ebin
      \paper Oral communication \yr November 1997\endref
\ref\no[F]\by H. Friedrich
      \paper Evolution equations for gravitating ideal fluid bodies in general relativity
      \jour  Phys. Rev. D \vol 57\yr 1998\endref
\ref\no[FN]\by H. Friedrich and G. Nagy\paper The initial boundary value problem
      for Einstein's vacuum field equation\jour Commmun. Math. Phys.
      \vol 201\pages 619-655\yr 1999\endref
\ref\no [Ha]\by R. Hamilton \paper Nash-Moser Inverse Function Theorem
      \jour  Bull. Amer. Math. Soc. (N.S.)  7  \yr 1982  \pages  65--222\endref
      \ref\no [H1] \by H\"ormander\paper The boundary problem of Physical geodesy
    \jour Arch. Rational Mech. Anal.  62  \yr 1976\pages  1--52. \endref
\ref\no [H2] \by H\"ormander\paper Implicit function theorems
    \publ Lecture Notes (Stanford )\yr 1977 \endref
\ref\no [H3] \by H\"ormander\paper The analysis of Linear Partial
      Differential Operators III\publ Springer Verlag\endref
\ref\no [K1]\by S. Klainerman\paper On the Nash-Moser-H\"ormander scheme
       \jour unpublished lecture notes\endref
\ref\no [K2]\by S. Klainerman\paper Global solutions of nonlinar wave equations
       \jour Comm. Pure Appl. Math.  33  \yr 1980\pages  43--101.\endref
\ref\no [L1] \by H. Lindblad
      \paper Well posedness for the linearized motion of an incompressible
      liquid with free surface boundary.\jour Comm. Pure Appl. Math.,
       \yr 2003 \endref
\ref\no [L2] \bysame \paper Well posedness for the linearized motion
      of a compressible liquid with free surface boundary.\jour
      Comm. Math. Phys. \yr 2003  \endref
\ref\no [L3] \bysame\paper Well posedness for the motion of an incompressible
     liquid with free surface boundary.\jour preprint:
   http://xxx.lanl.gov/abs/math.AP/0402327,  to appear in the Annals of Math.\endref
\ref \no [Na] \by V.I. Nalimov, \paper The Cauchy-Poisson Problem (in Russian),
      \jour Dynamika Splosh. Sredy 18\yr 1974,\pages 104-210\endref
\ref \no [Ni] \by T. Nishida\paper A note on a theorem of Nirenberg
      \jour J. Diff. Geometry \vol 12 \yr 1977\pages 629-633\endref
\ref\no[R] \by A. D. Rendall \paper
      The initial value problem for a class of general relativistic fluid bodies
     \jour J. Math. Phys. \pages 1047-1053\yr 1992\endref
\ref\no [SY] \by Schoen and Yau\book Lectures on Differential Geometry
        \publ International Press\yr 1994\endref
\ref \no [S] \by E. Stein \paper Singular Integrals and
      differentiability properties of
      functions\publ Princeton University Press\yr 1970\endref
\ref \no [W1]\by S. Wu
      \paper Well-posedness in Sobolev spaces of the full water wave problem in 2-D
      \jour Invent. Math. \vol 130\pages 39-72\yr 1997\endref
\ref \no [W2]\by S. Wu
      \paper Well-posedness in Sobolev spaces of the full water wave problem in 3-D
      \jour J. Amer. Math. Soc. \pages 445-495\vol 12\yr 1999\endref
\ref \no [Y]\by H. Yosihara
       \paper Gravity Waves on the Free Surface of an Incompressible Perfect Fluid
      \publ Publ. RIMS Kyoto Univ. \pages 49-96\vol 18\yr 1982\endref
\endRefs
